\documentclass[pdftex,12pt,a4paper]{amsart}
\usepackage[pdftex]{graphicx}
\DeclareGraphicsExtensions{.png,.pdf}
\usepackage{thumbpdf}
\usepackage{comment}

\graphicspath{{figs/}}
\usepackage[english]{babel}
\usepackage[utf8]{inputenc}
\usepackage{amssymb}
\usepackage{amsmath}
\usepackage{amsthm}
\usepackage{color}
\usepackage{url}
\usepackage{bm}

\newcommand{\R}{{\mathbb R}}
\newcommand{\T}{{\mathbb T}}
\newcommand{\Z}{{\mathbb Z}}
\newcommand{\supertmi}[1]{\{#1\}_{i=0}^{m-1}}
\newcommand{\supertm}[3]{\{#1_{#2}\}_{#2=0}^{#3-1}}
\newcommand{\dparcial}[1]{{\textstyle\frac\partial{\partial#1}}}

\newcommand{\bmi}{\bm{i}}

\def\ee{\mathrm{e}}

\def\deltatau{{\Delta\tau}}
\def\deltalambda{{\Delta\lambda}}
\def\deltah{{\Delta h}}
\def\deltaT{{\Delta T}}

\theoremstyle{remark}
\newtheorem{remark}{Remark}[subsection]
\newtheorem{algo}[remark]{Algorithm}

\begin{document}

\title[Computation of  invariant tori]{
Flow map parameterization methods for invariant tori in Hamiltonian systems
}
\author{Alex Haro}
\address{Departament de Matem\`atiques i Inform\`atica
and BGSMath, Universitat de Barcelona,
Gran Via 585, 08007 Barcelona, Spain.}
\email{alex@maia.ub.es}
\author[J.M. Mondelo]{J.M. Mondelo}
\address{Departament de Matem\`atiques \& CERES-IEEC \& BGSMath,
 Universitat Aut\`onoma de Barcelona. Av. de l'Eix Central, Edifici C, 
 08193 Bellaterra (Barcelona), Spain.}
\email{jmm@mat.uab.cat}
\thanks{
A.H. is supported by the grants PGC2018-100699-B-I00 (MCIU-AEI-FEDER, UE),
2017 SGR 1374 (AGAUR), MSCA 734557 (EU Horizon 2020), and
MDM-2014-0445 (MINECO), and by the NSF under Grant No.~1440140 to found his residence at 
MSRI in Berkeley, California, during the Fall 2018 semester. 
\\
J.M. Mondelo has been supported by the MINECO-AEI grants
MTM2014-52209-C2-1-P, MTM2016-80117-P, MTM2017-86795-C3-1-P.
\\
We would also like to acknowledge the work of the free software community in
providing all the software tools we have used, which are listed at the
beginning of Section~\ref{sec:testex}.
}

\date{}

\maketitle
 
\begin{abstract}
The goal of this paper is to present a methodology for the computation of
invariant tori in Hamiltonian systems combining flow map methods,
parameterization methods, and symplectic geometry.  While flow map methods
reduce the dimension of the tori to be computed by one (avoiding Poincar\'e
maps), parameterization methods reduce the cost of a single step of the
derived Newton-like method to be proportional to the cost of a FFT.
Symplectic properties lead to some magic cancellations that make the methods
work.  The multiple shooting version of the methods are applied to the
computation of invariant tori and their invariant bundles around librational
equilibrium points of the Restricted Three Body Problem. The invariant bundles
are the first order approximations of the corresponding invariant manifolds,
commonly known as the whiskers, which are very important in the dynamical
organization and have important applications in space mission design.

\smallskip

\noindent{\bf Keywords.} Invariant tori; parameterization method; KAM theory;
RTBP; Lissajous orbits.
\end{abstract}

\section{Introduction}

Hamiltonian systems are frequently found in physical and engineering
applications, from where challenging problems continuously emerge, of both
theoretical and practical nature. The developement of efficient methods for
computing invariant tori carrying quasi-periodic motion is a driving force in
the applications of Hamiltonian systems (see e.g.~\cite{Warnock91} for early
references), in areas such as plasma physics, semiclassical quantum theory,
accelerator theory, magnetohydrodynamics, oceanography and, of course,
celestial mechanics. While Lagrangian (maximal dimensional) invariant tori are
important in stability studies, partially hyperbolic invariant tori are also
important in studies of diffusion and chaos in Hamiltonian systems. An
important problem in astrodynamics is the design of station keeping orbits
lying on partially hyperbolic invariant tori around collinear libration points
in RTBP approximations, for which the stable manifolds are sort of entry
lanes, and the unstable manifolds are the exit lanes (see \cite{2003DuFa} for
a survey of early libration points missions, and the web pages of space
agencies for many newer ones). This will be the guiding problem of this paper
to fix a framework. 

To date, one of the most succesful approaches to compute invariant tori falls
in the category of {\em numerical Fourier methods}, in which parameterizations
of tori are given by (truncated) Fourier expansions, and the arising
discretized invariance equations (using e.g. collocation) are solved by
numerical methods such as Newton's method (see e.g.
\cite{Chan83,HuangKM97,CastellaJ00,SchilderVSO06,AdomaitisKL07} for several
variants of this approach in differents contexts). In spite of the relatively
simple formulation of the approach, the main practical drawback is what we
refer to as the {\em large matrix problem} \cite{HaroL06b}: the computational
bottleneck produced in solving the large dimension of the systems of equations
at each iteration step, whose time cost is $O(N^3)$ and the memory cost is
$O(N^2)$, where $N$ is the number of Fourier coefficients of the
approximations. We emphasize that the number of Fourier coefficients has to
do with the dimension of the tori to be computed and their regularity.
Remarkably, the numerical Fourier method introduced in \cite{GomezM01}
mitigates the curse of dimensionality by reducing the dimension of the tori by
1, in looking for invariant tori for time-$T$ flow maps (where $1/T$ is one of
the frequencies of the motion on the torus) instead of looking for invariant
tori for Poincar\'e maps. As of now, it is a well-stablished method that has
proven to be among the most adequate in computing partially hyperbolic
invariant tori around collinear points in the RTBP, by reducing the problem to
computing invariant curves of time-$T$ flow maps (see \cite{2018aBaOlSche} for
a review). As we see, avoiding the large matrix problem is already important
for models such as the RTBP, but is crucial when facing the computation of
higher dimensional tori either in higher dimensional problems (e.g. non
restricted problems), or in the non-autonomous (periodic and quasi-periodic)
improvements of RTBP (such as the elliptic case, the bicircular case, de
quasi-bicircular case, or other models that come from three or more body
problems). As of now, this still has not been attempted in a systematic
manner.

The object of this paper is twofold. First, to overcome the large matrix
problem in \cite{GomezM01} and second, but not less important, to establish a
mathematical framework for the analysis of the derived algorithms. This will
be performed by changing the discretization strategy and linking the approach
to the so called parameterization method for invariant manifolds, a general
strategy for proving existence of invariant manifolds in a constructive way,
so that the methods of proof lead to algorithms of computation, and can be
applied to many different contexts (see \cite{CabreFL03a,CabreFL03b,CabreFL05}
for the foundational papers for invariant manifolds attached to fixed points
and \cite{HaroCFLM16} for a review). More specifically, the algorithms
presented here are inspired by non-perturbative KAM strategies
\cite{Llave01,GonzalezJLV05,FontichLS09,HuguetLS12,LuqueV11} and by symplectic
geometry \cite{GonzalezHL13,HaroCFLM16}, that are applied to look for
invariant tori of flow maps in the spirit of the methodology introduced in
\cite{GomezM01}. The algorithms consist in performing Newton-like steps on
the invariance equations at the functional level (rather than directly at the
numerical level). To do so, the geometrical and dynamical properties of the
problem lead to the construction of a frame especially chosen in order to make
block triangular the linearized invariance equations. Using FFT to switch the
representations of invariant tori from samples to Fourier coefficients and
vice-versa makes the time cost of each step $O(N\log N)$ and the memory cost
$O(N)$, where $N$ is the number of either Fourier coefficients or samples used
to represent the tori. This is a significant improvement of these {\em
functional Fourier methods} with respect {\em numerical Fourier methods}
mentioned above. See \cite{HaroCFLM16} for some benchmarks comparing
parameterization method-like methodologies in the context of invariant tori in
skew-product systems \cite{HaroL06a,HaroL06b,HaroL07}, and also
\cite{CallejaL09,CanadellH17b,Huguet08,JorbaO09} for other contexts.

We will present several algorithms of computation and continuation of
invariant tori, including the isoenergetic case (i.e. invariant tori at a
fixed energy level). For the sake of simplicity, we will focus on partially
hyperbolic invariant tori of dimension $n-1$ of $n$ degrees of freedom
Hamiltonians, but many of the ideas can be extrapolated to other cases
(including Lagrangian tori and partially elliptic invariant tori), even to
non-autonomous Hamiltonian systems (e.g. using the common trick of adding
extra degrees of freedom). The methods we present not only compute the
invariant torus, but also the invariant stable, unstable and center bundles at
the same time. The stable and unstable bundles are the first order
approximations of the corresponding stable and unstable invariant manifolds,
commonly known as the whiskers. The center bundle provides the tangent
directions to the normally hyperbolic invariant cylinder containing the family
of tori being computed. We will come to the a priori unexpected realization
that, when invariant tori and invariant bundles are computed at the same time,
the algorithms are much more efficient than it they are designed for computing
invariant tori only. This is another improvement of the standard approach, in
which the stable and unstable bundles are computed after the torus is computed
\cite{GomezM01,2001Jo}, and avoids additional $O(N^3)$ time and $O(N^2)$
memory costs. The algorithms implement multiple shooting, in order to cope
with the instability that comes from the hyperbolic part. In summary, the
computational bottleneck of the flow map parameterization methods presented in
this paper is no longer the solution of the invariance equations but the
(unavoidable) numerical integration needed in order to evaluate the flow maps.
This is a task that can be performed easily in parallel. 

Last but not least, the convergence of the algorithms presented here could be
proved using KAM methods, under appropriate non-degeneracy conditions that we
make explicit (related to well-known Kolmogorov and isoenergetic conditions)
and smallness of the error of the approximate solutions. We do not pursue
this analysis here. We refer to
\cite{GonzalezJLV05,FontichLS09,LuqueV11,GonzalezHL13,CallejaCL13a,CanadellH17a,HaroCFLM16}
for proofs of several KAM results in a posteriori format, based on the
parameterization method. See also \cite{FiguerasHL17} for a methodology to
perform computer-assisted-proofs based on a posteriori format KAM theorems.
Following the standard practice in numerical analysis, we have tested the
algorithms with well-known computations, as the ones appearing in
\cite{GomezM01}. We will see that, already in this case, the algorithms are
much faster and let one reach unexplored regions and compute tori that are
about to break.

\subsection*{Summary of the paper.}
Section \ref{sec:setting} provides some geometrical background and introduces
notation for the invariance equations to be solved and the parameterizations
to be computed. Section \ref{sec:algo} progressively introduces the necessary
algorithms for the solution of multiple cohomological equations and the
computation of frames, in order to perform Newton steps and finally
continuation ones. The section ends with important comments on an actual
implementation. Section \ref{sec:application} is devoted to the application of
the algorithms to the computation of the family of partially hyperbolic KAM
tori born from the equilibrium point $L_1$ of the Earth-Moon RTBP. A
performance comparison is made with previous large-matrix methodology. Some
dynamical and geometric observables are introduced and graphically
represented, in order to discuss global properties of the family. The
graphical evolution of a few specific sub-families of tori is also shown in
order to illustrate the interaction with other families of objects of the
center manifold of $L_1$. Section~\ref{sec:conclusion} presents some
conclusions. The paper is ended with two appendices.
Appendix~\ref{ap:Poincare map} deals with the equivalence of the Poincar\'e
map and the flow map methods, and Appendix~\ref{ap:quadratiically small
averages} provides the proofs of some cancellations (coming from geometrical
properties) that are crucial for the design of the algorithms, as well as for
eventual proofs of their convergence using KAM methods.

\section{Setting}
\label{sec:setting}

\subsection{Notations}

In this paper we assume that all objects are real analytic.

Let $\T^\ell=\R^\ell/[0,1]^\ell$ be the standard $\ell$-torus. With a slight
abuse of notation, we identify a function $\xi:\T^\ell\to \R$ with a funcion
$\xi:\R^\ell\to\R$ that is $1$-periodic in each of its variables, the
components of  $\theta= (\theta_1,\dots,\theta_\ell)$. The average of $\xi$ is 
\[
	 \langle\xi\rangle:=\int_{\T^\ell}\xi(\theta)d\theta.
\]
We denote the Fourier coefficients of $\xi$ as $\{\hat\xi_k\}_{k\in\Z^\ell}$,
which are given by
\[
   \hat \xi_k = \int_{\T^\ell} \xi(\theta) \ee^{-\bmi 2\pi k \theta} d\theta
   .
\]
Then
\[
	\xi(\theta)= \sum_{k\in\Z^\ell} \hat\xi_k \ee^{\bmi 2\pi k \theta},
\]
where $k\theta:= \sum_{i= 1}^\ell k_i \theta_i$ and $\bmi$ denotes the
imaginary unit. The Fourier coefficients go to zero exponentially fast when
$|k|:=\sum_{i= 1}^\ell |k_i|$ goes to infinity.

\subsection{Symplectic structures and Hamiltonian systems}

We assume we are given an open set $U\subset \R^{2n}$ endowed with an {\em
exact symplectic structure} whose matrix representation is an antisymmetric
matrix map $\Omega:U\to \R^{2n\times 2n}$, which is invertible, and  it is
given by
\[
	\Omega(z) = {\rm D} a(z)^\top - {\rm D} a(z),
\]
where $a:U \to \R^n$, for which the transpose $a(z)^\top$  is the matrix
representation of the {\em action form} at the point $z\in U$.  The matrix map
$\Omega$ induces a {\em symplectic product} at each $z\in U$. For the sake of
simplicity, we will also assume we are given an almost complex structure
compatible with the symplectic structure, meaning a map $J:U\to \R^{2n\times
2n}$ that is involutive ($J(z)^2= -I_{2n}$), symplectic ($J(z)^\top \Omega(z)
J(z)= \Omega(z)$) and  such that the matrix map $G: U\to \R^{2n\times 2n}$
defined by $G(z)= -\Omega(z) J(z)$ induces an {\em scalar product} at each
$z\in U$.

The prototypical example is the standard symplectic structure, given by 
\[
	\Omega_0(z) = \begin{pmatrix} O_n & -I_n \\ I_n & O_n \end{pmatrix}, 
\]
where $I_n$ is the $n\times n$ identity matrix, and $O_n$ is the $n\times n$
zero matrix.  In this case, 
\[
	a_0(z)= \tfrac12 \begin{pmatrix} O_n & I_n \\ -I_n & O_n \end{pmatrix} z,\ 
	J_0(z)=    \begin{pmatrix} O_n & -I_n \\ I_n & O_n \end{pmatrix},\ 
	G_0(z)=  \begin{pmatrix} I_n & O_n \\ O_n & I_n \end{pmatrix}. 
\]

Given a function $H:U\to \R$, we obtain an {\em autonomous Hamiltonian system}
with $n$ degrees of freedom,
\begin{equation}\label{eq:hamsys}
   \dot z=X_H(z):=\Omega(z)^{-1} {\rm D} H(z)^\top.
\end{equation}
We will denote by $\varphi:D \subset \R\times U\to U$ the Hamiltonian flow
associated to the {\em Hamiltonian vector field} $X_H$. We will often use the
notation $\varphi_t= \varphi(t; \cdot)$. For $t$ fixed, the time-$t$ map
$\varphi_t$ is exact symplectic (in the appropriate domains), meaning that is
symplectic,  that is 
\[
	{\rm D}\varphi_t(z)^\top \Omega(\varphi_t(z)) {\rm D}\varphi_t(z)= \Omega(z),
\]
and, moreover, 
\[
	a(\varphi_t(z))^\top {\rm D}\varphi_t(z) - a(z)^\top = {\rm D}p_t(z)
\]
for a certain {\em primitive function} $p_t:U\to \R$, which in fact is given
by
\[
	p_t(z)= \int_0^t \left( a(\varphi_s(z))^\top X_H(\varphi_s(z)) -
        H(\varphi_s(z))  \right) ds.
\]
The exactness property leads to crucial cancellations that enable the
existence of invariant tori.  Moreover, it is well-known that invariant tori
carrying quasi-periodic dynamics have the special geometrical property of
being  isotropic.

\subsection{Isotropic tori and Calabi vectors}

Given a parameterization $K:\T^{m} \to U$ of an $m$-dimensional torus
$\mathcal  K$, we define its Calabi vector $C(K)$ as
\begin{equation}\label{eq:defCalabi}
	C(K)^\top= \int_{\T^{m}} a(K(\theta))^\top {\rm D}K(\theta) \ d\theta,
\end{equation}
where we emphasize the dependence on $K$. Its components are the Calabi
invariants, 
\begin{equation}\label{eq:defCalabicmp}
	C_i(K)= \int_{\T^{m}} a(K(\theta))^\top \tfrac{\partial K}{\partial \theta_i} (\theta) \ d\theta,
\end{equation}
for $i= 1,\dots,m$. The corresponding radii of $K$ are 
$
r_i(K)= \sqrt{\tfrac{|C_i(K)|}{\pi}},
$ 
for $i= 1,\dots,m$,
which measure the widths of the torus $\mathcal K$ (w.r.t. the
parameterization $K$). Notice that Calabi invariants (and the radii) are
invariant under the flow of a Hamiltonian vector field $X_H$, as we prove in
the following lines:
\[
\begin{split}
	C(\varphi_t(K(\theta)))^\top 
	&= \int_{\T^{m}} a(\varphi_t(K(\theta)))^\top {\rm D}(\varphi_t(K(\theta))) \ d\theta 
	\\
	& =   \int_{\T^{m}} a(\varphi_t(K(\theta)))^\top {\rm D}\varphi_t(K(\theta)) {\rm D}K(\theta) \ d\theta 
	\\
	&=   
	\int_{\T^{m}} a(K(\theta))^\top {\rm D}K(\theta) \ d\theta +  \int_{\T^{m}} {\rm D}p_t(K(\theta))  {\rm D}K(\theta)\ d\theta
	\\
	& = C(K)^\top,
\end{split}
\]
where $p_t$ is the primitive function of $\varphi_t$, and we apply that
$p_t\circ K$ is $1$-periodic in all its variables, so its differential has
zero average.

\begin{remark}\label{rm:calabi-chfreqs}
Given a torus automorphism $A:\T^{m} \to \T^{m}$, where $A\in {\rm
GL}_{m}(\Z)$ (i.e.~$A\in \Z^{m\times m},\det A= \pm 1$), the
repa\-ra\-me\-tri\-za\-tion $K\circ A:\T^{m}\to U$ of $\mathcal K$ has Calabi
vector $C(K\circ A)= A^\top C(K)$.
\end{remark}

In the previous constructs we have only used the symplectic properties of
phase space, but not the possible geometrical properties of the tori. The
torus $\mathcal K$ is isotropic if its parameterization satisfies
\[
  {\rm D}K(\theta)^\top \Omega(K(\theta)) {\rm D}K(\theta)= 0,
\]
for any $\theta\in \T^{m}$. In such a case, we may define, for $i= 1,\dots,m$,
\begin{equation}\label{eq:defCalabicmp-new}
	C_i(K)= \int_0^1 a(K(\theta))^\top \tfrac{\partial K}{\partial \theta_i} (\theta) \ d\theta_i
\end{equation}
by taking any fixed $(\theta_1,\dots, \theta_{i-1}, \theta_{i+1}, \dots,
\theta_{m})$ (hence, giving a generator of the torus). It is not difficult to
check that the definition does not depend on such a choice (just compute the
derivatives with respect to $\theta_j$ with $j\neq i$) and, hence, equals the
definition given in \eqref{eq:defCalabicmp}.

\begin{remark}
In the numerical computation of invariant tori it is useful to monitor these
geometrical quantities in order to detect shrinking of the tori. These are
geometrical observables we use along the computations.
\end{remark}

\subsection{Invariance equations for invariant tori}

A parameterization $\hat K:\T^{d}\longrightarrow U$ of an invariant
$d$-dimensional torus $\hat {\mathcal K}$ with {\em frequency vector}
$\hat\omega\in\R^d$ satisfies the invariance equation
\begin{equation}
\label{eq:invtrnm0}
   \varphi_t\bigl(\hat K(\hat\theta)\bigr)
   =
   \hat K(\hat\theta+t\hat\omega),
\end{equation}
for all $\hat \theta \in \T^d$ and $t\in\R$. The infinitesimal version of
\eqref{eq:invtrnm0} is 
\begin{equation}
\label{eq:invtrnm1}
	X_H(\hat K(\hat\theta))= {\rm D} \hat K (\hat \theta)\hat\omega.
\end{equation}
The frequency vector $\hat\omega$ is assumed to be (at least) non-resonant or
ergodic, that is, $\hat k \cdot \hat \omega\neq 0$ for any $\hat k\in
\Z^{d}\setminus\{0\}$.  It is well-known that the Hamiltonian $H$ is constant
on an invariant torus: $H(\hat K(\hat \theta))= h$ for all $\hat \theta\in
\T^d$, for a certain energy $h$.

\begin{remark}
\label{rm:undeterminacy1}
Equation~\eqref{eq:invtrnm1} determines $\hat K$ up to a phase: if $\hat
K(\hat \theta)$ is a solution, then, for any $\hat \alpha\in \R^{d}$, $\hat
K_{\hat\alpha}(\hat \theta):= \hat K(\hat \theta+\hat \alpha)$ is also a
solution, that parameterizes the same torus $\hat {\mathcal K}$.  We then have
$d$ degrees of freedom in the determinacy of the parameterization $\hat K$.
\end{remark}

\begin{remark}
\label{rm:undeterminacy-frequencies}
The frequency vector of $\hat{\mathcal K}$ is defined up to a unimodular
matrix $\hat A\in \Z^{d\times d}$  (with determinant $\pm 1$), so that to the
reparameterization $\hat K \circ \hat A:\T^d\to U$ corresponds the frequency
vector $\hat A^{-1} \hat\omega$.
\end{remark}

In order to reduce the dimension of the parameterization to be computed, and
also to avoid the use of Poincar\'e map, we borrow a trick from
\cite{GomezM01}. By writing $\hat\omega= \tfrac{1}{T}(\omega,1)$, where
$\omega\in\R^{d-1}$ and $T\in \R$, one looks instead for a parameterization
$K:\T^{d-1}\rightarrow U$ of a $(d-1)$-dimensional torus $\mathcal K$ inside
the starting one $\hat {\mathcal K}$ satisfying 
\begin{equation}\label{eq:invtrnm2}
   \varphi_T\bigl(K(\theta)\bigl)=K(\theta+\omega).
\end{equation}
We will refer to $T$ as the period or {\em flying time} of the torus $\mathcal
K$ inside $\hat{\mathcal K}$, and $\omega$ as its {\em rotation vector}.  From
$K$  satisfying \eqref{eq:invtrnm2} we recover the parameterization $\hat K$
satisfying \eqref{eq:invtrnm1}  via the flow through
\begin{equation}
\label{hK}
   \hat K(\hat\theta)
      = \varphi_{\theta_d T}\bigl(
            K(\theta-\theta_d\omega)
         \bigr),
\end{equation}
where $\hat\theta=(\theta,\theta_d)\in \T^{d-1}\times \T$.  From $\hat K$ we
get a $K$ just defining $K(\theta)= \hat K(\theta,0)$.  The rotation vector
$\omega\in \R^{d-1}$ is non-resonant, meaning that $k\cdot \omega\notin\Z$,
for any $k\in\Z^{d-1}\setminus\{0\}$. Notice that $H(K(\theta))= h$ for all
$\theta\in \T^{d-1}$.

\begin{remark}
\label{rm:underterminacy2}
Equation~\eqref{eq:invtrnm2} determines $K$ up to a phase: if $K(\theta)$ is a
solution, then, for any $\alpha\in \R^{d-1}$, $K(\theta+\alpha)$ is also a
solution that parameterizes the same torus $\mathcal K$ inside $\hat {\mathcal
K}$. But Equation~\eqref{eq:invtrnm2} determines $K$ up to a time translation
of $\mathcal K$ inside $\hat {\mathcal K}$: if $K(\theta)$ is a solution then,
for any $\bar \alpha\in \R$, $\varphi_{\bar \alpha} (K(\theta))$ is also a
solution. We then again have $d$ degrees of freedom in the determinacy of the
parameterization $K$.
\end{remark}

\begin{remark}
\label{rm:underterminacy-rotation-vector}
The rotation vector of $\mathcal K$ is defined up to a unimodular matrix $A\in
\Z^{(d-1)\times (d-1)}$, in such a way that the rotation vector corresponding
to the reparameterization $K \circ A:\T^{d-1}\to U$ is $A^{-1}\omega$.
\end{remark}

In the problem of existence (and computation) of invariant tori $\mathcal K$
with rotation vector $\omega$ one may consider two different cases:
\begin{itemize}
\item
   isochronous case: $T$ is fixed, and the torus $K$ and the energy $h$ are
   the unknows;
\item
   isoenergetic case: $h$ is fixed, and the torus $K$ and the flying time $T$
   are the unknowns.
\end{itemize}
In summary, the equations to face are:
\begin{align}
   \varphi_T\bigl(K(\theta)\bigl)-K(\theta+\omega)& = 0, \label{eq:invtrnm}\\
   \int_{\T^{d-1}} H(K(\theta)) \ d\theta - h & = 0, \label{eq:energy}
\end{align}
in which we fix either $T$ or $h$ accordingly. Notice that this formulation
also makes natural to consider either $T$ or $h$ as continuation parameters.

\begin{remark}
\label{rm:Poincare map}
As it is very well-known, using Poincar\'e map is another way of reducing the
dimension of the problem. Both the Poincar\'e map and the time-$T$ flow map
approaches are equivalent (see Appendix~\ref{ap:Poincare map}). From the
theoretical point of view, the time-$T$ approach has the advantage that one
deals with exact symplectic mapings in both the isochronous and isoenergetic
cases, while in the Poincar\'e map approach one produces symplectic mappings
once one reduces also to a fixed energy level (so one works in the
isoenergetic case). From the numerical point of view, both approaches of
course involve numerical integration, but computing Poincar\'e maps and their
differentials is a bit more time consuming task (both in terms of coding and
execution) than computing time-$T$ flows and their differentials. Moreover,
since we are planning to perform multiple shooting methods, the use of
multiple Poincar\'e maps can increase the difficulty of appropriately locating
them and the complexity of the algorithms.
\end{remark}

\subsection{Partially hyperbolic invariant tori of dimension $n-1$}

In this paper, we focus in algorithms for computing partially hyperbolic
invariant tori of dimension $n-1$, i.e. $d= n-1$, that is, with stable and
unstable bundles of rank 1.  The algorithms we present are very easy to adapt
to Lagrangian tori, that is $d= n$, and other lower dimensional partially
hyperbolic tori, that is $d<n$. Invariant tori with elliptic directions could
be also considered, with the aid of  parameters.  The case we consider appears
very often in applications, as the one presented in this paper. In
Section~\ref{sec:testex} we compute partially hyperbolic invariant tori around
the Lagrangian points of the Restricted Three Body Problem, that is $n= 3$,
$d= 2$. 
 
Hence, assume that $d= n-1$. A bundle $\hat{\mathcal W}$ of rank 1 (with base
$\hat{\mathcal K}$) parameterized by a map $\hat W:\T^{n-1}\rightarrow\R^{2n}$
satisfying the equation
\begin{equation}\label{eq:invbdnm1}
   {\rm D}\varphi_t\bigl(\hat K(\hat\theta)\bigr)\hat W(\hat\theta)
   =
   e^{t\chi}\hat W(\hat\theta+t\hat\omega),
\end{equation}
with $\chi\in\R$, is invariant under the linearized flow of $X_H$.  If
$\chi<0$ then $\hat{\mathcal W}= \hat{\mathcal W}^s$ is the stable bundle, and
if $\chi>0$ then $\hat{\mathcal W}= \hat{\mathcal W}^u$ is the unstable
bundle.

In order to decrease the dimension by one, we proceed again by considering
time-$T$ maps. With $\hat\omega= \tfrac{1}{T}(\omega,1)$, we look for a
parameterization $W:\T^{n-2}\rightarrow\R^{2n}$ of a line bundle  ${\mathcal
W}$ of the torus ${\mathcal K}$ such that 
\begin{equation}\label{eq:invbdnm2}
   {\rm D}\varphi_T\bigl(K(\theta)\bigr)W(\theta)
   =
   e^{T \chi}\ W(\theta+\omega).
\end{equation}
Again, we recover $\hat W$ by taking 
\begin{align*}
   \hat W(\hat\theta)
      &= e^{-\theta_d T\chi}
         {\rm D}\varphi_{\theta_d T}\bigl(
            K(\theta-\theta_d\omega)
         \bigr)W(\theta-\theta_d\omega),
\end{align*}
where $\hat\theta=(\theta,\theta_d)$.

\begin{remark}
Notice that with this formulation we assume that the line bundle
$\hat{\mathcal W}$ is trivial (and oriented). Non-oriented line bundles can
fit this formulation through the double-covering trick (see
e.g.~\cite{HaroL07}).
\end{remark}

\section{Multiple shooting algorithms}
\label{sec:algo}

As we have seen in previous section, the equations to solve involve numerical
integration of orbits from points on tori up to a certain $T$.  Since $T$
cannot be chosen to be small, dynamical instability can make the numerical
solution of \eqref{eq:invtrnm}, \eqref{eq:energy} and \eqref{eq:invbdnm2}
difficult. In our experience, with values of $e^{T\chi}$ of the order of
thousands we are not able to observe quadratic convergence of Newton iterates.
Also, continuation steps become very small. The effects of dynamical
instability can be avoided by reducing integration time through the use of
multiple shooting. Multiple shooting is classically introduced in the
numerical solution of boundary value problems for ordinary differential
equations (see e.g.~\cite{stoer}), but it can also be used in the computation
of general invariant objects.

\subsection{Multiple shooting invariance equations}

Instead of looking for a parameterization $K_0$ of a single
$(n-2)$-dimensional torus $\mathcal K_0$ inside the $(n-1)$-dimensional torus
$\hat{\mathcal K}$, we look for parameterizations $\{K_i\}_{i=0}^{m-1}$ of $m$
tori, $K_i:\T^{n-2}\to \R^{2n}$ satisfying the equations
\begin{equation}\label{eq:invtrms}
   \varphi_{T/m}\bigl(K_i(\theta)\bigr)
   -
   K_{i+1}(\theta+\tfrac\omega{m}) = 0
   ,
\end{equation}
for $i=0,\dots,m-1$. In the above equation, and in the following, we assume
the $i$ subindex (here in $K_i$) is defined modulo $m$. In particular,
$K_m=K_0$. Note that, if $\{K_i\}_{i=0}^{m-1}$ is solution of
\eqref{eq:invtrms}, any $K_i$ is solution of \eqref{eq:invtrnm}.  In the
following, we will refer to $\{K_i\}_{i=0}^{m-1}$ as a multiple torus, and
similar notations will be used for other multiple objects (bundles, frames,
etc.)

There are several ways to explicit the energy level of the torus
$\hat{\mathcal K}$, and the one we consider here is that of the average of the
energy on the first torus of  the chain:
\begin{equation}
\label{eq:energyms}
\int_{\T^{n-2}} H(K_0(\theta)) \ d\theta - h  = 0.
\end{equation}

We will also use multiple shooting for computing the invariant bundle
$\hat{\mathcal W}$, and thus look for $\lambda$ and a multiple bundle
$\{W_i\}_{i=0}^{m-1}$, where $W_i:\T^{n-2}\to \R^{2n}$, satisfying
\begin{equation}\label{eq:invbdms}
   {\rm D}\varphi_{T/m}\bigl(K_i(\theta)\bigr)W_i(\theta)
   -
   \lambda W_{i+1}(\theta+\tfrac\omega{m})= 0,
\end{equation}
for $i=0,\dots,m-1$.  Again,  if $\{W_i\}_{i=0}^{m-1}$ is solution of
\eqref{eq:invbdms}, any $W_i$ is solution of \eqref{eq:invbdnm2}, with $e^{T
\chi}= \lambda^m$.

Following the philosophy of the parameterization method, we will look for
adapted frames in which Newton's method for  the invariant equations
\eqref{eq:invtrms} and \eqref{eq:invbdms} can be performed through a sequence
of cohomological equations that are diagonal in Fourier space. In this way,
the large matrix systems that appear when doing direct Fourier discretization
of invariance equations are avoided, and all the computational effort goes in
the numerical integrations necessary to perform the change to the frame and to
obtain the new equations, and in the transformations from Fourier space to
grid space and vice versa, which are done through FFT.

\subsection{Multiple cohomological equations}
\label{sec:coho}

As we have just mentioned, the use of adapted frames is in the basis of our
algorithms, but also is in the core of KAM theory.  This section is devoted to
their application to the analysis of multiple cohomological equations.

In the following, we will fix a non-resonant rotation vector $\omega\in
\R^\ell$.  (In the context of this paper, $\ell= d-1= n-2$.)

\subsubsection{Small divisors equations}

The first equation we consider is the following {\em small divisors
cohomological equation},
\begin{equation}\label{eq:smdiv}
   \xi(\theta)-\xi(\theta+\omega)=\eta(\theta),
\end{equation}
where $\eta:\T^\ell\rightarrow\R$ is given and $\xi:\T^\ell\rightarrow\R$ is
to be found. A necessary condition to solve this equation is that the average
of the right hand side is zero: $\langle\eta\rangle= 0$.  The solution of
\eqref{eq:smdiv} is (formally) straightorward in terms of Fourier
coefficients: if
\begin{equation}\label{eq:expFxieta}
   \xi(\theta)
   =
   \sum_{k\in\Z^\ell}{\hat \xi}_ke^{\bmi 2\pi k\theta}
   ,\quad
   \eta(\theta)
   =
   \sum_{k\in\Z^\ell}{\hat \eta}_ke^{\bmi 2\pi k\theta}
   ,
\end{equation}
then the solutions of \eqref{eq:smdiv} are (formally) given by
\[
   \hat\xi_0\in\R\mbox{ arbitrary}
   ,\quad
   \hat\xi_k=\frac{\hat\eta_k}{1-e^{\bmi 2\pi  k\omega}}\mbox{ for $k\neq0$}.
\]
Note that the divisors $1-e^{\bmi 2\pi  k\omega}$ become arbitrarily small.
For analytic $\eta$, the convergence of the series $\xi$ is ensured by
stronger non-resonance properties of rotation vector $\omega$, such as the
so-called Diophantine condition: from now on, we assume there exist
$\gamma>0$, $\tau\geq d$ such that $|k\omega- q|\geq \gamma |k|_1^{-\tau}$ for
all $k\in\Z^\ell\setminus\{0\}, q\in \Z$. We will denote by ${\mathcal R}
\eta(\theta)$ the only $\xi(\theta)$ solution of 
\begin{equation}\label{eq:smdiv1}
   \xi(\theta)-\xi(\theta+\omega)=\eta(\theta) - \langle\eta\rangle,
\end{equation}
with zero average.

The multiple version of the small divisors cohomological equation
\eqref{eq:smdiv} is: given functions $\{\eta_i\}_{i=0}^{m-1}$,
$\eta_i:\T^\ell\rightarrow\R$, we want to find functions
$\{\xi_i\}_{i=0}^{m-1}$, $\xi_i:\T^\ell\rightarrow\R$, satisfying
\begin{equation}\label{eq:smdivms}
   \xi_i(\theta)-\xi_{i+1}(\theta+\tfrac\omega{m})  =\eta_i(\theta)
   ,\quad
   i=0\dots,m-1.
\end{equation}
As already mentioned, we assume the subindex $i$ to be defined modulo $m$, so
that, in particular $\xi_m=\xi_0$.  A telescopic sum turns \eqref{eq:smdivms}
into a small divisors cohomological equation:
\begin{equation}
   \sum_{i=0}^{m-1}\eta_i(\theta+\tfrac im\omega)
   =
      \xi_0(\theta)-\xi_0(\theta+\omega).
      \label{eq:smdivtmtele}
\end{equation}
Hence, the necessary condition to solve the multiple equation is the average
condition
\[
   \langle \eta \rangle:= \frac{1}{m} \sum_{i=0}^{m-1}\langle\eta_i\rangle = 0.
\]
Here, and in the following we will use the notation $\langle g \rangle$ for
the mean of the averages of a set of functions $\{g_i\}_{i=0}^{m-1}$.  Notice
that any solution of \eqref{eq:smdivtmtele} is given by 
\begin{equation}
\label{eq:xi0}
        \xi_0(\theta)= {\hat \xi_{0,0}} + \sum_{i=0}^{m-1} {\mathcal R}\eta_i
        (\theta+\tfrac{i}m\omega).
\end{equation}

Once the average ${\hat \xi_{0,0}}= \langle\xi_0\rangle$ is fixed,  the
remaining $\xi_i$ can be determined from \eqref{eq:smdivms} by the relations
\[
   \xi_{i}(\theta+\tfrac1m\omega)= \xi_{i-1}(\theta) - \eta_{i-1}(\theta), 
   \quad i= 1\dots,m-1.
\]
Roundoff propagation is reduced by computing every $\xi_i$ independently from
\[
   \sum_{i=0}^{m-1} \eta_{j+i}(\theta+\tfrac im\omega)
   =
      \xi_j(\theta)-\xi_j(\theta+\omega).
\]
Then,
\[
   \xi_j(\theta) 
   = \hat \xi_{j,0}
     + \sum_{i=0}^{m-1} {\mathcal R}\eta_{j+i}(\theta+\tfrac{i}m\omega),
\]
with $\hat\xi_{j,0}= \hat \xi_{{j-1},0} - \hat\eta_{j-1,0}$.  The number of
translations to be done to $\eta_0,\dots,\eta_{m-1}$ (and thus the coding and
run-time overhead of this second approach) is reduced by using
\begin{equation}
\label{eq:xij}
   \xi_j(\theta+\tfrac jm\omega) 
   = \hat \xi_{j,0}
     + \sum_{i=0}^{m-1} {\mathcal R}\eta_{j+i}(\theta+\tfrac{j+i}m\omega),
\end{equation}
instead. Anyway, the overhead of this second approach is negligible against
the cost of numerical integration (an orbit per Fourier coefficient, as we
will see).

\subsubsection{Non-small divisors equations} 
The other kind of equation we will consider is a {\em non-small divisors
cohomological equation},
\begin{equation}\label{eq:nsmdiv}
   \lambda\xi(\theta)-\mu\xi(\theta+\omega)=\eta(\theta),
\end{equation}
where $\lambda,\mu\in\R$, $|\lambda|\neq|\mu|$, $\eta:\T^\ell\rightarrow\R$
are given and $\xi:\T^\ell\rightarrow\R$ is to be found.  Its formal solution
is straightforward in terms of Fourier coefficients: using the notation of
\eqref{eq:expFxieta}, for all $k\in\Z$, 
\[
   \hat\xi_k=\frac{\hat\eta_k}{\lambda-\mu e^{\bmi 2\pi  k\omega}}.
\]
The multiple version of the non-small divisors cohomological equation
\eqref{eq:nsmdiv} is: given $\lambda,\mu\in\R$, $|\lambda|\neq|\mu|$,
$\{\eta_i\}_{i=0}^{m-1}$, $\eta_i:\T^\ell\rightarrow\R$, find
$\{\xi_i\}_{i=0}^{m-1}$ satisfying
\begin{equation}\label{eq:nsmdivms}
   \lambda\xi_i(\theta)-\mu\xi_{i+1}\bigl(\theta+\tfrac\omega m\bigr)
   =
   \eta_i(\theta)
\end{equation}
for $i=0,\dots,m-1$.  As before, and also to reduce roundoff propagation, a
telescopic sum turns \eqref{eq:nsmdivms} into a non-small divisors
cohomological equation for each $\xi_j$:
\[
   \lambda^m\xi_j(\theta)-\mu^m\xi_j(\theta+\omega)
   =
   \sum_{i=0}^{m-1}\mu^i\lambda^{m-1-i}
      \eta_{j+i}(\theta+\tfrac{i}{m}\omega).
\]
These $m$ independent equations (that can be solved in parallel) are all of
the same type as \eqref{eq:nsmdiv}. Also as before, the number of translations
to be done to $\eta_0,\dots,\eta_{m-1}$ is reduced by subtituting $\theta$ by
$\theta+\tfrac jm\omega$ in the previous equation.

\subsection{Computation of adapted frames}
\label{sec:frm}

In the spirit of the parameterization method, for a multiple torus
$\{K_i\}_{i=0}^{m-1}$ and a multiple bundle $\{W_i\}_{i=0}^{m-1}$, we look for
a multiple frame $\{P_i\}_{i= 0}^{m-1}$, $P_i:\T^{n-2}\to \R^{2n\times 2n}$,
such that
\begin{equation}
\label{eq:dflwframe}
   P_{i+1}(\theta+\tfrac\omega{m})^{-1}
      {\rm D}\varphi_{T/m}\bigl(K_i(\theta)\bigr)
      P_i(\theta)
   =
      \left(
         \begin{array}{c|c}
            \Lambda & S_i(\theta) \\
            \hline
                    & \raisebox{0pt}[2.2ex][0ex]{$\Lambda^{-\top}$}
         \end{array}
      \right),
\end{equation}
for $i= 0,\dots, m-1$, with
\begin{equation}\label{eq:defLmb}
   \Lambda= \begin{pmatrix} I_{n-1} & 0 \\ 0 & \lambda \end{pmatrix}
   ,\quad
   S_i(\theta)=
   \begin{pmatrix} S_i^1(\theta) & 0 \\ 0 & 0 \end{pmatrix}
\end{equation}
where $I_{n-1}$ is the $(n-1)\times(n-1)$ identity matrix, $S_i^1(\theta)$ is
a symmetric $(n-1)\times(n-1)$ matrix (referred to as {\em torsion} matrix),
and each $0$ (and empty block) stands for a zero matrix of the corresponding
dimensions. We will see this is possible if $\{K_i\}_{i=0}^{m-1}$ and
$\{W_i\}_{i=0}^{m-1}$ are solutions of equations \eqref{eq:invtrms} and
\eqref{eq:invbdms}, respectively. In our algorithm, the key point is that
\eqref{eq:dflwframe} is approximately true if $\{K_i\}_{i=0}^{m-1}$ and
$\{W_i\}_{i=0}^{m-1}$ are approximate solutions.

First, we define the multiple subframe $\{L_i\}_{i=0}^{m-1}$, $L_i:\T^{n-2}\to
\R^{2n\times n}$, by
\begin{equation}\label{eq:defLi}
   L_i(\theta)=
   \begin{pmatrix}
      {\rm D}K_i(\theta) & X_H(K_i(\theta)) & W_i(\theta)
   \end{pmatrix}
   ,
\end{equation}
By differentiating
\eqref{eq:invtrms} with respect to $\theta$, we obtain
\begin{equation}\label{eq:truco1}
   {\rm D}\varphi_{T/m}\bigl(K_i(\theta)\bigr){\rm D}K_i(\theta)={\rm D}K_{i+1}(\theta+\tfrac\omega{m}).
\end{equation}
By differentiating $\varphi_{T/m}\bigl(\varphi_t(K_i(\theta)\bigr) =
\varphi_t\bigl(K_{i+1}(\theta+\tfrac\omega{m})\bigr)$ with respect to $t$ and
taking $t=0$,
\begin{equation}\label{eq:truco2}
   {\rm D}\varphi_{T/m}\bigl(K_i(\theta)\bigr)X_H\bigl(K_i(\theta)\bigr)
   =
   X_H\bigl(K_{i+1}(\theta+\tfrac\omega{m})\bigr)
   .
\end{equation}
From \eqref{eq:truco1}, \eqref{eq:truco2}, \eqref{eq:invbdms}, and the
definition of $L_i(\theta)$ in equation \eqref{eq:defLi}, we have
\begin{equation}
\label{eq:invariance_lagrangian_bundle}
   {\rm D}\varphi_{T/m}\bigl(K_i(\theta)\bigr)L_i(\theta)
   =
   L_{i+1}(\theta+\tfrac\omega{m}) 
   \Lambda
   .
\end{equation}

It is well-know  that symplecticity properties imply that each subframe
$L_i(\theta)$ is Lagrangian, meaning that \begin{equation}
\label{eq:lagrangian_condition}
	L_i(\theta)^\top \Omega(K_i(\theta)) L_i(\theta)= 0.
\end{equation}

\begin{remark}
We emphasize conditions \eqref{eq:invariance_lagrangian_bundle} and
\eqref{eq:lagrangian_condition} also work for matrices of the form
\[
	L_i'(\theta) = L_i(\theta) \begin{pmatrix} A& 0 \\ 0 & b \end{pmatrix},
\]
where $A\in \R^{(n-1)\times (n-1)}$ and $b\in \R$ are constant and invertible.
In particular, one can scale frames in order to mitigate possible
degeneracies.
\end{remark}

The goal is now completing each Lagrangian subframe $L_i$ to a symplectic
frame $P_i$, by juxtaposing a complementary Lagrangian frame $N_i$. There are
several ways to do so (see \cite{HaroCFLM16}). Here, with the aid of the
compatible almost complex structure $J$, we define 
\begin{equation}\label{eq:defNi}
   \hat N_i(\theta)= J\bigl(K_i(\theta)\bigr)
               L_i(\theta)G_i(\theta)^{-1}
   ,\quad
   G_i(\theta)= L_i(\theta)^\top G(K_i(\theta)) L_i(\theta)
   .
\end{equation}
Now, the matrix 
\begin{equation}
\label{eq:defhatP}
   \hat P_i(\theta)
   =
   \begin{pmatrix}     
      L_i(\theta) & \hat N_i(\theta) 
   \end{pmatrix}   ,
\end{equation}
is symplectic, that is to say 
\[
	\hat P_i(\theta)^\top \Omega(K_i(\theta)) \hat P_i(\theta) = \Omega_0.
\]
Hence, 
\begin{equation}
\label{eq:first_reduction}
   \hat P_{i+1}(\theta+\tfrac\omega m)^{-1} {\rm D}\varphi_{T/m}(K_i(\theta))\hat P_i(\theta)
   =
   \left(
   \begin{array}{c|c}
   \Lambda & \hat S_i(\theta) \\
   \hline
   O_n	          & \raisebox{0pt}[2.2ex][0pt]{$\Lambda^{-\top}$}
   \end{array}
   \right)
\end{equation}
is symplectic (with respect to $\Omega_0$), where $\hat S_i(\theta)$ is a
$n\times n$ matrix given by
\begin{equation}\label{eq:defSi}
   \hat S_i(\theta)
   =
   \hat N_{i+1}(\theta+\tfrac\omega m)^\top
   \Omega(K_{i+1}(\theta+\tfrac\omega m))
   {\rm D}\varphi_{T/m}(K_i(\theta))
   \hat N_i(\theta).
\end{equation}
From symplecticity, it follows that $\hat S^i(\theta) \Lambda^\top= \Lambda
\hat S^i(\theta)^\top$,

In order to have \eqref{eq:dflwframe}, we perform a new change of frame by
considering, for $i=0,\dots,m-1$, matrices
\begin{equation}\label{eq:defQ}
   Q_i(\theta)=
   \left(
      \begin{array}{c|c}
         I_n & B_i(\theta) \\
         \hline
             & I_n
      \end{array}
   \right)
      ,
 \end{equation}
with $B_i(\theta)^\top=B_i(\theta)$, so that they are symplectic. Then,
\[
   Q_i(\theta)^{-1}
   =
   \left(
      \begin{array}{c|c}
         I_{n}     & -B_i(\theta) \\
      \hline
           & I_{n}
      \end{array}
   \right)
   .
\]
For the new frame
\begin{equation}\label{eq:defPi}
   P_i(\theta)=\hat P_i(\theta)Q_i(\theta),
\end{equation}
we have
\[
      P_{i+1}(\theta+\tfrac\omega{m})^{-1}
         {\rm D}\varphi_{T/m}(K_i(\theta))
         P_i(\theta)
   =
      \left(
         \begin{array}{c|c}
            \Lambda & S_i(\theta) \\
         \hline
                            & \raisebox{0pt}[2.2ex][0pt]{$\Lambda^{-\top}$}
         \end{array}
      \right)
   ,
\]
with
\begin{equation}\label{eq:Si}
   S_i(\theta)=
               \Lambda B_i(\theta)
               +\hat S_i(\theta)
               -B_{i+1}(\theta+\tfrac\omega{m})\Lambda^{-\top}
   .
\end{equation}
By splitting the matrix $\hat S_i(\theta)$, in blocks of sizes
$(n-1)\times(n-1)$, $(n-1)\times1$, $1\times(n-2)$ and $1\times1$, as in
\begin{equation}\label{eq:splittinghatSi}
\hat S_i(\theta)
   =
   \begin{pmatrix}
      \hat S_i^1(\theta) & \hat S_i^2(\theta) \\
      \hat S_i^3(\theta) & \hat S_i^4(\theta)
   \end{pmatrix},
\end{equation}
and using analogous splittings for $S_i(\theta)$ and $B_i(\theta)$,
formula~\eqref{eq:Si} reads
\begin{align}
   S_i^1(\theta) &= \hat S_i^1(\theta) + \phantom{\lambda} B_i^1(\theta)- B_{i+1}^1(\theta+\tfrac\omega{m}) ,\label{eq:Si1}\\
   S_i^2(\theta) &= \hat S_i^2(\theta)
                     +\phantom{\lambda} B_i^2(\theta)-B_{i+1}^2(\theta+\tfrac\omega{m}) \lambda^{-1},
                     \label{eq:Si2}\\
   S_i^3(\theta) &= \hat S_i^3(\theta)
                     +\lambda B_i^3(\theta)-B_{i+1}^3(\theta+\tfrac\omega{m}),
                     \label{eq:Si3}\\                
   S_i^4(\theta) &= \hat S_i^4(\theta)
                     +\lambda B_i^4(\theta)-B_{i+1}^4(\theta+\tfrac\omega{m}) \lambda^{-1}.
                     \label{eq:Si4}
\end{align}
We then take $B_i^1(\theta)= I_{n-1}$ for each $i= 0,\dots, m-1$, so that
$S_i^1(\theta)= \hat S_i^1(\theta)$, and look for $B_i^3(\theta)^\top=
B_i^2(\theta)$, $B_i^4(\theta)$ so that $S_i^3(\theta)^\top= S_i^2(\theta)
\lambda= 0$ and $S_i^4(\theta)= 0$.  In summary, we need to find
$\{B_i^2(\theta)\}_{i=0}^{m-1}$ and $\{B_i^4(\theta)\}_{i=0}^{m-1}$ such that
\begin{align}
   B_i^2(\theta)-B_{i+1}^2(\theta+\tfrac\omega{m}) \lambda^{-1}
      &= -\hat S_i^2(\theta),
      & i&=0,\dots,m-1,
      \label{eq:bi}
      \\
   \lambda B_i^4(\theta)-B_{i+1}^4(\theta+\tfrac\omega{m})\lambda^{-1}
      &= -\hat S_i^4(\theta),
      & i&=0,\dots,m-1.
      \label{eq:Di}
\end{align}
The solution of these multiple shooting cohomological equations has been
discussed in Section~\ref{sec:coho}.

All the previous developments are summarized in the algorithm that follows.

\begin{algo}\label{algo:frame}
Given $\{K_i\}_{i=0}^{m-1}$, $\{W_i\}_{i=0}^{m-1}$ satisfying
\eqref{eq:invtrms}, \eqref{eq:invbdms}, compute the multiple frame
$\{P_i\}_{i=0}^{m-1}$ and the corresponding reduced expression of $\{{\rm
D}\varphi_{T/m}\circ K_i\}_{i=0}^{m-1}$, as given by \eqref{eq:dflwframe}, by
following these steps:
\begin{enumerate}
\item
   Compute $\{L_i\}_{i=0}^{m-1}$ from \eqref{eq:defLi}.
\item Compute 
   $\{\hat N_i\}_{i=0}^{m-1}$ from \eqref{eq:defNi}, in order to get 
   $\{\hat P_i\}_{i=0}^{m-1}$ from \eqref{eq:defhatP}.
\item\label{en:frnumi}
   Compute $\{\hat S_i\}_{i=0}^{m-1}$ from \eqref{eq:defSi}, and 
   obtain $\supertmi{S_i^1}$ as $\supertmi{\hat S_i^1}$, according to the
   matrix splitting in \eqref{eq:splittinghatSi}.
\item
   Compute $\supertmi{B_i^2}$ from  \eqref{eq:bi}, and $\supertmi{B_i^3}$ from  $B_i^3= (B_i^2)^\top$.
\item 
   Compute $\supertmi{B_i^4}$ from \eqref{eq:Di}. 
\item 
   Define $\supertmi{B_i^1}$ as $B_i^1=I_{n-1}$ for $i= 0\dots m-1$.
\item
   Compute $\supertmi{P_i}$ from \eqref{eq:defPi}, with $\{ Q_i\}_{i=
   0}^{m-1}$ given by \eqref{eq:defQ}.
\end{enumerate}
\end{algo}

\begin{remark}
Observe that step \ref{en:frnumi} is the only one that requires numerical
integration. The other operations are diagonal either in Fourier space or in
grid space.
\end{remark}

\begin{remark}
\label{remark:constant torsion}
We can extend the arguments and reduce the torsion matrices $\{S_i^1\}_{i=
0}^{m-1}$ to constant coefficients by considering multiple small divisors
cohomological equations \eqref{eq:Si1}. To do so, we define 
\begin{align}
\label{eq:def hS1}
 	S^1_0 & = \langle S^1 \rangle= \frac1m \sum_{i= 0}^{m-1} \langle \hat S_i^1(\theta)\rangle 
\end{align}
and solve
\begin{align}
\label{eq:newbi1}
       B_i^1(\theta)-B_{i+1}^1(\theta+\tfrac\omega{m})
      &= S^1_0 -\hat S_i^1(\theta)& i&=0,\dots,m-1.
\end{align}
as we discussed in Section~\ref{sec:coho}. With the choice, the torsion
matrices are  $S^1_i(\theta)= S^1_0$, for all $i =0,\dots,m-1$.  Hence, step
(6) of Algorithm~\ref{algo:frame} can be replaced by 
\begin{itemize}
\item[(6')]
   Compute $\supertmi{B_i^1}$ from \eqref{eq:newbi1}, with $S_i^1=  S^1_0$
   given from \eqref{eq:def hS1}.
\end{itemize}
\end{remark}

\subsection{Description of a Newton step}
\label{sec:nwt}

The goal of this subsection is to develop the formulation necessary to perform
Newton steps in the (multiple shooting) invariance equations for the tori,
\eqref{eq:invtrms}, the energy level, \eqref{eq:energyms}, and the bundles,
\eqref{eq:invbdms}. This will be done by solving the linearization of these
equations around a known approximation expressed in the frame
\eqref{eq:defPi}. We will consider isochronous and isoenergetic cases. 

\subsubsection{A Newton step on the torus}
\label{sec:nwttr}

Let us consider  a multiple torus $\{K_i\}_{i=0}^{m-1}$, a multiple bundle
$\{W_i\}_{i=0}^{m-1}$ and $\lambda\neq 0$ satisfying equations
\eqref{eq:invtrms}, \eqref{eq:energyms} and \eqref{eq:invbdms} approximately,
for a given $T$ (and fixed $\omega$). Let $\{E_i\}_{i=0}^{m-1}$ be the error
in the invariant equations, so , $E_i:\T^{n-2}\to \R^{2n}$ is defined by   
\begin{equation}\label{eq:errtr}
   E_i(\theta):=
      \varphi_{T/m}\bigl(K_i(\theta)\bigr)
      -K_{i+1}(\theta+\tfrac\omega{m})
   ,
\end{equation}
for $i=0,\dots,m-1$.  We also consider the energy error for a given $h$:
\[
 E_h= \langle H(K_0) \rangle - h.
\]
We plan  to give rather explicit formulas for the corrections of tori,
$\{\Delta K_i\}_{i=0}^{m-1}$, and also of the energy, $\deltah$, and the
flying time $\deltaT$, in the isochronous case (for which $\deltaT= 0$) and
the isoenergetic case (for which $\deltah= 0$).

With $\{P_i\}_{i=0}^{m-1}$ the frame defined in \eqref{eq:defPi}, we write the
correction of the tori in the form $\Delta K_i(\theta)=
P_i(\theta)\xi_i(\theta)$.  Expanding by Taylor up to first order the
invariance equation
\begin{align*}
   \varphi_{(T+\deltaT)/m}\bigl(
      K_i(\theta)+P_i(\theta)\xi_i(\theta)
   \bigr)
   \hspace*{9em} \\
   -K_{i+1}(\theta+\tfrac\omega{m})-P_{i+1}(\theta+\tfrac\omega{m})
                               \xi_{i+1}(\theta+\tfrac\omega{m})
   &= 0
\end{align*}
around the approximated tori and flying time, and neglectic second order error
terms, we get the equation for one step of Newton's method:
\begin{align*}
   {\rm D}\varphi_{T/m}\bigl(K_i(\theta)\bigr)P_i(\theta)\xi_i(\theta)
   \hspace*{13em} \\
   -P_{i+1}\bigl(\theta+\tfrac\omega m\bigr)
           \xi_{i+1}\bigl(\theta+\tfrac\omega m\bigr)
           + X_H\bigl(K_{i+1}\bigl(\theta+\tfrac\omega m\bigr)\bigr)\deltatau
   &= -E_i(\theta),
\end{align*}
where $\deltatau= \tfrac{\deltaT}{m}$.

Multiplying the previous equations by $P_{i+1}(\theta+\tfrac\omega m)^{-1}$,
if the frame $\{P_i\}_{i=0}^{m-1}$ satisfied \eqref{eq:dflwframe} exactly, we
would obtain
\begin{equation}\label{eq:corrtrmtr}
   \begin{pmatrix} \Lambda & S_i(\theta) \\ 0 & \Lambda^{-\top} \end{pmatrix}
   \xi_i(\theta)
   -\xi_{i+1}(\theta+\tfrac\omega{m})
   +e\Delta\tau
   =
   \eta_i(\theta)
\end{equation}
with $\Lambda,S_i(\theta)$ defined as in \eqref{eq:defLmb},
\begin{equation}\label{eq:defetatr}
   \eta_i(\theta)=-P_{i+1}(\theta+\textstyle\tfrac\omega m)^{-1}E_i(\theta),
\end{equation}
and
\[
   \xi_i(\theta)=
      \begin{pmatrix}
         \xi_i^1(\theta)\\ \xi_i^2(\theta)\\ \xi_i^3(\theta)\\ \xi_i^4(\theta)
      \end{pmatrix}
   ,\quad
   e = \begin{pmatrix} e_{n-1}\\ 0 \\ 0 \\ 0 \end{pmatrix}
   ,\quad
   \eta_i(\theta)=
      \begin{pmatrix}
         \eta_i^1(\theta)\\ \eta_i^2(\theta)\\ \eta_i^3(\theta)\\
         \eta_i^4(\theta)
      \end{pmatrix}
   ,
\]
being $e_{n-1}= (0,\dots, 1)^\top\in \R^{n-1}$. Notice we implicitly consider
$\xi_i, \eta_i:\T^{n-2} \to \R^{n-1}\times \R \times \R^{n-1} \times \R$ and
enumerate the corresponding block components accordingly.  Actually, since
$\{K_i\}_{i=0}^{m-1}$, $\{W_i\}_{i=0}^{m-1}$ satisfy equations
\eqref{eq:invtrms}, \eqref{eq:invbdms} approximately, the frame also satisfies
\eqref{eq:dflwframe} approximately, so we would need to add an error term to
the matrix in \eqref{eq:corrtrmtr}. This error term can be disregarded,
because when multiplied by $\xi_i(\theta)$ becomes of second order. Taking
this into accout, we can rewrite \eqref{eq:corrtrmtr} as a system of equations
in order to obtain
\begin{align}
   \xi_i^1(\theta)+S_i^1(\theta)\xi_i^3(\theta)
         -\xi_{i+1}^1(\theta+\tfrac\omega m) + e_{n-1} \deltatau
      &= \eta_i^1(\theta),\label{eq:corrtr1} \\
   \lambda\xi_i^2(\theta)-\xi_{i+1}^2(\theta+\tfrac\omega m) 
   \phantom{+  e_{n-1} \deltatau}
      &= \eta_i^2(\theta),\label{eq:corrtr2} \\
   \xi_i^3(\theta)-\xi_{i+1}^3(\theta+\tfrac\omega m)
    \phantom{+  e_{n-1} \deltatau}
      &= \eta_i^3(\theta),\label{eq:corrtr3} \\
   \lambda^{-1}\xi_i^4(\theta)-\xi_{i+1}^4(\theta+\tfrac\omega m)
    \phantom{+  e_{n-1} \deltatau }
      &= \eta_i^4(\theta),\label{eq:corrtr4}
\end{align}
for $i=0,\dots,m-1$.

Equations \eqref{eq:corrtr2}, \eqref{eq:corrtr4} can be solved as multiple
non-small divisors cohomological equations of the type \eqref{eq:nsmdivms}.
Equation \eqref{eq:corrtr3} is a multiple small divisors cohomological
equation of the type \eqref{eq:smdivms}. Solving it as such implies assuming
that
\[
   \langle \eta^3 \rangle= 0.
\]
In Appendix~\ref{ap:quadratiically small averages} we prove that exact
symplectic properties of the flow imply that $\langle \eta^3 \rangle$ is in
fact quadratically small, so the previous assumption is coherent with a Newton
step.  We denote  $\{\bar\xi_i^3\}_{i=0}^{m-1}$ to be the solution of
\eqref{eq:corrtr3} with $\langle\bar\xi_0^3\rangle=0$ and define
\[
   \xi_i^3(\theta)=\xi_{0,0}^{3}+\bar\xi_i^3(\theta),
\]
for $i=0,\dots,m-1$ with $\xi_{0,0}^{3}$ free. Notice that
$\{\xi_i^3\}_{i=0}^{m-1}$ is a solution of \eqref{eq:corrtr3} for any
$\xi_{0,0}^{3}$, which will be fixed later. Substituting
$\{\xi_i^3\}_{i=0}^{m-1}$ in \eqref{eq:corrtr1}, we obtain
\begin{equation}\label{eq:corrtr1xi30}
   \xi_i^1(\theta)-\xi_{i+1}^1(\theta+\tfrac\omega m)
   =
   \eta_i^1(\theta)-S_i^1(\theta)\bar\xi_i^3(\theta)-S_i^1(\theta)\xi_{0,0}^3 -  e_{n-1} \deltatau
   .
\end{equation}
For this last equation to be solved as a multiple shooting small divisors
cohomological equation \eqref{eq:smdivms}, we need the sum of averages of the
right hand sides to be zero, which gives
\[
  \sum_{i= 0}^{m-1} \langle S_i^1\rangle\xi_{0,0}^3 + m e_{n-1} \deltatau
   =
   \sum_{i= 0}^{m-1} \langle\eta_i^1-S_i^1\bar\xi_i^3\rangle
   ,
\]
which is a linear system for $\xi_{0,0}^3$ and $\deltatau$ equivalent to 
\[
\langle S^1 \rangle \xi_{0,0}^3 + e_{n-1} \deltatau= 
\langle \eta^1 - S^1  \bar \xi^3 \rangle.
\]

Now, considering the error in energy as
$E_h=\int_{\T^{n-2}}H(K_0(\theta))d\theta-h$, by substitution of the corrected
first torus and energy in \eqref{eq:energyms} and linearization, it follows
that
\[
	\langle {\rm D}H(K_0(\theta)) P_0(\theta) \xi_0(\theta) - \deltah
        \rangle= - E_h.
\]
Moreover, by using that the frame $P_0(\theta)$ is approximately symplectic,
\[
	{\rm D}H(K_0(\theta)) P_0(\theta) = 
	-X_H(K_0(\theta))^\top \Omega(K_0(\theta)) P_0(\theta) 
	\simeq \begin{pmatrix} 0_{n-1}^\top & 0 & e_{n-1}^\top & 0 \end{pmatrix}.
\]
Hence, neglecting second order error terms we get 
\begin{equation}\label{eq:eninfrm}
	e_{n-1}^\top \xi_{0,0}^3 - \deltah = -E_h.
\end{equation}
By collecting the linear equations for $\xi_{0,0}^3$, $\deltatau$ and
$\deltah$ we get the $(n-1)\times n$ system
\[
\begin{pmatrix} 
\langle S^1\rangle & e_{n-1} & 0 \\
e_{n-1}^\top                       &  0         & -1
\end{pmatrix}
\begin{pmatrix}
\xi_{0,0}^3 \\ \deltatau \\ \deltah
\end{pmatrix}
= 
\begin{pmatrix}
\langle \eta^1 - S^1  \bar\xi^3 \rangle
\\
-E_h
\end{pmatrix}
\]
It suffices system matrix has rank $n-1$ to get one-parameter families of
solutions.  We consider two cases:
\begin{itemize}
\item
   Isochronous case: $\deltatau= 0$, and $\deltah= E_h + e_{n-1}^\top
   \xi_{0,0}^3$, where 
   $\xi_{0,0}^3$ solves linear equation
   \begin{equation}\label{eq:sislinisoc}
            \langle S^1\rangle  \xi_{0,0}^3  = 
            \langle \eta^1 - S^1  \bar\xi^3\rangle
   \end{equation}
   provided that 
   \[
      \det \langle S^1\rangle \neq 0 .
   \]
   This is the isochronous twist condition, which corresponds to Kolmogorov
   condition.  (In practice, the energy equation is not considered, so
   $\deltah$ is not computed.)
\item
   Isoenergetic case: $\deltah= 0$, and $\xi_{0,0}^3, \deltatau$ solve the
   linear equation
   \begin{equation}\label{eq:sislinisoe}
      \begin{pmatrix} 
      \langle S^1\rangle & e_{n-1} \\
      e_{n-1}^\top                       &  0         
      \end{pmatrix}
      \begin{pmatrix}
      \xi_{0,0}^3 \\ \deltatau 
      \end{pmatrix}
      = 
      \begin{pmatrix}
      \langle \eta^1 - S^1  \bar\xi^3 \rangle
      \\
      -E_h
      \end{pmatrix},
   \end{equation}
   provided that 
   \[
      \det \begin{pmatrix} 
      \langle S^1\rangle & e_{n-1} \\
      e_{n-1}^\top                       &  0         
      \end{pmatrix}
       \neq 0 .
   \]
   This is the isoenergetic twist condition. 
\end{itemize}
 
Once one computes $\xi_{0,0}^3$, $\{\xi_3^i\}_{i=0}^{m-1}$ is fully determined
and one can solve \eqref{eq:corrtr1xi30}. The general solution is 
\begin{equation}\label{eq:xi1fromzravg}
   \xi_i^1(\theta)=\xi_{0,0}^{1}+\bar\xi_i^1(\theta),
\end{equation}
for $i=0,\dots,m-1$, where $\xi_{0,0}^{1}$ is free and
$\{\bar\xi_i^1\}_{i=0}^{m-1}$ is the solution of \eqref{eq:corrtr1xi30} with
$\langle\bar\xi_0^1\rangle=0$.  The freedom to choose $\xi^1_{0,0}$ has to do
with the phase and time underterminacy of the parameterization of the first
(and then all) tori (see Remark~\ref{rm:underterminacy2}).  A simple choice is
$\xi_{0,0}^{1}= 0$.

All the process described is summarized in the algorithm that follows.

\begin{algo}\label{algo:nwttr}
(Newton step on a multiple torus, isochronous or isoenergetic case)
Let $\{K_i\}_{i=0}^{m-1}$, $\{W_i\}_{i=0}^{m-1}$, $\lambda$ satisfy equations
\eqref{eq:invtrms},
\eqref{eq:energyms}, \eqref{eq:invbdms} approximately. Obtain the corrected
tori and, in the isoenergetic case, also the corrected flying time 
by following these steps:
\begin{enumerate}
\item
   Compute $\supertm P i m$, $\supertm{S^1}i m$, following
   Algorithm \ref{algo:frame}.
\item
   Compute the multiple error $\supertm E i m$ from \eqref{eq:errtr}.
\item
   Compute the right-hand side of the cohomological equations $\supertm\eta i
   m$ from \eqref{eq:defetatr}.
\item
   Solve \eqref{eq:corrtr2}, \eqref{eq:corrtr4} as non-small divisors
   multiple cohomological equations, in order to obtain $\supertm{\xi^2} i m$,
   $\supertm{\xi^4} i m$.
\item
   Solve \eqref{eq:corrtr3} as small divisors multiple cohomological equation,
   in order to obtain its zero-average solution $\supertm{\bar\xi^3} i m$.
\item
   In the isochronous (resp.~isoenergetic) case, compute $\langle S^1\rangle$
   and the right-hand side of the linear system \eqref{eq:sislinisoc}
   (resp.~\eqref{eq:sislinisoe}) and solve it in order to obtain $\xi_{0,0}^3$
   (resp.~$\xi_{0,0}^3,\Delta\tau$).
\item
   Solve \eqref{eq:corrtr1} as small divisors multiple cohomological equation
   in order to obtain $\supertm{\bar\xi^1} i m$, and obtain
   $\supertm{\xi^1} i m$ from \eqref{eq:xi1fromzravg} by choosing
   $\xi_{0,0}^1=0$.
\item
   Compute the corrected tori as
   $\{K_i(\theta)+P_i(\theta)\xi_i(\theta)\}_{i=0}^{m-1}$. In the isoenergetic
   case, obtain also the corrected flying time as $T+m\Delta\tau$.
\end{enumerate}
\end{algo}

\subsubsection{A Newton step on the bundle}
\label{sec:nwtbd}

Consider again parameterizations $\{K_i\}_{i=0}^{m-1}$ of $(n-2)$-dimensional
tori inside a larger $(n-1)$-dimensional torus, and associated
parameterizations $\{W_i\}_{i=0}^{m-1}$ of $(n-2)$-dimensional bundles inside
the bundle of the larger torus, that satisfy equations \eqref{eq:invtrms} and
\eqref{eq:invbdms} approximately.  (In the implementation, the $K_i$'s are the
ones we have improved in the previous step.) Denote the error in the invariant
equations of the $W_i$ as
\begin{equation}\label{eq:errbd}
   E_i^W(\theta)
   :=
   {\rm D}\varphi_{T/m}\bigl(K_i(\theta)\bigr)W_i(\theta)
   -\lambda W_{i+1}(\theta+\tfrac\omega m),
\end{equation}
for $i=0,\dots,m-1$.  Recalling the frame $\{P_i\}_{i=0}^{m-1}$ defined in
\eqref{eq:defPi}, we would like to find corrections $P_i(\theta)\xi_i(\theta)$
of the bundles $W_i(\theta)$ and a correction $\deltalambda$ for the
eigenvalue $\lambda$ such that the corrected bundles
$W_i(\theta)+P_i(\theta)\xi_i(\theta)$ and corrected eigenvalue
$\lambda+\delta$ make \eqref{eq:invbdms} vanish at first order. The invariance
equation on the corrected bundles and eigenvalue is
\begin{align*}
   \lefteqn{
      {\rm D}\varphi_{T/m}\bigl(K_i(\theta)\bigr)
            \bigr(W^i(\theta)+P_i(\theta)\xi_i(\theta)\bigr)
   }\hspace*{4em}
   \\
   & -(\lambda+\deltalambda)\bigl(
         W_{i+1}(\theta+\tfrac\omega m)
         +P_{i+1}(\theta+\tfrac\omega m)\xi_{i+1}(\theta+\tfrac\omega m)
   \bigr)
   = 0.
\end{align*}
Expanding the parentheses and neglecting errors of second order, as the ones
with a factor
$\deltalambda\ \xi_{i+1}(\theta+\tfrac\omega m)$, the previous equation becomes
\begin{align*}
   \lefteqn{
       {\rm D}\varphi_{T/m}\bigl(K_i(\theta)\bigr) P_i(\theta) \xi_{i}(\theta)
   }\hspace*{4em}
   \\
   & -\lambda P_{i+1}(\theta+\tfrac\omega m)\xi_{i+1}(\theta+\tfrac\omega m)
     -\deltalambda\ W_{i+1}(\theta+\tfrac\omega m) = -E_i^W(\theta).
\end{align*}
Multiplying by $P_{i+1}(\theta+\tfrac\omega m)^{-1}$, using
\eqref{eq:dflwframe} and neglecting second order error terms, we obtain
\begin{equation}\label{eq:corrbdmtr}
   \begin{pmatrix} \Lambda & S_i(\theta) \\ 0 & \Lambda^{-\top} \end{pmatrix}
   \xi_i(\theta)
   -\lambda\xi_{i+1}(\theta+\tfrac\omega m)
   -e\Delta\lambda
   =
   \eta_i(\theta),
\end{equation}
with $\Lambda,S_i(\theta)$ defined as in \eqref{eq:defLmb}, 
\begin{equation}\label{eq:defetabd}
   \eta_i(\theta)=-P_{i+1}(\theta+\tfrac\omega m)^{-1}E_i^W(\theta),
\end{equation}
and
\[
   \xi_i(\theta)=
      \begin{pmatrix}
         \xi_i^1(\theta)\\ \xi_i^2(\theta)\\ \xi_i^3(\theta)\\ \xi_i^4(\theta)
      \end{pmatrix}
   ,\quad
   e = \begin{pmatrix} 0\\ 1 \\ 0 \\ 0 \end{pmatrix}
   ,\quad
   \eta_i(\theta)=
      \begin{pmatrix}
         \eta_i^1(\theta)\\ \eta_i^2(\theta)\\ \eta_i^3(\theta)\\
         \eta_i^4(\theta)
      \end{pmatrix}
   .
\]
As before, we implicitly consider $\xi_i, \eta_i:\T^{n-2} \to \R^{n-1}\times
\R \times \R^{n-1} \times \R$ and enumerate the corresponding block components
accordingly.  Rewriting \eqref{eq:corrbdmtr} as a system of equations, we
obtain
\begin{align}
   \xi_i^1(\theta)+S_i^1(\theta)\xi_i^3(\theta)
         -\lambda\xi_{i+1}^1(\theta+\tfrac\omega m)
         \phantom{\ \,-\deltalambda}
      &= \eta_i^1(\theta),\label{eq:corrbd1} \\
   \lambda\xi_i^2(\theta)-\lambda\xi_{i+1}^2(\theta+\tfrac\omega m)
   -\deltalambda
      &= \eta_i^2(\theta),\label{eq:corrbd2} \\
   \xi_i^3(\theta)-\lambda\xi_{i+1}^3(\theta+\tfrac\omega m)
      \phantom{\ \,-\deltalambda}
      &= \eta_i^3(\theta),\label{eq:corrbd3} \\
   \lambda^{-1}\xi_i^4(\theta)-\lambda\xi_{i+1}^4(\theta+\tfrac\omega m)
      \phantom{\ \,-\deltalambda}
      &= \eta_i^4(\theta),\label{eq:corrbd4}
\end{align}
for $i=0,\dots,m-1$. Equations \eqref{eq:corrbd3} and \eqref{eq:corrbd4} can
be solved as multiple non-small divisors cohomological equations of the form
\eqref{eq:nsmdivms}.  Once $\xi_i^3(\theta)$ is known, \eqref{eq:corrbd1} is
also solved as multiple non-small divisors cohomological equation.  For
\eqref{eq:corrbd2} to be solved as multiple small divisors cohomological
equation of the form \eqref{eq:smdivms}, we need that
\[
   \sum_{i=0}^{m-1}\langle\eta_i^2(\theta)+\deltalambda\rangle=0,
\]
which is achieved by taking $\deltalambda=-\langle\eta^2\rangle$.  In doing
so, $\langle \xi^2_0 \rangle$ remains free.  In this case, the freedom to
choose this average is related to the underdeterminacy of selecting the
lengths of  $\{W_i\}_{i= 0}^{m-1}$ (as, say, the undeterminacy of selecting
the lenght of an eigenvector of a given eigenvalue of a matrix). A simple
choice is to take $\langle \xi^2_0 \rangle= 0$.

The previous steps are summarized in algorithm that follows.

\begin{algo}\label{algo:nwtbd}
(Newton step on a multiple bundle) Let $\{K_i\}_{i=0}^{m-1}$,
$\{W_i\}_{i=0}^{m-1}$, $\lambda$ satisfy equations \eqref{eq:invtrms},
\eqref{eq:energyms}, \eqref{eq:invbdms} approximately. Obtain the corrected
bundle and eigenvalue by following these steps:
\begin{enumerate}
\item
   Compute $\supertm P i m$, $\supertm{S^1}i m$, following
   Algorithm \ref{algo:frame}.
\item
   Compute the multiple error $\supertm{E^W}i m$ from \eqref{eq:errbd}.
\item
   Compute the right-hand side of the cohomological equations $\supertm\eta i
   m$ from \eqref{eq:defetabd}.
\item
   Solve \eqref{eq:corrbd3}, \eqref{eq:corrbd4}, \eqref{eq:corrbd1} for
   $\supertm{\xi^3}i m$, $\supertm{\xi^4}i m$, $\supertm{\xi^1}i m$,
   respectively, as multiple non-small divisors cohomological equations.
\item
   Take $\Delta\lambda=-\langle\eta^2\rangle$.
\item
   Take $\supertm{\xi^2}i m$ as the solution with $\langle\xi^2_0\rangle=0$
   of \eqref{eq:corrbd2} as small divisors multiple cohomological equation.
\item
   Compute the corrected multiple bundle as
   $\{W_i(\theta)+P_i(\theta)\xi_i(\theta)\}_{i=0}^{m-1}$ and the corrected
   eigenvalue as $\lambda+\Delta\lambda$.
\end{enumerate}
\end{algo}

\subsection{Algorithms for continuation}
\label{sec:continuation}

In this section we explain a me\-tho\-do\-lo\-gy to continue invariant tori
with respect to parameters. We will first consider continuation with respect
to time $T$ (the flying time, related with one of the frequencies of the
torus) and with respect to the energy $h$. The methodology can easily be
adapted for continuation of external parameters (appearing on the
Hamiltonian). In fact, we will only present an algorithm to compute the
tangent of the continuation curve with respect to the continuation parameter.
This provides a first order aproximation for the seed for the new value, as it
is common practice in numerical continuation (see e.g.~\cite{allgower-georg}).

Note that both approaches produce the same family of invariant tori, this is,
the one that corresponds to the fixed value of the rotation  vector $\omega$.
The difference is that the first approach is better to aim to a torus with a
specific value of return time $T$, whereas the second is better to aim to a
specific energy $h$. Aiming to a specific energy is useful to produce
isoenergetic Poincar\'e sections, which is a common way to represent the
center manifold of the collinear points of the RTBP (see
e.g.~\cite{1999JoMa,ESA3,GomezM01}). In this respect, we also consider at the
end of this section the continuation of the objects with respect to the
rotation vector $\omega$, in the isoenergetic case. Notice that in this case
the family of objects is parameterized by a Cantor set of parameters, the
rotation vector, and the derivatives to be computed are in Whitney sense. 

\subsubsection{Continuation with respect to $T$}
\label{sec:cntT}
In order to be able to perform continuation with respect to $T$, our goal now
is to compute the derivatives with respect to $T$ of the parameterizations of
tori and bundles that solve equations \eqref{eq:invtrms} and
\eqref{eq:invbdms}. Assume that these equations define implicitly
$\{K_i\}_{i=0}^{m-1}$, $\{W_i\}_{i=0}^{m-1}$ as functions of $T$. In order not
to burden the notation, we do not write the dependence on $T$ of $K_i$, $W_i$,
$\lambda$, but we denote by $\partial_TK_i(\theta)$, $\partial_TW_i(\theta)$,
$\partial_T\lambda$ their corresponding derivatives with respect to $T$.

By differentiating \eqref{eq:invtrms} with respect to $T$ we obtain
\begin{align*}
   \lefteqn{
      {\rm D}\varphi_{T/m}\bigl(K_i(\theta)\bigr)\partial_T K_i(\theta)
   }\hspace*{3em}\\
   &
   -\partial_T K_{i+1}(\theta+\tfrac\omega m) + \tfrac1m X_H\bigl(K_{i+1}(\theta+\tfrac\omega m)\bigr)
   =0.
\end{align*}
Considering the derivatives $\partial_T K_i(\theta)$ in the frame, this is,
assuming $\partial_T K_i(\theta)= P_i(\theta) \xi_i(\theta)$, the previous
equation is rewritten as
\begin{align*}
   \lefteqn{
      {\rm D}\varphi_{T/m}\bigl(K_i(\theta)\bigr)P_i(\theta) \xi_i(\theta)
   }\hspace*{3em}\\
   &
   -P_{i+1}(\theta+\tfrac\omega m)\xi_{i+1}(\theta+\tfrac\omega m)
   =
   -\tfrac1m X_H\bigl(K_{i+1}(\theta+\tfrac\omega m)\bigr).
\end{align*}
Multiplying by $P_{i+1}(\theta+\tfrac\omega m)^{-1}$, we obtain equation
\eqref{eq:corrtrmtr} but without the $e\Delta\tau$ term and with a different
definition for $\eta_i$, namely
\[
   \eta_i(\theta)=
      -\tfrac 1m P_{i+1}(\theta+\tfrac\omega m)^{-1}
      X_H\bigl(K_{i+1}(\theta+\tfrac\omega m)\bigr)
   = 
   -\tfrac1m \begin{pmatrix}
      e_{n-1} \\
      0 \\
      0_{n-1} \\
      0
   \end{pmatrix},
\]
for $i=0,\dots,m-1$.  It can be solved as described in subsection
\ref{sec:nwttr}, under isochronous nondegeneracy condition.

By differentiating \eqref{eq:invbdms} with respect to $T$, we obtain
\begin{align*}
   \lefteqn{
      \Bigl(
         \dparcial T \Bigl(
         {\rm D}\varphi_{T/m}\bigl(K_i(\theta)\bigr) \Bigr)
      \Bigr)W_i(\theta)
      +{\rm D}\varphi_{T/m}\bigl(K_i(\theta)\bigr)\partial_TW_i(\theta)
   }\hspace*{6em}
   \\
   &
   -(\partial_T\lambda)W_{i+1}(\theta+\tfrac\omega m)
   -\lambda\ \partial_TW_{i+1}(\theta+\tfrac\omega m)
   =0.
\end{align*}
The above expression is rewritten as
\begin{align*}
   \lefteqn{
      {\rm D}\varphi_{T/m}\bigl(K_i(\theta)\bigr) \partial_TW_{i}(\theta)
   }\hspace*{2em}
\\
   &
   -\lambda \ \partial_T W_{i+1}(\theta+\tfrac\omega m)
   -\partial_T\lambda \ W_{i+1}(\theta+\tfrac\omega m)
   =
   -E^W_i(\theta)
\end{align*}
with
\begin{equation}\label{eq:errcntTbd}
\begin{split}
   E^W_i(\theta)  & =
      \biggl(
         \dparcial T\Bigl(
            {\rm D}\varphi_{T/m}\bigl(K_i(\theta)\bigr)
         \Bigr)
      \biggr)W_i(\theta) \\ 
   & = 
      \tfrac\lambda{m}
    {\rm D} X_H\bigl(K_{i+1}(\theta+\tfrac\omega{m})\bigr)
      W_{i+1}(\theta+\tfrac\omega m)
\\
   & \phantom{=}
   +{\rm D}^2\varphi_{T/m}\bigl(K_i(\theta)\bigr)
   \Bigl[
      W_i(\theta), \partial_TK_i(\theta)
   \Bigr].
\end{split}
\end{equation}
Here $D^2\varphi_T\bigl(K_i(\theta)\bigr)\bigl[\cdot,\cdot\bigr]$ is the
bilinear form given by the second differential of $\varphi_T$ evaluated at
$K_i(\theta)$.  These equations are of the same type we have considered in
Section~\ref{sec:nwtbd}, by taking frames and writing $
\partial_TW_{i}(\theta)= P_i(\theta) \xi_i(\theta)$, $\partial_T \lambda=
\deltalambda$.

The steps to follow in order to solve the two systems of multiple
cohomological equations from which $\supertm{\partial_TK}i m$,
$\supertm{\partial_TW}i m$, $\partial_T\lambda$ can be obtained are summarized
in the algorithm that follows. In order not to burden the notation,
$\supertm\xi i m$ is used to denote the solution of both systems.

\begin{algo}\label{algo:cntT}
(Continuation step with respect to $T$) Let $\{K_i\}_{i=0}^{m-1}$,
$\{W_i\}_{i=0}^{m-1}$, $\lambda$ be implicit functions of $T$ through
equations \eqref{eq:invtrms}, \eqref{eq:energyms}, \eqref{eq:invbdms}. Find
$\supertm{\partial_T K}i m$, $\supertm{\partial_T W}i m$,
$\partial_T\lambda$\index{$\partial_T\lambda$} through the following steps:
\begin{enumerate}
\item
   Compute $\supertm P i m$, $\supertm{S^1}i m$, following
   Algorithm \ref{algo:frame}.
\item
   Take $\xi_i^2(\theta)= \xi_i^4(\theta)= 0$, $i=0,\dots,m-1$.
\item
   Compute $\langle S^1\rangle$ and $\xi_{0,0}^3= -\tfrac1m \langle S^1
   \rangle^{-1} e_{n-1}$.
\item
   Find $\supertm{\xi^1}i m$ as the solution with $\langle\xi^1_0\rangle=0$ of
   the small divisors multiple cohomological equation
   \[
           \xi_i^1(\theta) - \xi_{i+1}^1(\theta+\tfrac\omega{m}) = 
           -\tfrac1m e_{n-1}-S_i^1(\theta) \xi_{0,0}^3.
   \]
\item
   Obtain $\supertm{\partial_TK}i m$ as
   $\partial_TK_i(\theta)=P_i(\theta)\xi_i(\theta)$.
\item\label{en:cntTnumi}
   Compute $\supertm{E^W}i m$ from \eqref{eq:errcntTbd}.
\item
   Compute $\supertm\eta i m$ from \eqref{eq:defetabd}.
\item
   Solve \eqref{eq:corrbd3}, \eqref{eq:corrbd4}, \eqref{eq:corrbd1} for
   $\supertm{\xi^3}i m$, $\supertm{\xi^4}i m$, $\supertm{\xi^1}i m$,
   respectively, as multiple non-small divisors cohomological equations.
\item
   Take $\partial_T\lambda=-\langle\eta^2\rangle$.
\item
   Take $\supertm{\xi^2}i m$ as the solution with $\langle\xi^2_0\rangle=0$
   of \eqref{eq:corrbd2} as small divisors multiple cohomological equation.
\item
   Compute $\supertm{\partial_TW}i m$ as
   $\partial_TW_i(\theta)=P_i(\theta)\xi_i(\theta)$.
\end{enumerate}
\end{algo}

\begin{remark}
Step (2) of Algorithm~\ref{algo:cntT} (see also later
Algorithm~\ref{algo:cnth} and Algorithm~\ref{algo:cntomega}), set to zero the
components of the derivatives of the parameterizations of the tori in the
hyperbolic directions. This is geometrically very natural, since partially
hyperbolic tori are presented in families (for instance contained in a center
manifold of an equilibrium point, or in a normally hyperbolic cylinder), and
the derivatives are tangent to such families (and that invariant objects).
\end{remark}

\subsubsection{Continuation with respect to $h$}
\label{sec:cnh}

In this section we will consider the continuation of invariant tori with
respect to the energy $h$, so we have to compute the derivatives with respect
to $h$ of the parameterizations of tori, bundles and flying time $T$ that
solve equations \eqref{eq:invtrms}, \eqref{eq:energyms}, \eqref{eq:invbdms}.
Assume that these equations define implicitly $T$,
$\{K_i(\theta)\}_{i=0}^{m-1}$, $\{W_i(\theta)\}_{i=0}^{m-1}$  and $\lambda$ as
functions of $h$. We will follow the same criterion in the notation as in the
previous section, omiting explicitly the dependence on $h$ of these objects,
but we denote by $\partial_h T$, $\partial_h K_i(\theta)$, $\partial_h
W_i(\theta)$, $\partial_h \lambda$ their corresponding derivatives with
respect to $h$.

By differentiating \eqref{eq:invtrms} and \eqref{eq:energyms}
with respect to $h$ we obtain
\begin{align*}
   \lefteqn{
      {\rm D}\varphi_{T/m}\bigl(K_i(\theta)\bigr)\partial_h K_i(\theta)
   }\hspace*{3em}\\
   &
      -\partial_h K_{i+1}(\theta+\tfrac\omega m) + 
       \tfrac1m X_H\bigl(K_{i+1}(\theta+\tfrac\omega m)\bigr) \partial_h T
      =0
\end{align*}
and 
\[
	\langle {\rm D}H(K_0(\theta)) \partial_h K_0(\theta) \rangle - 1= 0.
\]
Considering the derivatives $\partial_h K_i(\theta)$ in the frame, this is,
assuming
$\partial_h K_i(\theta)= P_i(\theta) \xi_i(\theta)$, and denoting
$ \partial_h \tau= \tfrac1m \partial_h T$, the previous equations are rewritten as
\begin{align*}
\lefteqn{
   {\rm D}\varphi_{T/m}\bigl(K_i(\theta)\bigr)P_i(\theta) \xi_i(\theta)
}\hspace*{3em} \\
&
   -P_{i+1}(\theta+\tfrac\omega m)\xi_{i+1}(\theta+\tfrac\omega m)
   +X_H\bigl(K_{i+1}(\theta+\tfrac\omega m)\bigr) \partial_h \tau = 0.
\end{align*}
and
\[
   \langle {\rm D}H(K_0(\theta)) P_0(\theta)\xi_0(\theta) - 1 \rangle= 0.
\]
Multiplying by $P_{i+1}(\theta+\tfrac\omega m)^{-1}$, we obtain equation
\eqref{eq:corrtrmtr} with $\Delta\tau= \partial_h \tau$ and $\eta_i= 0$.  We
also obtain equation \eqref{eq:eninfrm} but with $1$ in place of $\Delta h$
and $0$ in place of $E_h$, that is $e_{n-1}^\top \xi_{0,0}^3= 1$.  This can be
solved as described in subsection \ref{sec:nwttr}, under isoenergetic
nondegeneracy condition. 

By differentiating \eqref{eq:invbdms} with respect to $h$, we obtain
\begin{align*}
   \lefteqn{
      {\rm D}\varphi_{T/m}\bigl(K_i(\theta)\bigr) \partial_h W_{i}(\theta)
   }\hspace*{3em} \\
&
   -\lambda \ \partial_h W_{i+1}(\theta+\tfrac\omega m)
   -\partial_h\lambda \ W_{i+1}(\theta+\tfrac\omega m)
   =
   -E^W_i(\theta)
\end{align*}
with
\begin{equation}\label{eq:errcnthbd}
\begin{split}
   E^W_i(\theta) & = 
      \biggl(
         \dparcial h\Bigl(
            {\rm D}\varphi_{T/m}\bigl(K_i(\theta)\bigr)
         \Bigr)
      \biggr) W_i(\theta) \\ 
      & = 
      \tfrac\lambda{m}
    {\rm D} X_H\bigl(K_{i+1}(\theta+\tfrac\omega{m})\bigr)
      W_{i+1}(\theta+\tfrac\omega m) \partial_h T
   \\
   &\phantom{=} 
   +{\rm D}^2\varphi_{T/m}\bigl(K_i(\theta)\bigr)
   \bigl[
      W_i(\theta), \partial_h K_i(\theta)
   \bigr]
\end{split}
\end{equation}

These equations are of the same type we have considered in
Section~\ref{sec:nwtbd}, by taking frames and writing $
\partial_hW_{i}(\theta)= P_i(\theta) \xi_i(\theta)$, $\partial_h \lambda=
\deltalambda$.

The steps to follow in order to solve the two systems of multiple
cohomological equations from which $\supertm{\partial_hK}i m$, $\partial_hT$,
$\supertm{\partial_hW}i m$, $\partial_h\lambda$ can be obtained are summarized
in the algorithm that follows. As in Section~\ref{sec:cntT}, we use
$\supertm\xi i m$ to denote the solution of both systems.

\begin{algo}\label{algo:cnth}
(Continuation step with respect to $h$) Let $\{K_i\}_{i=0}^{m-1}$, $T$,
$\{W_i\}_{i=0}^{m-1}$, $\lambda$ be implicit functions of $h$ through
equations \eqref{eq:invtrms}, \eqref{eq:energyms}, \eqref{eq:invbdms}. Find
$\supertm{\partial_h K}i m$, $\partial_h T$, $\supertm{\partial_h W}i m$,
$\partial_h\lambda$ through the following steps:
\begin{enumerate}
\item
   Compute $\supertm P i m$, $\supertm{S^1}i m$, following
   Algorithm \ref{algo:frame}.
\item
   Set $\xi_i^2(\theta)= \xi_i^4(\theta)= 0$ and $\xi_i^3(\theta)=
   \xi_{0,0}^3$ for $i= 0,\dots, m-1$, where 
   \[
      \begin{pmatrix}
      \xi_{0,0}^3 \\ \partial_h \tau
      \end{pmatrix}
      = 
      \begin{pmatrix}
      \langle {S^1} \rangle & e_{n-1} \\
      e_{n-1}^\top & 0 
      \end{pmatrix}^{-1}
      \begin{pmatrix}
      0_{n-1} \\
      1
      \end{pmatrix}
   \]	
\item
   Take $\{\xi_i^1\}_{i= 0}^{m-1}$ as the solution of
   \[
        \xi_i^1(\theta) - \xi_{i+1}^1(\theta+\tfrac\omega{m}) = 
        - e_{n-1} \partial_h \tau -S_i^1(\theta) \xi_{0,0}^3,
   \]
   with $\langle \xi_0^1 \rangle= 0$, and set $\partial_h T= m\ \partial_h
   \tau$.
\item
   Obtain $\supertm{\partial_hK}i m$ as
   $\partial_hK_i(\theta)=P_i(\theta)\xi_i(\theta)$.
\item\label{en:cnthnumi}
   Compute $\supertm{E^W}i m$ from \eqref{eq:errcnthbd}.
\item
   Compute $\supertm\eta i m$ from \eqref{eq:defetabd}.
\item
   Solve \eqref{eq:corrbd3}, \eqref{eq:corrbd4}, \eqref{eq:corrbd1} for
   $\supertm{\xi^3}i m$, $\supertm{\xi^4}i m$, $\supertm{\xi^1}i m$,
   respectively, as multiple non-small divisors cohomological equations.
\item
   Take $\partial_h\lambda=-\langle\eta^2\rangle$.
\item
   Take $\supertm{\xi^2}i m$ as the solution with $\langle\xi^2_0\rangle=0$
   of \eqref{eq:corrbd2} as small divisors multiple cohomological equation.
\item
   Compute $\supertm{\partial_hW}i m$ as
   $\partial_hW_i(\theta)=P_i(\theta)\xi_i(\theta)$.
\end{enumerate}
\end{algo}

\subsubsection{Continuation with respect to $\omega$, in the isoenergetic case}

Another natural continuation problem is, given a certain energy level $h$,
continue the invariant tori on such an energy level and their invariant
bundles with respect to the frequencies. Assume that  equations
\eqref{eq:invtrms}, \eqref{eq:energyms}, \eqref{eq:invbdms} define implicitly
$T$, $\{K_i(\theta)\}_{i=0}^{m-1}$, $\{W_i(\theta)\}_{i=0}^{m-1}$  and
$\lambda$ as functions of $\omega$. We will again omit explicitly the
dependence on the rotation vector $\omega$ of these objects, but we denote by
$\partial_\omega T$, $\partial_\omega K_i(\theta)$, $\partial_\omega
W_i(\theta)$, $\partial_\omega \lambda$ their corresponding derivatives with
respect to $\omega$. Since the domain of the rotation vector is a Cantor set,
these derivatives are understood formally (in fact, these are derivatives in
the sense of Whitney, but we will avoid technicalities here).

By differentiating \eqref{eq:invtrms} and \eqref{eq:energyms} with respect to
$\omega$ we obtain
\begin{align*}
   \lefteqn{
      {\rm D}\varphi_{T/m}\bigl(K_i(\theta)\bigr)\partial_\omega K_i(\theta)
   }\hspace*{3em}\\
   &
      -\partial_\omega K_{i+1}(\theta+\tfrac\omega m) + 
       \tfrac1m X_H\bigl(K_{i+1}(\theta+\tfrac\omega m)\bigr) \partial_\omega T
      = \tfrac{1}{m} {\rm D}K_{i+1}(\theta+\tfrac\omega m)
\end{align*}
and 
\[
     \langle {\rm D}H(K_0(\theta)) \partial_\omega K_0(\theta) \rangle= 0.
\]
Considering $\partial_\omega K_i(\theta)= P_i(\theta) \xi_i(\theta)$, denoting
$ \partial_\omega \tau= \tfrac1m \partial_\omega T$, 
and multiplying by $P_{i+1}(\theta+\tfrac\omega m)^{-1}$, we obtain equation
\eqref{eq:corrtrmtr} with $\Delta\tau= \partial_\omega \tau$ and 
\[
   \eta_i(\theta)= \tfrac{1}{m} \begin{pmatrix} I_{n-2} \\ O \end{pmatrix}.
\]
Moreover, equation \eqref{eq:eninfrm} reads $e_{n-1}^\top \xi_{0,0}^3= 0$.
This can be solved as described in subsection \ref{sec:nwttr}, under
isoenergetic nondegeneracy condition. 

By differentiating \eqref{eq:invbdms} with respect to $\omega$, we obtain
\begin{align*}
   \lefteqn{
      {\rm D}\varphi_{T/m}\bigl(K_i(\theta)\bigr) \partial_\omega W_{i}(\theta)
   }\hspace*{3em} \\
&
   -\lambda \ \partial_\omega W_{i+1}(\theta+\tfrac\omega m)
   -\partial_\omega\lambda \ W_{i+1}(\theta+\tfrac\omega m)
   =
   -E^W_i(\theta)
\end{align*}
with
\begin{equation}\label{eq:errcntomegabd}
\begin{split}
   E^W_i(\theta) & = 
      \biggl(
         \dparcial\omega
         \Bigl(
            {\rm D}\varphi_{T/m}\bigl(K_i(\theta)\bigr)
         \Bigr)
      \biggr) W_i(\theta) \\ 
      & = 
      \tfrac\lambda{m}
    {\rm D} X_H\bigl(K_{i+1}(\theta+\tfrac\omega{m})\bigr)
      W_{i+1}(\theta+\tfrac\omega m) \partial_\omega T
   \\
   &\phantom{=} 
   +{\rm D}^2\varphi_{T/m}\bigl(K_i(\theta)\bigr)
   \bigl[
      W_i(\theta), \partial_\omega K_i(\theta)
   \bigr] 
   \\
   &\phantom{=} 
   - \tfrac{\lambda}{m}  {\rm D}W_{i+1}(\theta+\tfrac\omega{m}).
\end{split}
\end{equation}
These equations are of the same type we have considered in
Section~\ref{sec:nwtbd}, by taking frames and writing $ \partial_\omega
W_{i}(\theta)= P_i(\theta) \xi_i(\theta)$, $\partial_\omega \lambda=
\deltalambda$.

Similarly as we proceed in previous sections, we present in an algorithm the
steps  to solve the two systems of multiple cohomological equations to compute
$\supertm{\partial_\omega K}i m$, $\partial_\omega T$,
$\supertm{\partial_\omega W}i m$, $\partial_\omega\lambda$. Notice that these
are equations for the partial derivatives with respect to the components of
$\omega$, that we formulate  in paralell.

\begin{algo}\label{algo:cntomega}
(Continuation step with respect to $\omega$, for energy $h$ fixed) Let
$\{K_i\}_{i=0}^{m-1}$, $T$, $\{W_i\}_{i=0}^{m-1}$, $\lambda$ be implicit
functions of $h$ through equations \eqref{eq:invtrms}, \eqref{eq:energyms},
\eqref{eq:invbdms}. Find $\supertm{\partial_\omega K}i m$, $\partial_\omega
T$, $\supertm{\partial_h W}i m$, $\partial_\omega\lambda$ through the
following steps:
\begin{enumerate}
\item
   Compute $\supertm P i m$, $\supertm{S^1}i m$, following
   Algorithm \ref{algo:frame}.
\item
   Set $\xi_i^2(\theta)= \xi_i^4(\theta)= 0$ and $\xi_i^3(\theta)=
   \xi_{0,0}^3$ for $i= 0,\dots, m-1$, where 
   \[
      \begin{pmatrix}
      \xi_{0,0}^3 \\ \partial_\omega \tau
      \end{pmatrix}
      = 
      \begin{pmatrix}
      \langle {S^1} \rangle & e_{n-1} \\
      e_{n-1}^\top & 0 
      \end{pmatrix}^{-1}
      \begin{pmatrix}
      \frac1m I_{n-2} \\ 0 \\ \hline
      0
      \end{pmatrix}
   \]	
\item
   Take $\{\xi_i^1\}_{i= 0}^{m-1}$ as the solution of
   \[
        \xi_i^1(\theta) - \xi_{i+1}^1(\theta+\tfrac\omega{m}) = \tfrac1m \begin{pmatrix}
   I_{n-2} \\ 0  \end{pmatrix}
        - e_{n-1} \partial_\omega\tau -S_i^1(\theta) \xi_{0,0}^3,
   \]
   with $\langle \xi_0^1 \rangle= 0$, and set $\partial_\omega T= m\ \partial_\omega \tau$.
\item
   Obtain $\supertm{\partial_\omega K}i m$ as
   $\partial_\omega K_i(\theta)=P_i(\theta)\xi_i(\theta)$.
\item\label{en:cntomeganumi}
   Compute $\supertm{E^W}i m$ from \eqref{eq:errcntomegabd}.
\item
   Compute $\supertm\eta i m$ from \eqref{eq:defetabd}.
\item
   Solve \eqref{eq:corrbd3}, \eqref{eq:corrbd4}, \eqref{eq:corrbd1} for
   $\supertm{\xi^3}i m$, $\supertm{\xi^4}i m$, $\supertm{\xi^1}i m$,
   respectively, as multiple non-small divisors cohomological equations.
\item
   Take $\partial_\omega\lambda=-\langle\eta^2\rangle$.
\item
   Take $\supertm{\xi^2}i m$ as the solution with $\langle\xi^2_0\rangle=0$
   of \eqref{eq:corrbd2} as small divisors multiple cohomological equation.
\item
   Compute $\supertm{\partial_\omega W}i m$ as
   $\partial_\omega W_i(\theta)=P_i(\theta)\xi_i(\theta)$.
\end{enumerate}
\end{algo}

\begin{remark}
In the implementation of the continuation with respect to rotation vector, one
have to select continuation steps for which the rotation vectors are
Diophantine.  In practice, since in the computer the rotation vectors are
rational, these has to be selected as resonant but of very high order. The
continuation can run into troubles when finding strong resonances, since these
are more difficult to jump.
\end{remark}

\subsection{Some comments about the implementation}
\label{sec:compimpl}

As it is common in implementation of the parameterization method in KAM-like
contexts (see \cite{HaroCFLM16} for an overview), the implementation of all
the previous algorithms relies on two numerical representations of all the
functions $\zeta:\T^d\rightarrow\R$ involved. In the grid representation, the
function is represented as a set of its values in a uniform grid of $\T^d$.
In the Fourier representation, the function is represented as a set of
approximate Fourier coefficients. Through the Discrete Fourier Transform
(DFT), that provides approximations of the Fourier coefficients, one
representation can be converted into the other. For an easier exposition, we
will assume $d=1$ in this section (which is the case in all the numerical
exploration of Section~\ref{sec:testex}). All the arguments generalize to
$d>1$, although an actual implementation is subtle (see \cite{HaroCFLM16} for
comments).  For detailed expositions on the DFT and its applications, see
e.g.~\cite{gasquet, brigham, nrcpp}.

Let $\zeta:\T\rightarrow\R$ be a function. Choose $N>0$. Its grid
representation is given by $\{\zeta_j\}_{j=0}^{N-1}$, with
$\zeta_j=\zeta(j/N)$. Its Fourier representation is given by a finite set of
(complex) approximate Fourier coefficients, $\{\tilde\zeta_k\}_{k=0}^{[N/2]}$,
where $[\cdot]$ denotes integer part. The two representations are related by
the DFT,
\[
   \tilde\zeta_k=\frac1N\sum_{j=0}^{N-1}\zeta_je^{-\bmi 2\pi  k\frac jN}.
\]
The DFT is a linear, one-to-one map between $\{\zeta_j\}_{j=0}^{N-1}$ and
$\{\tilde\zeta_k\}_{k=0}^{N-1}$. Namely, the expression for the inverse DFT is
\[
   \zeta_j=\sum_{k=0}^{N-1}\tilde\zeta_k e^{\bmi 2\pi  k\frac jN}
   .
\]
The DFT and its inverse are efficiently evaluated through a family of
algorithms known as Fast Fourier Transform (FFT), which allow computing
$\{\tilde\zeta_k\}_{k=0}^{N-1}$ from $\{\zeta_j\}_{j=0}^{N-1}$ in $O(N\log N)$
operations. The DFT is $N$-periodic in $k$, and, for real $\zeta$ (which is
our case), satisifes the Hermitian symmetry. This is, for $k\in\Z$,
\[
   \tilde\zeta_k=\tilde\zeta_{k+N}
   ,\quad
   \tilde\zeta_{-k}=(\tilde\zeta_k)^*
   ,
\]
where ${}^*$ denotes complex conjugate. This gives rise to redundancy in
$\{\tilde\zeta_k\}_{k=0}^{N-1}$, which is eliminated by truncating the Fourier
representation at $k=[N/2]$, as we have done. Another consequence of the
Hermitian symmetry is that $\tilde\zeta_0$ is real and, if $N$ is even,
$\tilde\zeta_{N/2}$ is also real. The DFT coefficients and the Fourier
coefficients are related by
\[
   \tilde\zeta_k=\hat\zeta_k
      +\sum_{l=1}^\infty(\hat\zeta_{k-lN}+\hat\zeta_{k+lN})
   .
\]
Actual bounds of the difference $\tilde\zeta_k-\hat\zeta_k$ can be obtained
from the Cauchy estimates, that ensure exponential decay in $|k|$ of
$|\hat\zeta_k|$ for analytic $\zeta$ (see e.g.~\cite{Henrici79}).  This is the
basis for a detailed error analysis of the algorithms, that we have not
pursued here (see e.g.~\cite{FiguerasHL17}), but a direct consequence of this
approximation is that, for the algorithms to work, $N$ has to be large enough
for $|\tilde\zeta_k-\hat\zeta_k|$ to be small. In the computations of
Section~\ref{sec:testex}, the minimum value of $N$ used is 32. Another
consequence is that the relative error of the DFT approximation of the Fourier
coefficients increases with $k$. Even with (large) DFT queues of the order of
the machine epsilon, we have observed instability in Newton iterates
(divergence after apparent convergence). In our implementation, we prevent it
by setting to zero the upper half of the DFT coefficients after each Newton
iteration.

All the steps of algorithms \ref{algo:frame}, \ref{algo:nwttr},
\ref{algo:nwtbd}, \ref{algo:cntT}, \ref{algo:cnth} can be done in $O(N)$
operations in either grid or Fourier representation. Consider, for instance,
Algorithm~\ref{algo:frame}. In the evaluation of \eqref{eq:defLi},
$W_i(\theta)$ is a function we already have, $X_H(K_i(\theta))$ is obtained
from $K_i(\theta)$ in $O(N)$ operations in grid form, as
\[
   \{X_H(K_i(j/N))\}_{j=0}^{N-1}
   ,
\]
and ${\rm D}K_i(\theta)$ obtained in $O(N)$ operations in Fourier form from
$K_i(\theta)$ as
\[
   (\widetilde{{\rm D}K_i})_k=\bmi 2\pi k(\widetilde{K_i})_k
   .
\]
Note that this last expression is actually an approximation: the equality is
satisfied by the Fourier coefficients of $K_i$ and ${\rm D}K_i$, not the DFT
ones. As another example, the evaluation of \eqref{eq:defSi} is done in
$O(N)$ operations in grid form if $\hat N_{i+1}(\theta+\frac\omega m)$,
$K_{i+1}(\theta+\frac\omega m)$ and ${\rm D}\varphi_{T/m}(K_i(\theta))\hat
N_i(\theta)$ are known in grid form. The computation of
$K_{i+1}(\theta+\frac\omega m)$ from $K_{i+1}(\theta)$ is done in $O(N)$
operations in Fourier form (again using approximate identities). The
computation of $D\varphi_{T/m}(K_i(\theta))\hat N_i(\theta)$ from
$K_i(\theta)$ and $\hat N_i(\theta)$ in grid form requires the numerical
integration of the differential equations \eqref{eq:hamsys} with his first
variationals applied to several vectors on $N$ trajectories.

The different steps of the algorithms stated consist of evaluating equations
like the ones just mentioned and solving multiple cohomological equations.
The solution of multiple cohomological equations discussed in
Section~\ref{sec:coho} is also done in $O(N)$ operations in Fourier form. As
mentioned, converting from grid to Fourier and vice-versa using FFT requires
$O(N\log N)$ operations. For the values of $N$ for which numerical integration
is feasible, $\log N$ is small enough for the cost of FFT to be considered
$O(N)$.  Numerical integration is present in step~\ref{en:frnumi} of
Algorithm~\ref{algo:frame} (Eq.~\eqref{eq:defSi}), step~\ref{en:cntTnumi} of
Algorithm~\ref{algo:cntT} (Eq.~\eqref{eq:errcntTbd}), step~\ref{en:cnthnumi}
of Algorithm~\ref{algo:cnth} (Eq.~\eqref{eq:errcnthbd}) and
step~\ref{en:cntomeganumi} of Algorithm~\ref{algo:cntomega}
(Eq.\eqref{eq:errcntomegabd}). Note that Eqs.~\eqref{eq:errcntTbd},
\eqref{eq:errcnthbd}, \eqref{eq:errcntomegabd} require numerical integration
of the second variational equations. Note also that all algorithms perform
numerical integration when computing the frame through
Algorithm~\ref{algo:frame}.

The cost of numerical integration is formally also $O(N)$, but for realistic
estimates it needs to be considered separatedly. Consider for instance
Equation~\eqref{eq:defSi} as discussed above. In the computations of
Section~\ref{sec:testex}, $n= 3$, so $\hat N_i(\theta)$ has $3$ columns, and
the system of differential equations to be numerically integrated consists of
$24$ equations.  Taking for example the Runge-Kutta-Felhberg method of orders
7 and 8 for numerical integration, that evaluates the vector field $13$ times,
since orbits in Section~\ref{sec:testex} take over $50$ integration steps, the
factor multiplying $N$ is already $15600$ times the number of operations
required by the evaluation of one differential equation. It turns out that,
compared to numerical integration, the remaining computational cost is almost
negligible. Numerical integration can be parallelized in a straightforward
manner by distributing the $N$ trajectories to be integrated among the threads
or processes available.

\newcommand{\pT}{\partial_T}
\newcommand{\ph}{\partial_h}
\newcommand{\cT}{{\mathcal T}}
\newcommand{\err}{\mathop{\rm err}\nolimits}
\newcommand{\nit}{n_{\it it}}
\newcommand{\ndes}{n_{\it des}}

An actual implementation of a continuation procedure in order to follow a
family of tori by keeping $\omega$ constant requires an strategy in order to
choose the number of samples $N$ and control the continuation step size as we
go along the family. We end this section
with a proposal of such a strategy in Algorithm~\ref{algo:cntstep}, that has
worked well in the numerical computations presented in
Section~\ref{sec:testex}. For shortness, we will represent a multiple torus,
its multiple bundle, flying time and eigenvalue as
\[
   \cT=(K_0,\dots,K_{m-1},W_0,\dots,W_{m-1},T,\lambda).
\]
For step size control, the following norm of this compound object is
considered:
\[
   \|\cT\|
   =
   \Bigl(
      T^2+\lambda^2
         +\sum_{i=0}^{m-1}\bigl(
            \langle\|K_i\|_2^2\rangle+\langle\|W_i\|_2^2\rangle
   \Bigr)^{1/2}
   ,
\]
where $\|K_i\|_2$ stands for the function $\theta\mapsto\|K_i(\theta)\|_2$,
$\|W_i\|_2$ is interpteted analogously, and the averages are approximated as
discrete averages of the grid values (i.e~as the $0$-th DFT coefficient of the
function averaged). The errors in the torus and the bundle are estimated as
\[
   \err(\cT)=\max_{\substack{
   0\leq i<m\\
   0\leq j<N
   }}\|E_i(j/N)\|_\infty
   ,\quad
   \err^W(\cT)=\max_{\substack{
   0\leq i<m\\
   0\leq j<N
   }}\|E^W_i(j/N)\|_\infty
   ,
\]
with $E(\theta)$, $E^W(\theta)$ defined as in Eqs.~\eqref{eq:errtr},
\eqref{eq:errbd}, respectively.

\begin{algo}\label{algo:cntstep}
   Assume we are given a multiple torus, bundle and associated parameters
   $\cT$, represented in grid form with $N$ samples. Assume we are also given
   a suggested continuation step $\alpha$, tolerances
   $\varepsilon,\varepsilon^W,\varepsilon_1,\varepsilon_2$ and integers
   $\ndes$, $n_\alpha$. Perform a continuation step in order to obtain a new
   $\cT$ along the corresponding family with the same $\omega$ as follows:
\begin{enumerate}
\item\label{en:cnt:lbl2}
   Set $\delta\leftarrow(\pT K_0,\dots,\pT K_{m-1},\pT W_0,\dots,\pT
   W_{m-1},1,\pT\lambda)$, using Algorithm~\ref{algo:cntT}.
\item
   Set $\Delta\cT\leftarrow\delta/\|\delta\|$.
\item
   Set $\cT'\leftarrow\cT+\alpha\Delta\cT$ and perform Newton steps
   (Algorithms~\ref{algo:nwttr}, \ref{algo:nwtbd}) on $\cT'$ until
   $\err(\cT')<\varepsilon$ and $\err^W(\cT')<\varepsilon^W$. If unsuccessful,
   try halving $\alpha$ up to $n_\alpha$ times. If unsuccessful, restart the
   algorithm by doubling $N$. Let $\nit$ be the number of iterates of the
   successful Newton iteration.
\item
   If $\err(\cT')<\varepsilon_2$, half $N$ as long as
   $\err(\cT')\leq \varepsilon_1$. Go to step \ref{en:cnt:lbl1}. 
\item
   If $\err(\cT')>\varepsilon_1$, double $N$ and perform Newton steps
   (Algorithms~\ref{algo:nwttr}, \ref{algo:nwtbd}) until
   $\err(\cT')\leq \varepsilon_1$. If unsuccessful, half $\alpha$ and restart
   the algorithm.
\item\label{en:cnt:lbl1}
   Set $\cT\leftarrow\cT'$ (i.e., accept the new torus). Perform continuation
   step size control as $\alpha\leftarrow\alpha\ndes/\nit$.
\end{enumerate}
In step \ref{en:cnt:lbl2}, Algorithm~\ref{algo:cnth} can be used instead, by
setting
\[
   \delta\leftarrow(\ph K_0,\dots,\ph K_{m-1},\ph W_0,\dots,\ph W_{m-1},\ph
   T,\ph\lambda).
\]
\end{algo}

\section{An application: computation of the Lissajous family of tori in the 
Restricted Three Body Problem}
\label{sec:application}

\newcommand{\sci}[2]{#1\times10^{#2}}
\newcommand{\Spec}{\mathop{\rm Spec}\nolimits}

\label{sec:testex}

In this section we apply the algorithms described to the computation of the
partially hyperbolic tori that emerge from the $L_1$ equilibrium point of the
circular, Spatial Restricted Three-body Problem (RTBP), in the Earth-Moon
case.  This family of tori is known as the Lissajous family by the
astrodynamics community, and plays a fundamental role in libration point
dynamics.  An outgrowth of our algorithms is the simultaneous computation of
the stable, unstable and center bundles of the invariant tori and other
observables, such as the Lyapunov multipliers and Calabi invariant, that
provide geometrical and dynamical information during the computation of the
families of tori.  In particular, we will obtain information about the quality
of hyperbolicity properties, that is useful to detect possible bifurcations
and breakdown phenomena \cite{CallejaF12,CanadellH17b,HaroL06b}.

Apart from the intrinsic importance of the example, our choice is motivated by
the fact that the Lissajous family has been extensively described in
\cite{GomezM01}, and has been used as a testbed for different algorithms by
other authors (see e.g.~\cite{2018aBaOlSche}).  Hence, one of our goals is  to
compare the performance of the algorithms described here with the ones used in
\cite{GomezM01}, on a more thorough numerical exploration of this family. 

All the numerical explorations presented here have been done by a program that
follows a family of tori with constant $\omega$ by performing continuation
steps through Algorithm~\ref{algo:cntstep} up to a maximum given number, plus
additional stopping criteria that will be specified below. The explorations
have been carried out on a Fujitsu Celsius R940 workstation, with two 8-core
Intel Xeon E5-2630v3 processors at 2.40GHz, running Debian GNU/Linux 9.11 with
the Xfce 4.12 desktop. The source code has been written in C, compiled with
GCC 6.3.0 and linked against the Glibc 2.24, LAPACK 3.7.0, FFTW 3.3.5 and
PGPLOT 5.2.2 libraries. The code uses OpenMP 4.0 extensions in order to
perform simultaneous numerical integrations in the continuation of one family
and also to perform numerical continuation of several families at once. The
figures have been generated with gnuplot 5.2.

\subsection{The Lissajous family of tori of the RTBP}

The circular, spatial Restricted Three-Body Problem (RTBP) describes the
motion of a particle of infinitesimal mass under the attraction of two massive
bodies known as primaries, with masses $m_1>m_2>0$. The primaries are assumed
to revolve uniformly in circles around their common center of mass. In a
rotating system of reference with the primaries in the horizontal coordinate
plane, known in astronomical terms as synodic, the primaries can be made to
lie at fixed positions in the $x_1$ axis. After a rescaling in space and time,
and defining the mass ratio $\mu=\tfrac{m_2}{m_1+m_2}$, the coordinates of the
primaries $m_1,m_2$ become $(\mu,0,0),(\mu-1,0,0)$, their masses become
$1-\mu$, $\mu$ respectively, and their period of revolution becomes $2\pi$.
The motion of the infinitesimal mass is then described by the autonomous
Hamiltonian system with Hamiltonian
\[
   H(x_1,x_2,x_3,p_1,p_2,p_3)=\frac12(p_1^2+p_2^2+p_3^2)
      -x_1p_2+x_2p_1-\frac{1-\mu}{r_1}-\frac\mu{r_2},
\]
where $r_1^2=(x_1-\mu)^2+x_2^2+x_3^2$, $r_2^2=(x_1-\mu+1)^2+x_2^2+x_3^2$. The
value of the hamiltonian will be denoted as ``the energy'' from now on.

The RTBP is shown to have 5 fixed points: the collinear ones, $L_1,L_2,L_3$,
due to Euler, and the triangular ones, $L_4,L_5$, due to Lagrange (see
e.g.~\cite{Szebehely}). Following the astrodynamical convention, we will
consider $L_1$ to be the point located between the primaries.  The $x$
coordinate of this point is $x_{L_1}=\mu-1+\gamma_1$, with $\gamma_1$ the
positive root of one of Euler's quintic equations,
\[
   \gamma_1^5-(3-\mu)\gamma_1^4+(3-2\mu)\gamma_1^3-\mu\gamma_1^2
      +2\mu\gamma_1-\mu = 0.
\]
The linear behaviour around $L_1$ is of the type
center$\times$center$\times$saddle. Namely, for the value of $\mu$ we use,
\[
   \Spec {\rm D}X_H(L_1) = \{
      \bmi 2\pi\omega_p^0, -\bmi2\pi\omega_p^0,
      \bmi 2\pi\omega_v^0, -\bmi2\pi\omega_v^0,
      \lambda^0, -\lambda^0
   \}.
\]
From now on we will focus our attention on this point, and we will consider
the primaries to be the Earth and the Moon, with mass parameter
$\mu=\sci{1.215058560962404}{-2}$, for which $x_{L_1}\approx-0.83692$,
$\omega_p^0\approx0.371529, \omega_v^0\approx0.361096,
\lambda^0\approx2.932056$.

Lyapunov's center theorem (see e.g.~\cite{meyer-hall-offin,siegel-moser})
ensures the existence of a family of periodic orbits (p.o.), known as the
planar (resp.~vertical) Lyapunov family, that fills a 2D manifold tangent to
the $\pm\bmi2\pi\omega^0_p$ (resp.~$\pm\bmi2\pi\omega^0_v$) eigenplane.  The
planar (resp.~vertical) denomination comes from the fact that the eigenvectors
of eigenvalues $\pm\bmi2\pi\omega^0_p$ (resp. $\pm\bmi 2\pi\omega^0_v$) have
zero $x_3,p_3$ (resp.~$x_1,x_2,p_1,p_2$) coordinates. Both families start at
the energy of $L_1$, that will be denoted as $h_0$, and evolve through higher
energies. Denote by $T_p^h$ (resp.~$T_v^h$) the period of the planar
(resp.~vertical) Lyapunov p.o.~of energy $h$. Denote also as $e^{\pm \bmi
2\pi\nu_p^h}$ (resp.~$e^{\pm \bmi 2\pi\nu_v^h}$) the multipliers of modulus
one of the monodromy matrix of the planar (resp.~vertical) Lyapunov p.o.~of
energy $h$, with $\nu_p^h$ (resp.~$\nu_v^h$) chosen in $[0,1/2]$, as found
when computing numerically. This is,
\[
   e^{\pm \bmi 2\pi\nu_j^h}\in\Spec {\rm D}\varphi_{T_j^h}(x_j^h)
   ,\quad
   \nu_j^h\in[0,1/2]
   ,\quad
   j=p,v
   ,
\]
where $x^h_p$ (resp.~$x^h_v$) is an initial condition in the planar
(resp.~vertical) periodic orbit of energy $h$. Lyapunov's center theorem also
ensures that
\[
   T_p^h\stackrel{h\rightarrow h_0}\longrightarrow1/\omega_p^0
   ,\quad
   T_v^h\stackrel{h\rightarrow h_0}\longrightarrow1/\omega_v^0
   ,
\]
and
\[
   \{e^{\pm\bmi2\pi\nu_p^h}\}
      \stackrel{h\rightarrow h_0}\longrightarrow
      \{e^{\pm\bmi2\pi\omega^0_v/\omega^0_p}\}
   ,\quad
   \{e^{\pm\bmi2\pi\nu_v^h}\}
      \stackrel{h\rightarrow h_0}\longrightarrow
      \{e^{\pm\bmi2\pi\omega^0_p/\omega^0_v}\}
   .
\]
From the numerical values of $\omega_p^0,\omega_v^0$, we have
\[
   \nu_p^h
      \stackrel{h\rightarrow h_0}\longrightarrow
      1-\omega^0_v/\omega^0_p
   ,\quad
   \nu_v^h
      \stackrel{h\rightarrow h_0}\longrightarrow
      \omega^0_p/\omega^0_v-1
   .
\]

The Lissajous family of tori mentioned above is made of the KAM tori generated
by the $4D$ central part of $L_1$, that are contained inside the 4D center
manifold of this point.  Denote by $\omega_p,\omega_v$ the frequencies of any
torus in the family, chosen as to have $\omega_p\rightarrow\omega_p^0$ and
$\omega_v\rightarrow\omega_v^0$ as the torus collapses to $L_1$. We will refer
to  $\hat\omega=(\omega_p,\omega_v)$ as the vector of natural frequencies.
Denote by $\hat K$ a parameterization of a torus $\mathcal {\hat K}$ of the
family satisfying the invariance equation \eqref{eq:invtrnm0}.  Define
\begin{equation}\label{eq:defnupnuvtori}
   \nu_p(\omega_p,\omega_v)=1-\omega_v/\omega_p
   ,\quad
   \nu_v(\omega_p,\omega_v)=\omega_p/\omega_v-1
   .
\end{equation}
As stated previously, we will not compute a parameterization $\hat K$ of the
whole torus, but of an invariant curve parameterized by $K$ inside it (recall
that we actually compute a collection of invariant curves). 

\begin{figure}[htbp]
\includegraphics{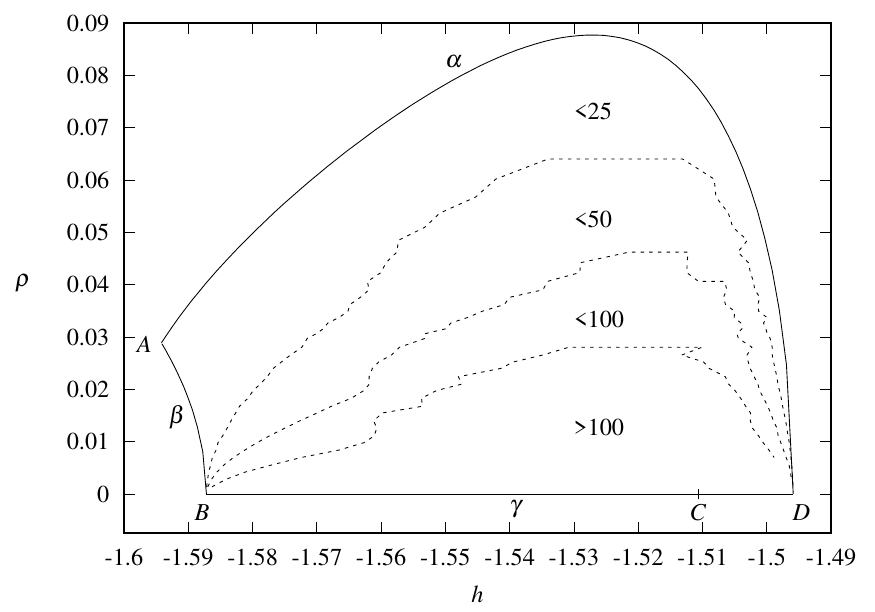}
\caption{\label{fig:enrho}(Adapted from \cite{GomezM01}) Energy-rotation
number representation of the Lissajous family of invariant tori of the RTBP
around $L_1$ for the Earth-Moon mass parameter.}
\end{figure}

Following  \cite{GomezM01}, we will use the energy $h$ and the vertical
rotation number $\rho= \nu_v$ as parameters in order to represent the tori of
the Lissajous family. It was numerically found that, when varying $h,\nu_v$ in
the region enclosed by the $\alpha,\beta,\gamma$ curves of
Fig.~\ref{fig:enrho}, they  uniquely determine a torus in the family.  The
$\alpha$ curve, with coordinates $(h,\nu_v^h)$, represents the vertical
Lyapunov family, from its birth at $L_1$ at energy $h=h_0$ (point $A$) to its
first 1:1 bifurcation, at energy $h=-1.49590$ (point $D$). The $\beta$ curve
represents the planar Lyapunov family from its birth to its first 1:1
bifurcation, at energy $h=h_B=-1.58718$ (point $B$, in which the Halo family
of p.o.~appears).  In order to have continuity at the point $A$, the vertical
coordinate of the points of the $\beta$ curve is not $\nu_p^h$ but the
$\nu_v(\omega_p,\omega_v)$ value of limiting nearby tori.  From
\eqref{eq:defnupnuvtori}, it is found to be
\[
   \frac{1}{1-\nu_p^h}-1.
\]
The $\gamma$ curve, which is the segment from point $B$ to point $D$,
corresponds to the separatrix between the Lissajous family of tori and other
families of quasi-periodic motion in the center manifold of $L_1$, that are
described in \cite{GomezM01}. 

As it has been mentioned in the introduction, the large
matrix\footnote{According to the nomenclature of \cite{HaroL06b}.} approach in
\cite{GomezM01} was to write $K(\theta)$ as a truncated Fourier series,
$K(\theta)= A_0+\sum_{k=1}^{N_f}\bigl(A_k\cos(2\pi k\theta)+B_k\sin(2\pi
k\theta)\bigr)$, an then turn Eq.~\eqref{eq:invtrnm2} into a finite non-linear
system of equations by imposing it at $1+2N_f$ equally spaced values of
$\theta$. Multiple shooting was also implemented: all the computations were
done with $m=2$. The computational bottleneck of this procedure is that large
values of $N_f$ give rise to large systems of equations. The different
sub-regions inside the $\alpha,\beta,\gamma$ curves in Fig.~\ref{fig:enrho},
that are not disjoint but nested, are labeled according to the value of $N_f$
obtained in the computations of \cite{GomezM01} for the tori inside them.  A
global upper limit of $100$ was chosen for $N_f$, so tori in the ''$>100$''
sub-region were actually not computed.

\subsection{On the generators of tori in the Lissajous family}
\label{ssection:vertical and planar}

As stated previously, we will not compute a parameterization $\hat K$ of the
whole torus, but of an invariant curve parameterized by $K$ inside it (recall
that we actually compute a collection of invariant curves), to which we will
refer to as a generator of the torus. The choice has consequences in the
determination of the geometrical observables and the selection of the
frequencies. We consider two cases, that we will distinguish as vertical
generators and planar generators.

We will denote as {\em vertical generator} of the invariant torus a
parameterized curve  of the form $K_v(\theta)=\hat K(\theta,\theta_v^*)$, for
some fixed $\theta_v^*\in[0,1]$.  A calculation shows that
\[
   \varphi_{1/\omega_v}\bigl(K_v(\theta))
   =
   K_v\bigl(\theta+\nu_v \bigr)
   ,
\]
which is Eq.~\eqref{eq:invtrnm2} with $T= 1/\omega_v$, $\omega= \nu_v$, with
$\nu_v$ defined as in \eqref{eq:defnupnuvtori}. Close to vertical Lyapunov
p.o., an invariant curve of the linearized flow around the Lyapunov
p.o.~satisfies the linearized version of the previous equation and thus
provides an approximate solution, in order to obtain a first torus and start
continuation.
   
When globalizing the invariant curve $K_v$ (for the time-$T$ flow) to the
invariant torus (for the vector field) via Eq.~\eqref{hK} we get 
\[
   \hat K_v(\hat \theta) = \hat K_v(\theta_1,\theta_2)= 
   \varphi_{\theta_2/\omega_v} (K_v(\theta_1-\theta_2 \nu_v)) = 
   \hat K(\hat A_v \hat \theta + \hat \theta_v^*), 
\]
where 
\begin{equation}
\label{def:hAvertical}
     \hat A_v= \begin{pmatrix} 1 & 1 \\ 0 &  1 \end{pmatrix},\ 
     \hat \theta_v^* = \begin{pmatrix} 0 \\ \theta_v^* \end{pmatrix}.
\end{equation}
This is a reparameterization of the invariant torus $\mathcal K$, for which
the frequencies are $\hat \omega_v= {\hat A_v}^{-1} \hat \omega =
(\omega_p-\omega_v, \omega_v)$.

The geometrical observables provided by the Calabi invariants of the two
parameterizations of the torus $\mathcal K$ are related by the identities
\[
     C(\hat K_v) = \begin{pmatrix} 1 & 0 \\ 1 & 1 \end{pmatrix} C(\hat K),
     \ C(\hat K) = \begin{pmatrix} 1 & 0 \\ -1 & 1 \end{pmatrix} C(\hat K_v),
\]
see Remark~\ref{rm:calabi-chfreqs}, and, hence, $C_1(\hat K)= C_1(\hat K_v)=
C(K_v)$ (which also follows from the definition of $K_v$) and $C_2(\hat K)=
-C_1(\hat K_v) + C_2(\hat K_v)$. Notice that $C_2(\hat K_v)=
C(\varphi_{\theta/\omega_v}(K_v(-\theta\nu_v)))= C(\varphi_{T
\theta}(K_v(-\theta\omega)))$. In summary, we can compute the Calabi
invariants of $\hat K$ from $\hat K_v$.
   
We will denote as {\em planar generator} of the invariant torus a
parameterized curve  of the form $K_p(\theta)= \hat K(\theta_p^*,-\theta)$,
for some fixed $\theta_p^*\in[0,1]$. A calculation shows then that
\[
   \varphi_{1/\omega_p}\bigl(K_p(\theta))
   = K_p\bigl(\theta+\nu_p\bigr)
   ,
\]
wich is Eq.~\eqref{eq:invtrnm2} with $T=1/\omega_p$, $\omega= \nu_p$, with
$\nu_p$ defined as in \eqref{eq:defnupnuvtori}. Close to a planar Lyapunov
orbit, an invariant curve of the linearized flow around the Lyapunov
p.o.~satisfies the linearized version of the previous equation, and thus
provides an approximate solution, in order to obtain a first torus and start
continuation.
   
When globalzing the invariant curve $K_p$ (for the time-$T$ flow) to the
invariant torus (for the vector field) via Eq.~\eqref{hK} we get in this case
\[
   \hat K_p(\hat \theta) = \hat K_p(\theta_1,\theta_2)= 
   \varphi_{\theta_2/\omega_p} (K_p(\theta_1-\theta_2 \nu_p)) = 
   \hat K(\hat A_p \hat \theta + \hat \theta^*_p), 
\]
where 
\begin{equation}
\label{def:hAplanar}
     \hat A_p= \begin{pmatrix} 0 & 1 \\ -1 &  1 \end{pmatrix},\ 
     \hat \theta_p^* = \begin{pmatrix}  \theta_p^* \\ 0 \end{pmatrix}.
\end{equation}
This is another reparameterization of the invariant torus $\mathcal K$, for
which the frequencies are 
$\hat \omega_p= {\hat A_p}^{-1} \hat \omega = (\omega_p-\omega_v, \omega_p)$.

The Calabi invariants of the two parameterizations of the torus $\mathcal K$
are related by the identities
\[
     C(\hat K_p) = \begin{pmatrix} 0 & -1 \\ 1 & 1 \end{pmatrix} C(\hat K),
     \ C(\hat K) = \begin{pmatrix} 1 & 1 \\ -1 & 0 \end{pmatrix} C(\hat K_p),
\]
see Remark~\ref{rm:calabi-chfreqs}, and, hence, $C_1(\hat K)= C_1(\hat K_p)
+ C_2(\hat K_p)$ and $C_2(\hat K)= -C_1(\hat K_p)= -C(K_p)$ (as it follows
from the definition of $K_p$). Notice that $C_2(\hat K_p)=
C(\varphi_{\theta/\omega_p}(K_p(-\theta\nu_p)))= C(\varphi_{T
\theta}(K_p(-\theta\omega)))$. In summary, we can compute the Calabi
invariants of $\hat K$ from $\hat K_p$.

\subsection{The numerical explorations}

We have performed two numerical explorations of the Lissajous family.  In the
first one we compute tori with $\rho>\nu_v^{h_0}$, whereas in the second one
we compute tori with $\rho<\nu_v^{h_0}$. Since the tori with
$\rho>\nu_v^{h_0}$ are included in the ``$<100$'' sub-region of
Fig.~\ref{fig:enrho}, in the first exploration we are able to compare the
performance of this paper's parameterization procedure against the one of the
large matrix approach of \cite{GomezM01}.

In our first exploration, we have chosen $74$ values of $\rho>\nu_v^{h_0}$,
equally spaced between $0.02944$ and $0.08754$, and ``nobilized''\footnote{A
noble number is one whose continued fraction expansion coefficients are equal
to one from a position on.} with an absolute tolerance of $\sci{1.6}{-4}$.
For each of these values of $\rho$, we have performed continuation of
invariant curves given by vertical generators $K_v$, for constant
$\omega=\rho$ and increasing $T$, starting from a curve on a narrow torus
around a vertical Lyapunov p.o.~and finishing by collapsing to another
vertical Lyapunov p.o.~of a higher energy. We have also simultaneously
computed their invariant bundles. To do so, we have performed continuation
with respect to $T$ using Algorithm~\ref{algo:cntstep}, making predictions
through Algorithm~\ref{algo:cntT}, and refining each prediction through
Algorithm~\ref{algo:nwttr} (isochronous case). The parameters used in
Algorithm~\ref{algo:cntstep} have been: $m= 4$, $\varepsilon=10^{-7}$,
$\varepsilon^W=10^{-5}$, $\varepsilon_1=10^{-8}$, $\varepsilon_2=10^{-12}$,
$n_{\it des}=4$, $n_\alpha=5$. The stopping criterion has been that, when
approaching the final vertical p.o., $|C_1(\hat K)| = |C(K_v)| <0.001$ (see
Section~\ref{ssection:vertical and planar}).  We have repeated this first
exploration using the large-matrix approach of \cite{GomezM01}, selecting the
parameters accordingly for a fair comparison.

\begin{figure}[htbp]
\includegraphics{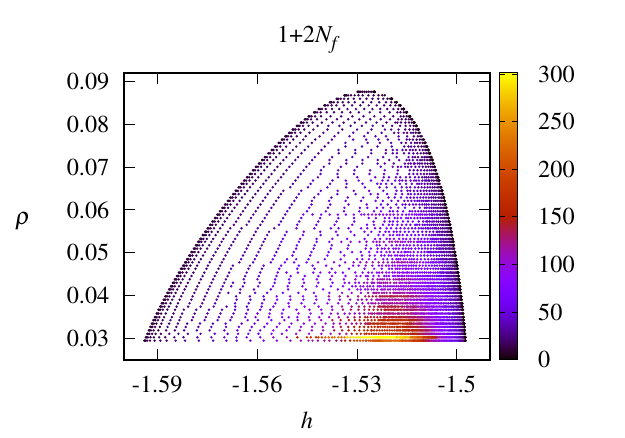}\includegraphics{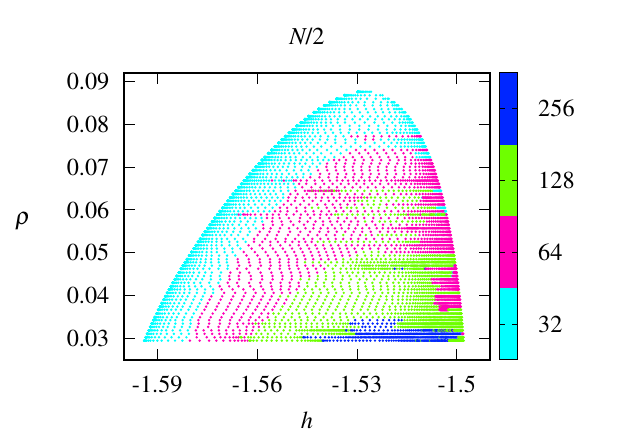}
\caption{\label{fig:rsftpevstrbdfl}Tori with $\rho>\nu_v^{h_0}$ of the
Lissajous family of the Earth-Moon RTBP around $L_1$, computed with the large
matrix procedure (left) and the parameterization one (right).}
\end{figure}

The results of this numerical exploration are shown in
Fig.~\ref{fig:rsftpevstrbdfl}. The left plot corresponds to the large matrix
approach, whereas the right plot corresponds to this paper's parameterization
one. Both plots show the total number of Fourier coefficients used in the
computation. For the left plot this is $1+2N_f$. For the right plot, it is
considered to be $N/2$, because of the DFT queue cleaning strategy mentioned
in Section~\ref{sec:compimpl}.  The total computing time\footnote{The
computing times given will always be qualified as ``total'', meaning the sum
of all the times used by all the threads. The actual wall-clock time is
roughly this time divided by the number of cores (16 in our case). This rule
is not followed exactly because of uneven load balancing: the continuation of
some constant-$\omega$ families takes longer than others.} of the large-matrix
approach is 68308 seconds, of which 22793 are spent in the computation of the
tori, whereas the rest are used in the computation of the stable and unstable
bundles through a slight modification of the method presented in
\cite{2001Jo}. The total computing time of the parameterization approach of
this paper (that includes tori and bundles) is 5992 seconds.  The total number
of tori computed is 4141 with the large-matrix method vs.~7008 with the
parameterization one.  This is due to the fact that, as can be appreciated in
Fig.~\ref{fig:rsftpevstrbdfl}, the continuation strategy of the large-matrix
procedure is able to use larger step sizes for tori with small values of
$N_f$. 

In a second exploration, we have computed invariant tori that are born from
planar Lyapunov orbits that have rotation number  $\rho<\nu_v^{h_0}$.  Namely,
we have chosen 31 values of $\rho$, equally spaced between 0.02785 and
0.00317, nobilized with an absolute tolerance of $\sci{1.6}{-4}$ and for each
of the values we have performed continuation of invariant curves given by
planar generators $K_p$, for constant $\omega=1-1/(1+\rho)$ and increasing
$T$, starting from a narrow torus around a planar Lyapunov p.o.  The
parameters used in the continuation algorithms have been the same as before.
The continuations have stopped either by reaching the computational limit or,
as in the first exploration, when $|C_1(\hat K)| = |C(K_v)| <0.001$ (see
Section~\ref{ssection:vertical and planar}).

The motivation of this second exploration is twofold. On the one hand, to
perform a ``stress test'' of our procedure by exploring phase space beyond the
computations of \cite{GomezM01}. On the other, to relate the behavior of
dynamical and geometrical observables to the destruction of invariant tori
(see the next section).  In this exploration, we have always achieved
convergence of Newton's method provided that the continuation step is small
enough and the number of points $N$ large enough.  The computational limits on
these two quantities, chosen as $10^{-5}$ and $8192$, respectively, have been
set in order to obtain reasonable run time and storage requirements in a
single workstation. Notice that working with such a large number of Fourier
coefficients is unfeasible with the large matrix approach, since it would
require the solution of non-linear sytems of equations with a size of nearly
$100000\times100000$.  A limit, equal to 10000, has also been put for the
maximum number of tori computed in each constant $\rho$ family of this second
exploration.  With all these limits, this exploration has run for a total time
of 27.0816 days and has generated a total of $130574$ tori that, each
compressed as a bz2 file, take up $127.51\,{\rm GiB}$ of disk
space.\footnote{Recall that the ``wall-clock'' time is roughly the total time
divided by 16. On the other hand, several strategies, that we have not pursued
here, can be used to reduce greatly these storage requirements, like using
binary files with single precision floating-point numbers, and not storing all
the tori but a grid of them fine enough in order to recover tori not in the
grid by interpolation (see e.g.~\cite{2012aMoBaGoOll}).}

\begin{figure}
\includegraphics{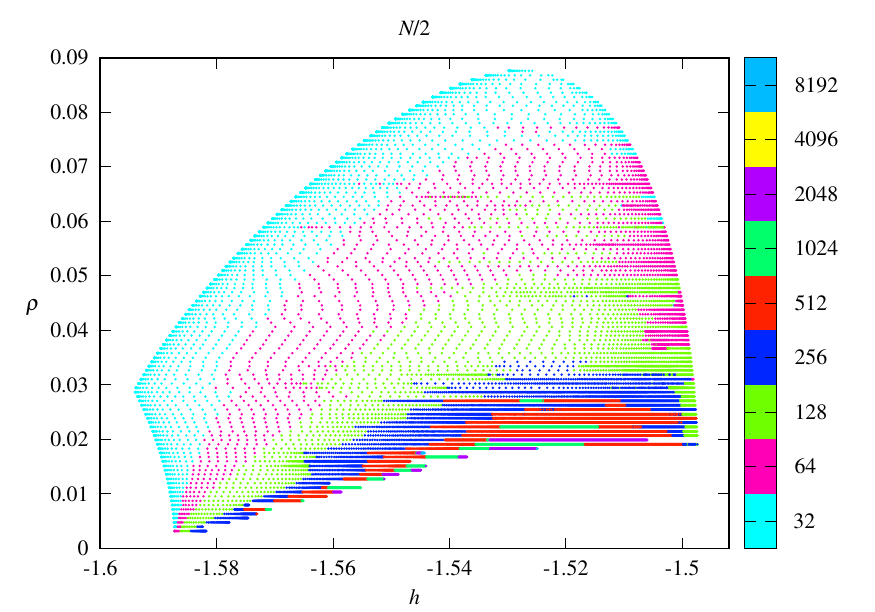}
\caption{\label{fig:totnp-red}Number of Fourier coefficients of all the tori
and invariant bundles computed of the Lissajous family of the Earth-Moon RTBP
around $L_1$, using this paper's parameterization strategy.}
\end{figure}

The results of this second exploration are shown in Fig.~\ref{fig:totnp-red},
that is analogous to Fig.~\ref{fig:rsftpevstrbdfl} right but including both
explorations. Many of the figures that follow will refer to the two
explorations as a whole.

\subsection{Dynamical and geometric observables}

In this section we will describe the behaviour of different dynamical and
geometric observables of the invariant tori of the Lissajous family that are
obtained during their computation. These observables provide insight on both
numerical and dynamical vicissitudes faced by the method. We will also display
the evolution of three constant $\rho$ families through the 3D representation of
some of their tori in configuration space. The results presented here
complement those in \cite{GomezM01} that,  by using iso-energetic Poincar\'e
sections, provides a detailed account of the evolution of the Lissajous family
of invariant tori, together with its interaction with other families of
invariant tori and periodic orbits.

For a parameterization $K= K_0:\T\to \R^6$ of the invariant curve of the time
$T$-flow with rotation number $\omega$, obtained with multiple shooting with
$m= 4$ steps, that is a generator of the 2D invariant torus given by
\eqref{hK}, the dynamical observables we consider here are:
\begin{itemize}
\item
   the flying time $T$ of the generator (notice that we get frequencies for
   the 2D-torus by the formula $\hat \omega= \tfrac{1}{T}(\omega,1)$);
\item
   the unstable Floquet multiplier $\Lambda^u= \lambda^{-m}$, that  provides
   information about hyperbolity properties of the  invariant curve, and from
   which one can obtain the Floquet exponent of the  2D-torus by $\chi=
   \tfrac{1}{T} \log \Lambda^u$.
\end{itemize}
Notice we select $m$ so that $\lambda^{-1}$ is not too big, in order to
mitigate numerical unstabilities. In our computations, $\lambda^{-1}$ runs in
the interval $[4,8]$.

The geometric observables we consider are: 
\begin{itemize}
\item
   The Calabi invariants $C_1(\hat K), C_2(\hat K)$ of the parameterization
   $\hat K(\theta_1,\theta_2)$ of the torus  with natural frequencies
   $\hat\omega= (\omega_p,\omega_v)$.  These invariants give insight on the
   size of the generators in area units. Section \ref{ssection:vertical and
   planar} provides formulae for their computation from $C(\hat K_p)$, $C(\hat
   K_v)$.
\item 
   The (minimum) distances between several pairings of bundles on the
   generator curve:  $T{\mathcal K}$, the tangent bundle of the generator
   curve (generated by $K'$); $X$, the bundle generated by the vector field on
   the curve; the stable and the unstable bundles $E^s$ and $E^u$,
   respectively; the central bundle, $E^c$, that has rank 4 and contains the
   tangent bundle to the generator, the vector field and, hence, the tangent
   bundle to the 2D torus. The distances we consider are:  $d(T{\mathcal
   K},X)$, to measure the transversality of the flow to the generator, and
   $d(E^s,E^u)$, $d(E^s,E^c)$, and $d(E^u,E^c)$, to measure the quality of
   hyperbolicity geometrical properties.
\end{itemize}
The bundles are generated by selected columns of the matrix map $P= P_0:\T\to
\R^{6\times 6}$, that we write as
\[
     P(\theta)= 
     \begin{pmatrix}
        K'(\theta) & X_H(K(\theta)) & W^s(\theta) & N^1(\theta) &
        N^2(\theta) & W^u(\theta)
     \end{pmatrix}
\]
for reference. Hence, at a point $K(\theta)$ of the invariant curve, the fiber
of the stable bundle $E^s$ is generated by $W^s(\theta)$, the fiber of the
unstable bundle $E^u$ is generated by $W^u(\theta)$, and fiber of the center
bundle $E^c$ is generated by $K'(\theta)$, $X_H(K(\theta))$, $N^1(\theta)$,
$N^2(\theta)$. Notice that the tangent bundle to the generator is generated by
$K'(\theta)$, and the tangent bundle of the 2D torus is generated by
$K'(\theta)$, $X_H(K(\theta))$.  There are several ways of defining distances
or angles between vector subspaces of a given normed vector space.  Here, the
vector space is $\R^6$, with the norm induced by the standard scalar product,
and the distance we consider between a vector subspace $E_1$ of dimension $1$
and another vector subspace $E_2$  is the length of the projection onto
$E_2^\perp$ of a unit vector in $E_1$ (the angle between $E_1$ and $E_2$ is
the $\arcsin$ of this lenght). Finally, we define the distance between two
bundles as the minimum distance between corresponding fibres of the bundles.

\begin{remark}
\label{remark:hyperbolicity breakdown}
We emphasize that the quality of hyperbolicity properties of the invariant
torus are not only given by the Floquet multipliers in the stable and unstable
directions, that have to be away from 1, but also by the positivity of the
angles between the stable, unstable and center directions. There are
mechanisms of breakdown of invariant tori that involve the degeneration of
some of these angles, that go to zero, while the stable and unstable Floquet
multipliers remain far from 1. See
\cite{CallejaF12,CanadellH14,CanadellH17b,FiguerasH15,HaroL06c,HaroL07}.
\end{remark}

\begin{remark}
Reversibility properties of the RTBP imply that stable and unstable bundles
can be obtained from each other using reversors, and that they have same
angles with the center manifold of the torus. These properties could also be
used to reduce the cost of the algorithms presented here (reducing, for
instace, the cost of generating the frame). We prefer not doing so for the
sake of generality. 
\end{remark}

\begin{figure}
\begin{center}
\includegraphics{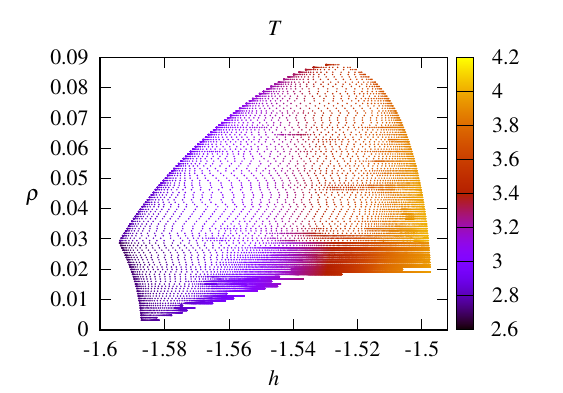}
\includegraphics{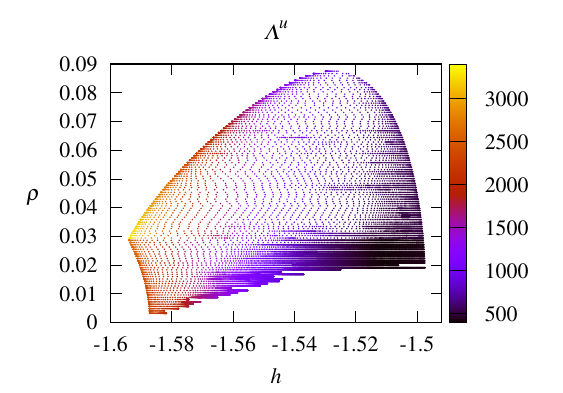}

\includegraphics{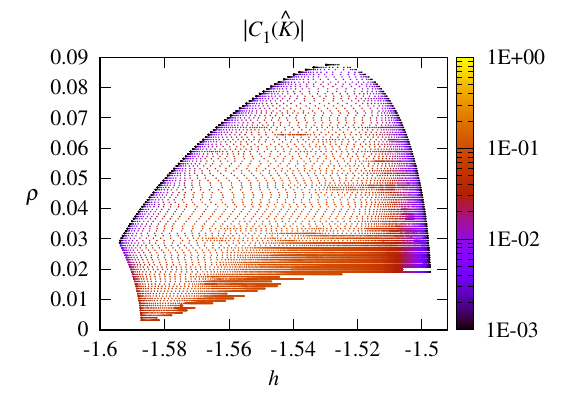}
\includegraphics{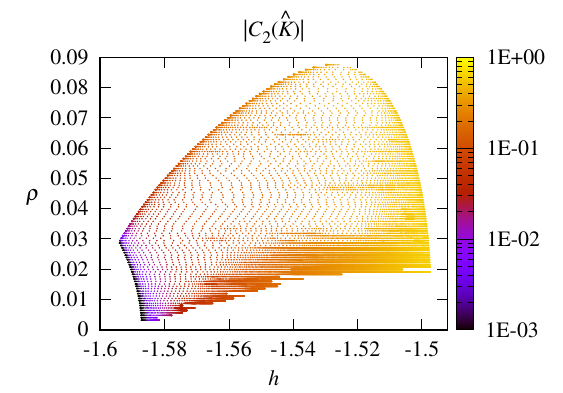}

\includegraphics{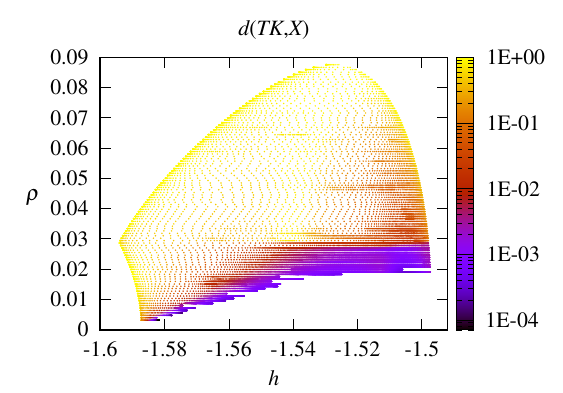}
\includegraphics{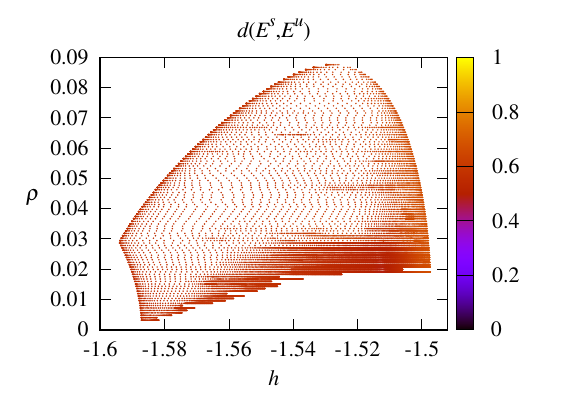}

\includegraphics{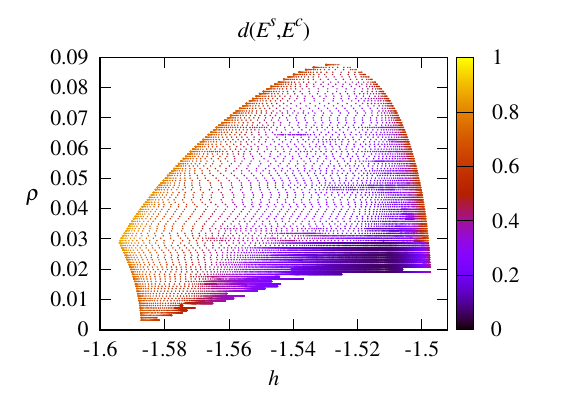}
\includegraphics{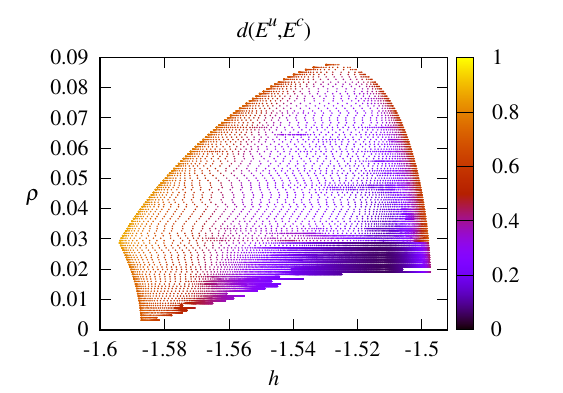}
\end{center}
\caption{\label{fig:geomobs}Dynamical and geometrical observables}
\end{figure}

By monitoring these observables during the continuation we can get insight
about dynamical and geometric properties of the torus (and its invariant
bundles), and detect numerical unstabilities caused by degeneracies of these
properties (such as the hyperbolicity, regularity of the frame, size of the
generator). For instance, we recall that  one stopping criterion is that
$|C_1(\hat K)| = |C(K_v)| < 0.001$, revealing that the torus is approaching a
periodic orbit. We have collected these observables from the two numerical
explorations exposed in the previous section, and the results are summarized
in Figure~\ref{fig:geomobs}.

In this massive computation we observe that:
\begin{itemize}
\item
   The unstable multiplier ranges from $413.205$ to $3344.26$, from which the
   spectral condition of  hyperbolicity of the torus is satisfied;
\item
   The distance between the stable and unstable bundles is bigger than
   0.492489, and the distante between the stable and center bundle, and the
   unstable and center bundle, is bigger than 0.0615721, from which the
   geometrical conditions of hyperbolicy is also satisfied.
\end{itemize} 

The continuations of families of Lissajous tori with smaller rotation numbers
stop because of the computational limit of 10000 tori for each constant-$\rho$
family. But, as we see from the behavior of the observables, this phenomenon
is not apparently due to the fact that hyperbolicity breaks down.  However,
the fact that the step size of the continuation becomes smaller and the number
of Fourier coefficients becomes larger reveals that the torus is losing
regularity (the analyticy strip of the complex domain of the parameterization
of the torus goes to zero), indicating an obstruction for the existence of the
torus and that it is breaking down.  There is another possible mechanism of
breakdown, it is what we call KAM breakdown.  These tori lie on the center
manifold of the $L^1$ point, $W^c(L_1)$, which is a $4D$ symplectic manifold.
So, inside $W^c(L_1)$, these tori are KAM tori, and the basic mechanism of
breakdown is the collision with resonances (the overlap criterion
\cite{Chirikov79}), which can be more geometrically described as the
obstruction produced by homoclinic and heteroclinic webs produced by the
invariant manifolds of unstable periodic orbits inside the center manifold
(the obstruction criterion in \cite{OlveraS87,LlaveO06}). So, it is very
likely that in this case the breakdown is produced by this phenomenon inside
the center manifold. We will come back to this issue later. 

\begin{figure}[htbp]
\vspace*{-.7cm}
\begin{tabular}{cc}
\includegraphics[width=.5\linewidth]{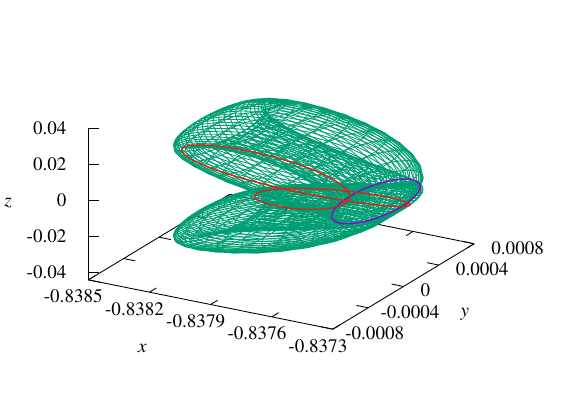} & \includegraphics[width=.5\linewidth]{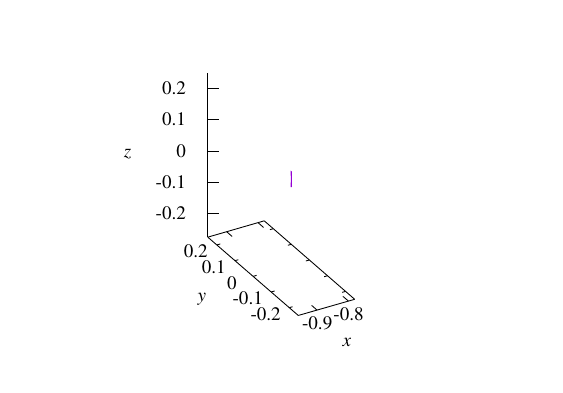} 
\vspace{-1.25cm}\\
\includegraphics[width=.5\linewidth]{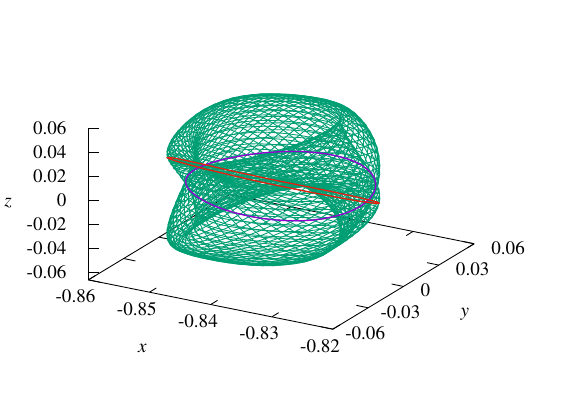} & \includegraphics[width=.5\linewidth]{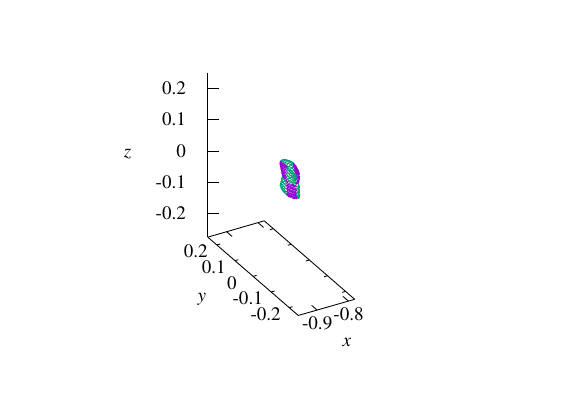} 
\vspace{-1.25cm} \\
\includegraphics[width=.5\linewidth]{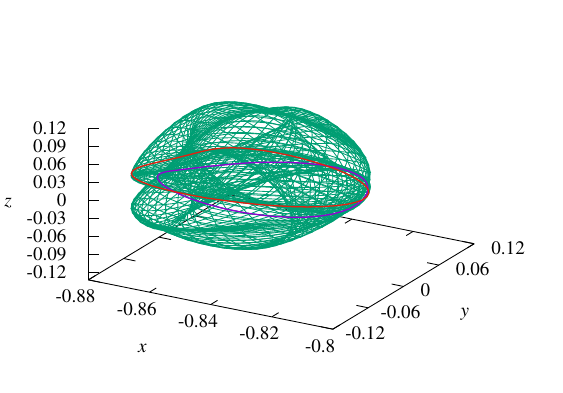} & \includegraphics[width=.5\linewidth]{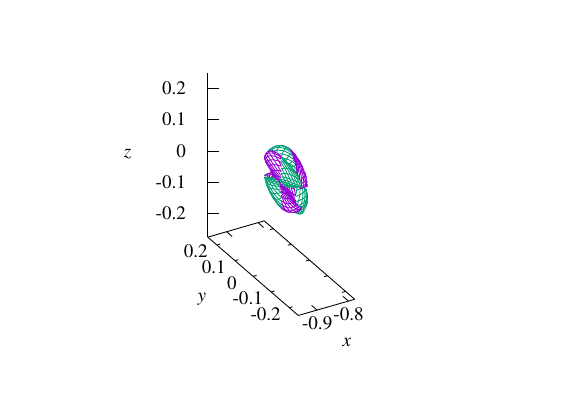} 
\vspace{-1.25cm} \\
\includegraphics[width=.5\linewidth]{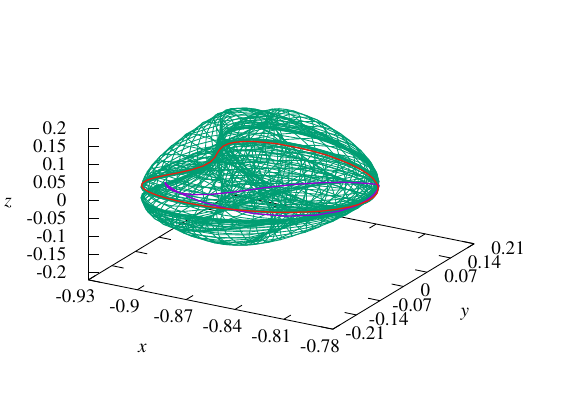} & \includegraphics[width=.5\linewidth]{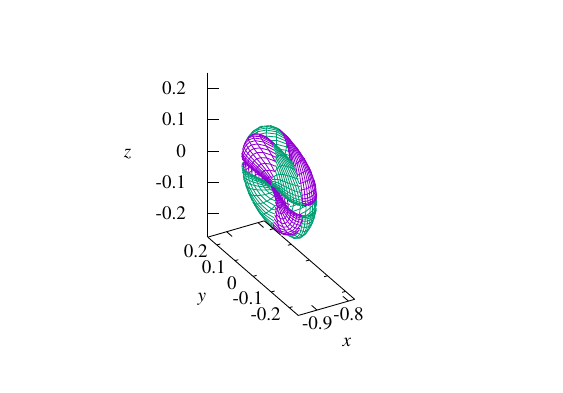} 
\vspace{-1.25cm} \\
\includegraphics[width=.5\linewidth]{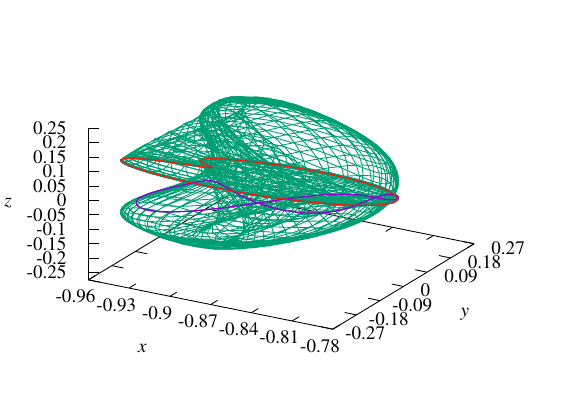} & \includegraphics[width=.5\linewidth]{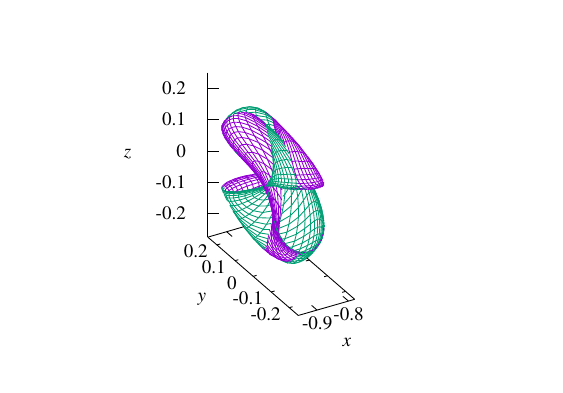} 
\vspace{-1.25cm} \\
\includegraphics[width=.5\linewidth]{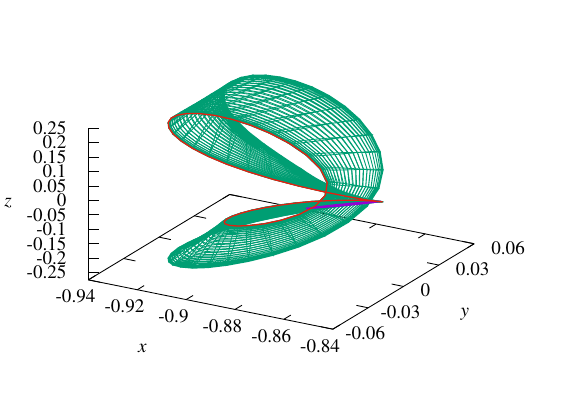} & \includegraphics[width=.5\linewidth]{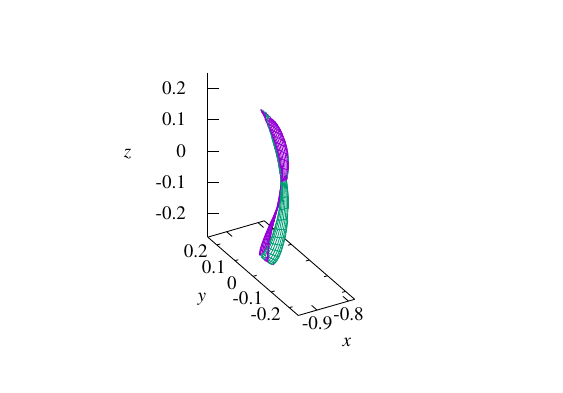} 
\end{tabular}
\caption{\label{fig:cnt008-70-wrgen}
Two views of a sample of invariant tori  with $\rho=0.031865$. 
Letf: Parameterized surfaces, including generators
(Blue: invariant curve $\{\hat K(\theta,0)\}_{\theta\in\T}$. Red: invariant
curve $\{\hat K(0,\theta)\}_{\theta\in\T}$). Right: with the same scale on all axes, 
and with opaque surfaces.}
\end{figure}

\begin{figure}[htbp]
\vspace*{-.7cm}
\begin{tabular}{cc}
\includegraphics[width=.5\linewidth]{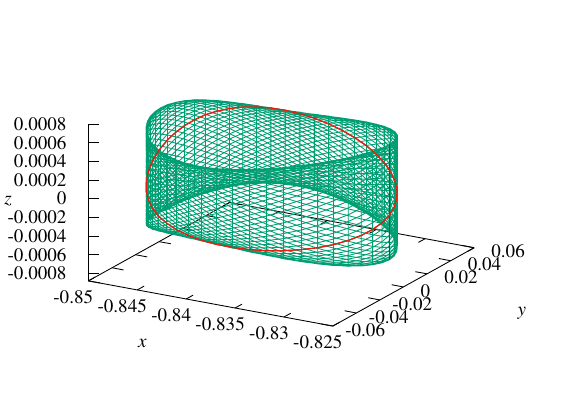} &
\includegraphics[width=.5\linewidth]{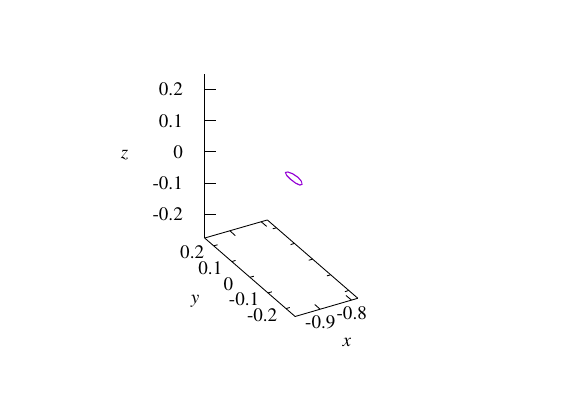} 
\vspace{-1.25cm}\\
\includegraphics[width=.5\linewidth]{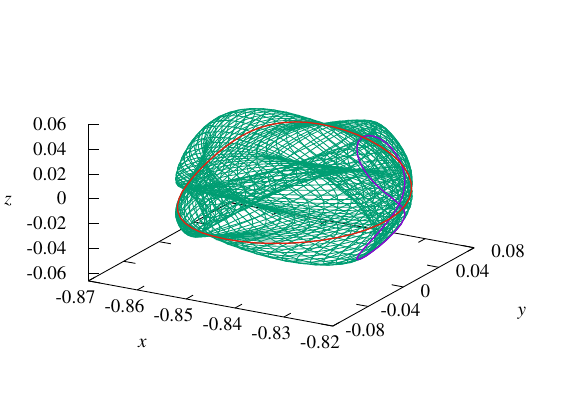} &
\includegraphics[width=.5\linewidth]{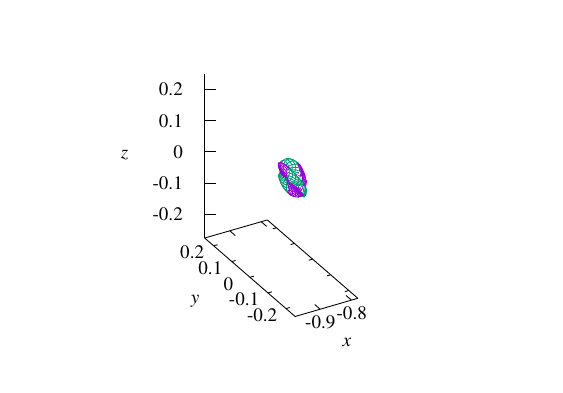} 
\vspace{-1.25cm} \\
\includegraphics[width=.5\linewidth]{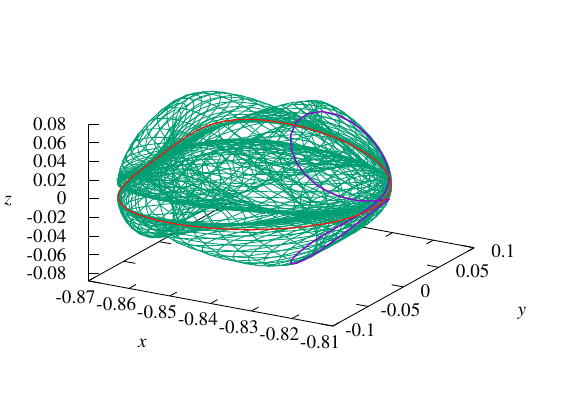} &
\includegraphics[width=.5\linewidth]{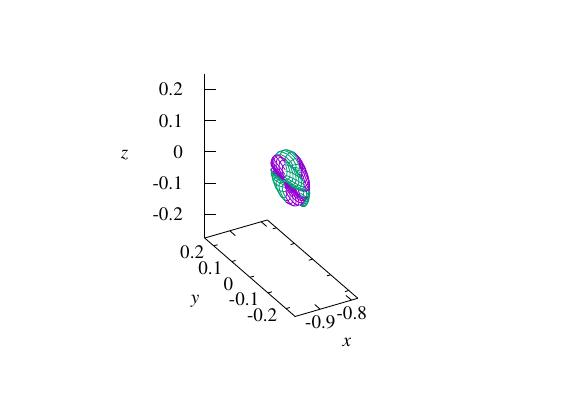} 
\vspace{-1.25cm} \\
\includegraphics[width=.5\linewidth]{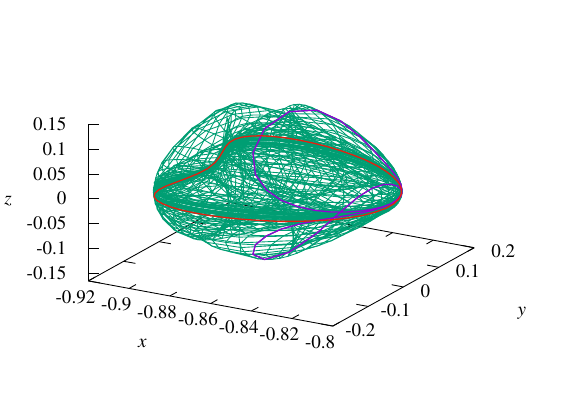} &
\includegraphics[width=.5\linewidth]{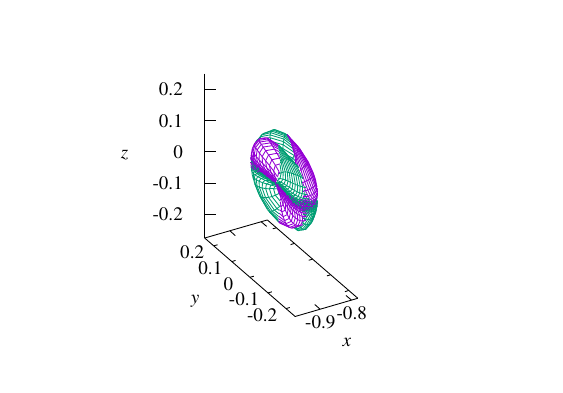} 
\vspace{-1.25cm} \\
\includegraphics[width=.5\linewidth]{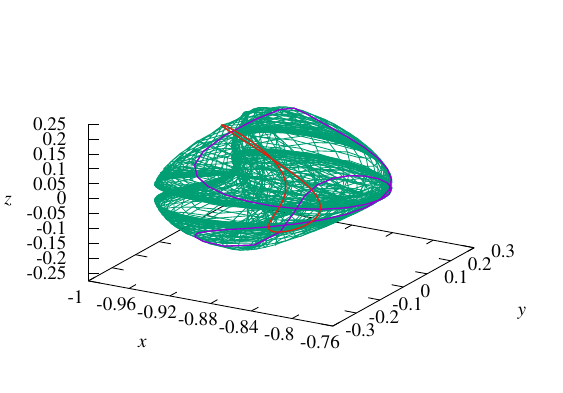} &
\includegraphics[width=.5\linewidth]{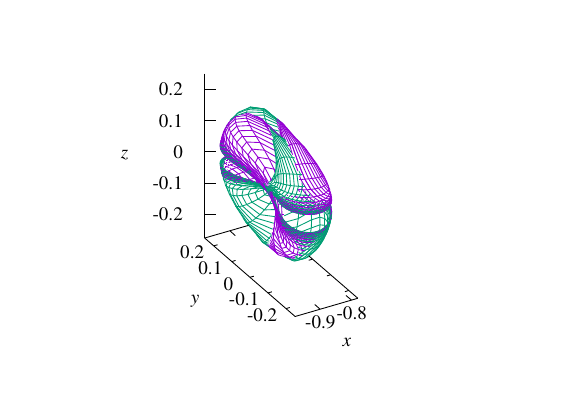} 
\vspace{-1.25cm} \\
\includegraphics[width=.5\linewidth]{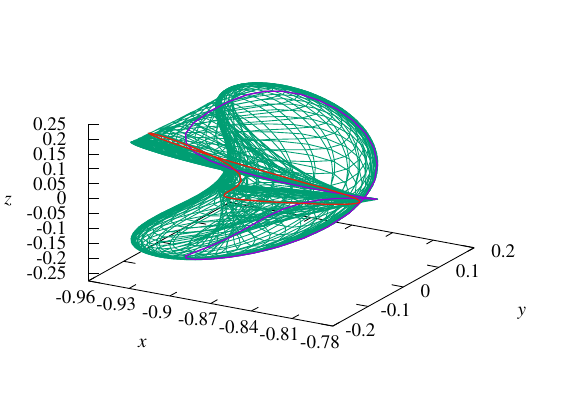} &
\includegraphics[width=.5\linewidth]{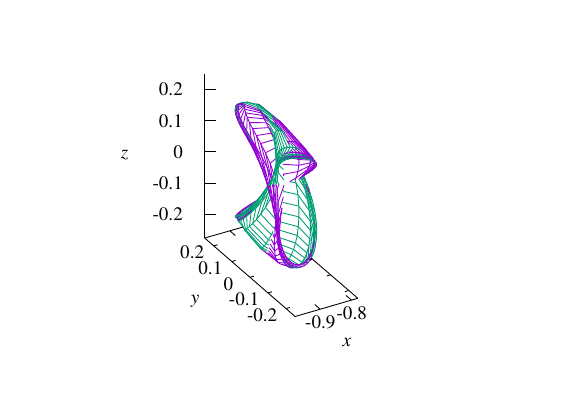} 
\end{tabular}
\caption{\label{fig:baix005-13-wrgen}Two views of a sample of invariant tori
with $\rho=0.019091$. Left: parameterized surfaces, including generators (Blue: invariant curve $\{\hat K(\theta,0)\}_{\theta\in\T}$. Red: invariant
curve $\{\hat K(0,\theta)\}_{\theta\in\T}$). Right: with the same scale on all
axes, and with opaque surfaces.}
\end{figure}

\begin{figure}[htbp]
\vspace*{-.7cm}
\begin{tabular}{cc}
\includegraphics[width=.5\linewidth]{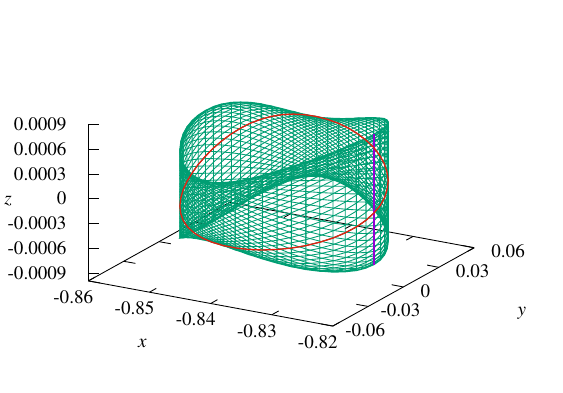} &
\includegraphics[width=.5\linewidth]{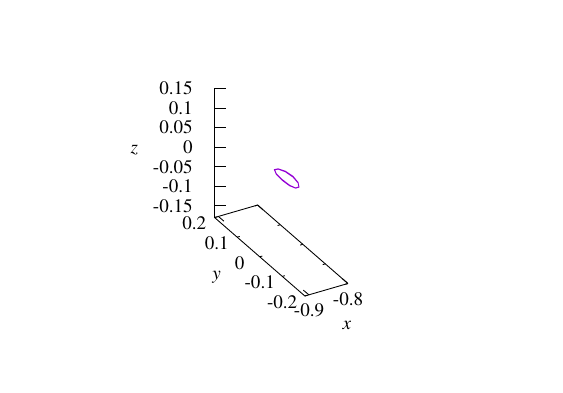} 
\vspace{-1.25cm}\\
\includegraphics[width=.5\linewidth]{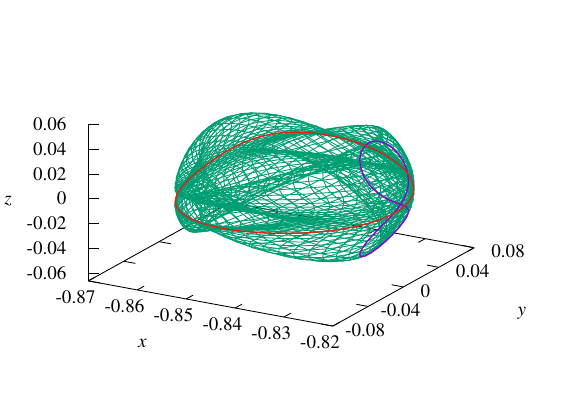} &
\includegraphics[width=.5\linewidth]{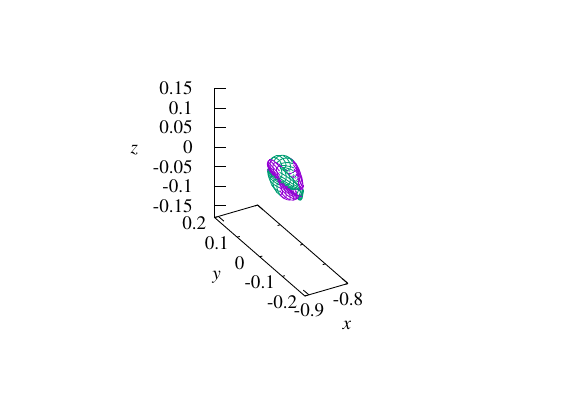} 
\vspace{-1.25cm} \\
\includegraphics[width=.5\linewidth]{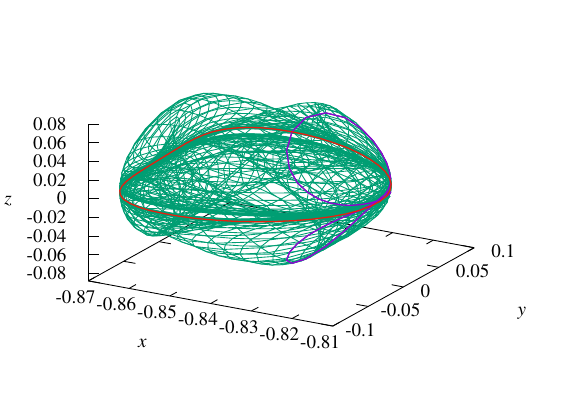} &
\includegraphics[width=.5\linewidth]{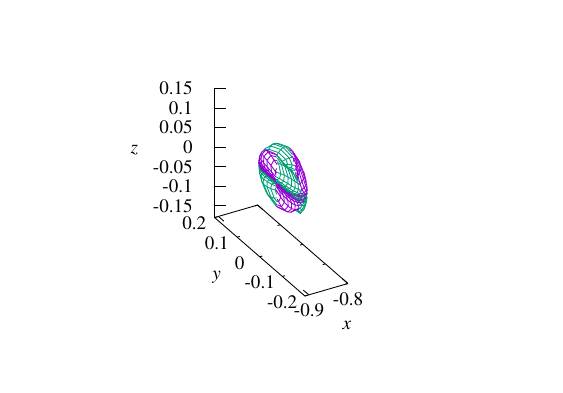} 
\vspace{-1.25cm} \\
\includegraphics[width=.5\linewidth]{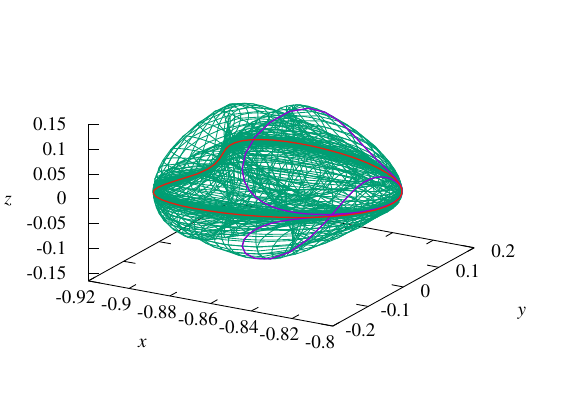} &
\includegraphics[width=.5\linewidth]{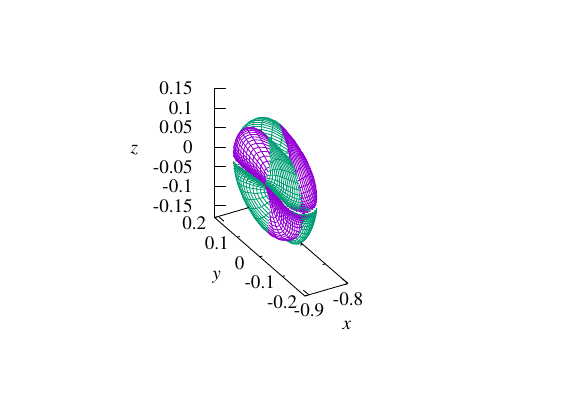} 
\end{tabular}
\caption{\label{fig:baix005-20-wrgen}Two views of a sample of invariant tori
with $\rho=0.013584$. Left: parameterized surfaces, including generators (Blue: invariant curve $\{\hat K(\theta,0)\}_{\theta\in\T}$. Red: invariant
curve $\{\hat K(0,\theta)\}_{\theta\in\T}$). Right: with the same scale on all
axes, and with opaque surfaces.}
\end{figure}

In the following, we will particularize the results for three families of
Lissajous tori, with rotation numbers $\rho=0.031865$, $\rho=0.019091$ and
$\rho=0.013584$. Figs.~\ref{fig:cnt008-70-wrgen}, \ref{fig:baix005-13-wrgen}
and \ref{fig:baix005-20-wrgen} show several samples of tori of these three
families, projected on the configuration space in different forms and views: 
\begin{itemize}
\item[(left)]
   as grids on the parameterized surfaces $\{\hat
   K(\theta_1,\theta_2)\}_{(\theta_1,\theta_2)\in\T^2}$ (see Eq.~\eqref{hK}),
   including two generators of the homothopy group of the torus, given by
   $\{\hat K(\theta,0)\}_{\theta\in\T}$,  in blue, and $\{\hat
   K(0,\theta)\}_{\theta\in\T}$, in red; 
\item[(right)]
   as opaque surfaces, with the same scale and range on all axes, with colors
   corresponding to the different sides of the surface, revealing
   self-intersections of the projections of tori on configuration space.
\end{itemize}
Actually, instead of Eq.~\eqref{hK}, the expression
\[
   \hat K(\theta_1,\theta_2)
   =
   \varphi_{(\theta_2-\frac jm)T}
      \Bigl(
         K_j\bigl(
            \theta_1-(\theta_2-\frac jm)\omega
         \bigr)
      \Bigr)
      ,
\]
with $j=[m\theta_2]$, has been used, in order to take advantage of multiple
shooting. 
In addition to these views, for each of these families we have plotted in Figs.~\ref{fig:cnt008-70-obs}, 
\ref{fig:baix005-13-obs} and \ref{fig:baix005-20-obs} the dynamical and geometrical observables as functions 
of the energy $h$.   We will describe our findings below.

The results for the family $\rho= 0.031865$ are summarized in
Figs.~\ref{fig:cnt008-70-wrgen} and \ref{fig:cnt008-70-obs}.  From
Fig.~\ref{fig:cnt008-70-wrgen} we appreciate how the family begins with a
small torus around a vertical Lyapunov p.o., that grows up to approximately
the size of the planar Lyapunov orbit of the same energy, and then starts
bending until it ``is about to close''.  Then it opens again and shrinks until
it collapses to a vertical Lyapunov p.o.~of higher energy. All this is done
while increasing in size, since energy also increases.
Fig.~\ref{fig:cnt008-70-obs} displays the observables for this family. Since
the family is born in a vertical Lyapunov p.o. and dies in another vertical
Lyapunov p.o. of a higher energy, the Calabi invariant of the generator starts
being 0 and finishes being 0. Notice that, in both cases, close to the p.o.
the Calabi invariant goes to zero assymptotically as a linear function of the
difference of the energy with the one of the  p.o. 

The results for the family $\rho= 0.019091$ are summarized in
Figs.~\ref{fig:baix005-13-wrgen} and \ref{fig:baix005-13-obs}.
Fig.~\ref{fig:baix005-13-wrgen} shows a sample tori starting in a planar
Lyapunov p.o. and, hence, their homothopy group generators are exchanged with
respect to Fig.~\ref{fig:cnt008-70-wrgen}.  The evolution with energy is
similar to the one of Fig.~\ref{fig:baix005-13-wrgen}, with three main
differences: the torus ``seems to close after bending'' (a zoom of the fourth
torus reveals that it does not actually close), there is a more important
accumulation of wireframe lines at the ``boundary that closes and opens'', and
the family starts from a planar Lyapunov p.o. instead of a vertical one.  This
is in fact the reason the Calabi invariant of the computed generators starts
being 0 (and, again, with an asymptotic linear behaviour) and increases till
the end of the continuation (converging to the Calabi invariant of the
vertical Lyapunov p.o.), as it is observed in Fig.~\ref{fig:baix005-13-obs}.

\begin{figure}[htbp]
\begin{tabular}{c}
\includegraphics{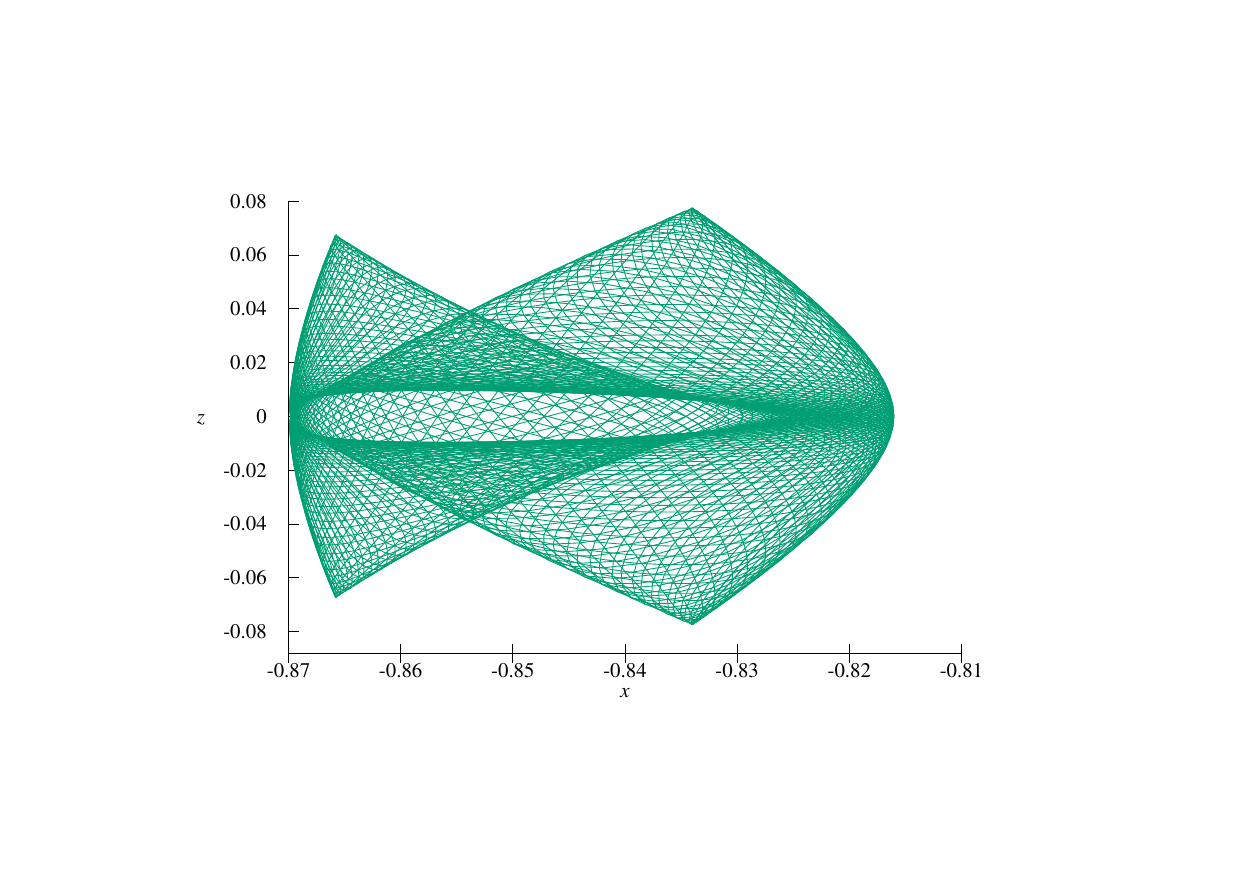}
\vspace{-2cm}	\\
\includegraphics{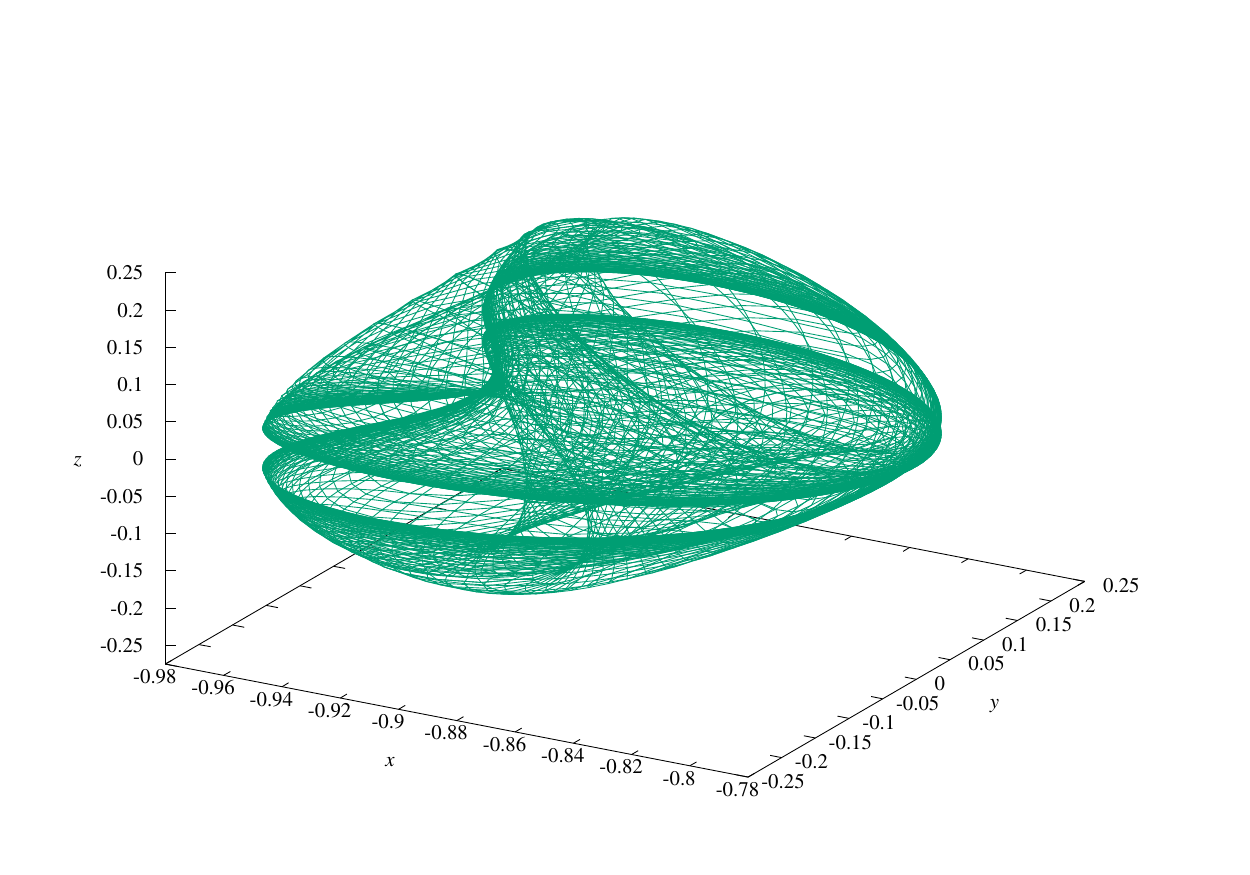}
\end{tabular}
\caption{\label{magnification}  Magnifications of projections of the third and fifth tori of 
Fig.~\ref{fig:baix005-13-wrgen}, corresponding to $\rho=0.019091$.}
\end{figure}

In the iso-energetic Poincar\'e section plots of \cite{GomezM01} (and, for
lower energy levels, in the ones of previous references like
\cite{ESA3,1999JoMa}), it is numerically seen that, from energy $h_B$ to
$h_C$, a double homoclinic connection inside $W^c(L_1)$ of the planar Lyapunov
family of p.o.~acts as separatrix from the Lissajous family and the quasi-Halo
family of tori. From energy $h_C$ to $h_D$, this role is taken by heteroclinic
connections (also inside $W^c(L_1)$) between the vertically symmetric families
of p.o.~that are born at the second 1:1 bifurcation of the planar Lyapunov
p.o. Plots of orbits of these last families of p.o.~can be found in Fig.~8 of
\cite{GomezM01} (they are known as ``axial'' by other autors,
e.g.~\cite{2003aDiDoPa}).  In Fig.~\ref{fig:baix005-13-wrgen}, it is observed
how the tori approach these connections. This fact is better appreciated in
Fig.~\ref{magnification}, that shows magnifications of projections of the
third and fifth tori of Fig.~\ref{fig:baix005-13-wrgen}.  The view of the
third torus has been chosen in order to stress the fact that the torus
represented approaches two different, vertically simmetric quasi-Halo tori, as
can be inferred from the Poincar\'e representations of the center manifold in
references \cite{GomezM01,ESA3,1999JoMa}.  When approaching these connections,
dynamics becomes slow and, since the parameterization $\hat K$ is tied to the
dynamics through the invariance equation \eqref{eq:invtrnm0}, it produces the
accumulation of wireframe lines of Fig.~\ref{fig:baix005-13-wrgen} and the
small values of $d(T{\mathcal K},X)$ of Fig.~\ref{fig:geomobs}. We believe
this kind of stiffness to be responsible for the drastic reduction of step
length of the second exploration and, to a lesser extent, of the first one. As
it has been commented, these connections are responsible for the destruction
of the families of invariant tori in the center manifold.

The phenomenon of breakdown is illustrated with the family with
$\rho=0.013584$. The results for this family are summarized in
Figs.~\ref{fig:baix005-20-wrgen} and \ref{fig:baix005-20-obs}. In
Fig.~\ref{fig:baix005-20-wrgen} we observe that the torus is increasingly
pinched, while  dynamical and geometrical observables in
Fig.~\ref{fig:baix005-20-obs} do not suggest the torus is being destroyed.
However, the fact the continuation step is getting very small and the number
of Fourier coefficients of the approximations is getting larger, related to
the mentioned ``pinching phenomenon", envisages the breakdown of the torus
inside the center manifold. 

\newcommand{\dist}{\mathop{\rm dist}\nolimits}

\begin{figure}
\includegraphics{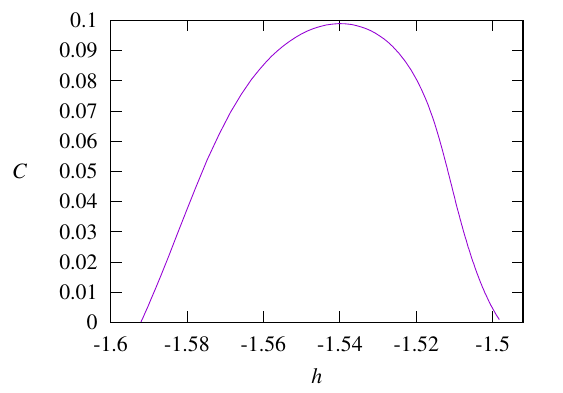}
\includegraphics{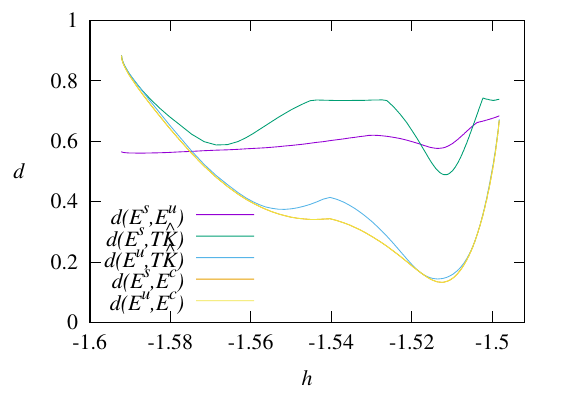}
\caption{\label{fig:cnt008-70-obs} Calabi invariant $C(K_v)=C_1(\hat
K)$, and distances between several pairings of bundles, for the family
$\rho=0.031865$.}
\end{figure}

\begin{figure}
\includegraphics{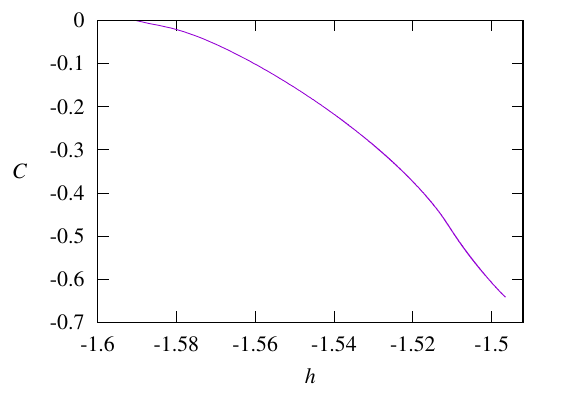}
\includegraphics{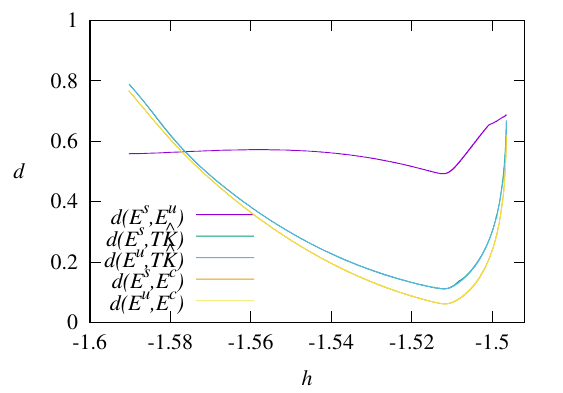}
\caption{\label{fig:baix005-13-obs}Calabi invariant $C(K_p)=-C_2(\hat
K)$, and distances between several pairings of bundles, for the family
$\rho=0.019091$.}
\end{figure}

\begin{figure}
\includegraphics{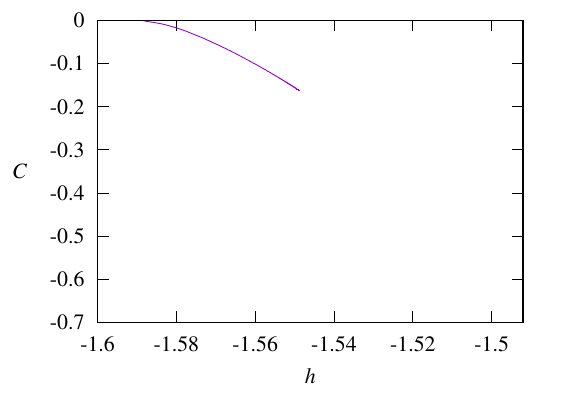}
\includegraphics{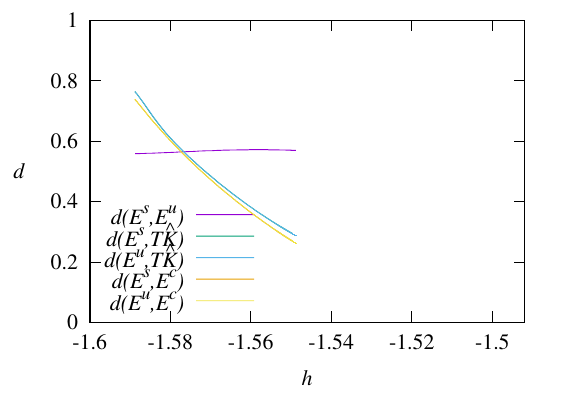}
\caption{\label{fig:baix005-20-obs} Calabi invariant $C(K_p)=-C_2(\hat
K)$, and distances between several pairings of bundles, for the family
$\rho=0.013584$.}
\end{figure}


\section{Conclusion}
\label{sec:conclusion}

We have presented in this paper a very efficient method to compute
invariant tori in Hamiltonian systems. To do so, we have first reduced
the dimensionaly of the objects, by considering invariant tori for flow
maps, and then taken advantage of the geometrical and dynamical
properties of invariant tori in Hamiltonian systems. The method also
provides online information on the linearized dynamics around the torus,
as well as on other geometrical properties. He have focused
our attention on partially hyperbolic invariant tori with rank-one stable
and unstable manifolds, and in this case the method provides not only
parameterizations of the tori but also of the linear approximations
of those manifolds. Tests have been performed for the computation of
invariant tori around libration points of the Circular Restricted
Three Body Problem, for which there is an extensive literature, so any
one could easily compare the performances.

As we have already mentioned, we only present algorithms of computation, based
on Newton's method.  Eventually,  a proof of convergence could be completed
using KAM techniques for obtaining results in a posteriori format.  Although
we have not attempted this, we do provide information about the magical
cancellations appearing in the linearized equations, that are key in the KAM
proofs, and, with much more effort, could be implemented as computer assisted
proofs \cite{FiguerasHL17}.

The increasing complexity of problems and applications has spurred the
research of this paper.  The algorithms presented here are the first ones in a
new generation of algorithms to compute invariant tori and their manifolds in
Hamiltonian systems. We plan to extend the methodologies to more complex
problems in the future, with an eye in the applications. 


\appendix

\section{From Poincar\'e map to time-$T$ map}
\label{ap:Poincare map}

In this section we will see how to obtain, from an invariant torus of a the
Poincar\'e map, an invariant torus of a time-$T$ map. For instance, the
invariant torus of the Poincar\'e map could have been computed from a center
manifold reduction around an equilibrium point at a certain fixed energy level
(see e.g. \cite{Simo98,Jorba99} for normal form methods and \cite{HaroCFLM16}
for direct parameterization methods, applied to the computation of the center
manifold of a colinear fixed point in the RTBP). 

Assume we are given a parameterization $K_P:\T^{d-1}\to \R^{2n}$ of a
$(d-1)$-dimensional torus ${\mathcal K}_P$ inside a $d$-dimensional torus
$\hat{\mathcal K}$, produced  by a Poincar\'e map $P$ associated  to a
transversal section to the vector field $X_H$. We also assume that the
rotation vector of ${\mathcal K}_P$ is $\omega$ (which is assumed to be
Diophantine). That is, we assume for all $\theta\in\T^{d-1}$
\[
	\varphi_{T_P(\theta)}(K_P(\theta))= K_P(\theta+\omega),
\]
where $T_P:\T^{d-1}\to \R$ gives for each $\theta\in\T^{d-1}$ the time for a
point $K_P(\theta)$ to return to the transversal section. The flying time
depends then on the point on the torus.

We want to find a parameterization $K:\T^{d-1}\to \R^{2n}$ of a
$(d-1)$-dimensional torus ${\mathcal K}$  for which the flying time is
constant. To do so, we look for $\tau:\T^{d-1}\to \R$ and $T$ such that
\[
        K(\theta)= \varphi_{\tau(\theta)}(K_P(\theta)),\
        \varphi_{T}(K(\theta))= K(\theta+\omega).
\]
Hence, since
\[
   \varphi_{T}(K(\theta)) = \varphi_{T+\tau(\theta)}(K_P(\theta))
\]
and 
\[
	K(\theta+\omega)= \varphi_{\tau(\theta+\omega)}(K_P(\theta+\omega)) = 
	\varphi_{\tau(\theta+\omega)+T_P(\theta)}(K_P(\theta)),
\]
we impose that
\[
   \tau(\theta) -  \tau(\theta+\omega)=T_P(\theta) - T.
\]
This is the well-known small divisors equation (discussed here in
Section~\ref{sec:coho}).  We adjust then $T$ to be the average of $T_P$, $T=
\langle T_P\rangle$, and solve for $\tau$.  Notice that $\tau$ is defined up
to a constant, the average that we take as $0$. This is natural since $K$ is
determined up to a time translation (see Remark~\ref{rm:underterminacy2}).

\section{Quadratically small averages}
\label{ap:quadratiically small averages}

In this section we will prove that, given a multiple torus
$\{K_i(\theta)\}_{i= 0}^{m-1}$, approximately invariant with errors
\[
	E_i(\theta)= \varphi_{T/m} (K_i(\theta)) - K_{i+1}(\theta+\tfrac\omega{m}),
\]
then the averages of $\{\eta_i^3(\theta)\}_{i= 0}^{m-1}$ given by
\[
	\eta_i^3(\theta) = 
	\begin{pmatrix} 
	-{\rm D}K_{i+1}(\theta+\tfrac\omega{m})^\top \Omega(K_{i+1}(\theta+\tfrac\omega{m})) E_i(\theta) 
	\\
	-X_H(K_{i+1}(\theta+\tfrac\omega{m})) \Omega(K_{i+1}(\theta+\tfrac\omega{m})) E_i(\theta)
	\end{pmatrix} =: 
	\begin{pmatrix}
	\eta_i^{31}(\theta)
	\\
	\eta_i^{32}(\theta)
	\end{pmatrix}	
\]
are quadratically small with respect to the errors (and their derivatives).
This is a crucial step in our algorithms, and also in a eventual proof of
their convergence using KAM methods. 

First, we start by proving that $\langle {\eta}^{32} \rangle$ is quadratically
small. In fact
\[
\begin{split}
\sum_{i= 0}^{m-1} \langle \eta_i^{32}(\theta) \rangle 
= &\sum_{i= 0}^{m-1} \langle {\rm D}H(K_{i+1}(\theta+\tfrac\omega{m})) E_i(\theta) \rangle\\
 = &\sum_{i= 0}^{m-1} \langle H(\varphi_{T/m}(K_i(\theta)) - H(K_{i+1}(\theta+\tfrac\omega{m})) \rangle \\
&  
-  \sum_{i= 0}^{m-1} \langle \int_{0}^1 (1-s)\ {\rm D}^2 H(K_{i+1}(\theta+\tfrac\omega{m})+s E_i(\theta)) [E_i(\theta),E_i(\theta)]\ ds \rangle \\ = 
& -  \sum_{i= 0}^{m-1} \langle \int_{0}^1 (1-s)\ {\rm D}^2 H(K_{i+1}(\theta+\tfrac\omega{m})+s E_i(\theta)) [E_i(\theta),E_i(\theta)]\ ds \rangle,
\end{split}
\]
which is quadratically small in the errors. 

Second, we prove that $\langle { \eta}^{31} \rangle$ is quadratically small.
We start using the exactness of the symplectic form:
\[
\begin{split}
 \sum_{i= 0}^{m-1} \langle \eta_i^{31}(\theta) \rangle 
= &
 -\sum_{i= 0}^{m-1} 
 \langle {\rm D}K_{i+1}(\theta+\tfrac\omega{m})^\top {\rm D}a(K_{i+1}(\theta+\tfrac\omega{m}))^\top E_i(\theta)\rangle \\
 & 
 +  \sum_{i= 0}^{m-1} 
 \langle  {\rm D}K_{i+1}(\theta+\tfrac\omega{m})^\top {\rm D}a(K_{i+1}(\theta+\tfrac\omega{m})) E_i(\theta) \rangle \\
= &
 \sum_{i= 0}^{m-1} 
\langle  {\rm D}E_i(\theta)^\top a(K_{i+1}(\theta+\tfrac\omega{m})) + 
{\rm D}K_{i+1}(\theta+\tfrac\omega{m})^\top (\Delta^1 a_i(\theta)-\Delta^2 a_i(\theta))\rangle,
\end{split}
\]
where we use 
\[
\begin{split}	
	0 = & \langle  {\rm D}\bigl( a(K_{i+1}(\theta+\tfrac\omega{m}))^\top E_i(\theta) \bigr)  \rangle \\
	   = & \langle E_i(\theta)^\top  {\rm D}(a(K_{i+1}(\theta+\tfrac\omega{m}))) + 
	                     a(K_{i+1}(\theta+\tfrac\omega{m}))^\top {\rm D}E_i(\theta)   \rangle,
\end{split}	
\]
and the definitions
\[
\begin{split}
	\Delta^1 a_i(\theta)= &
	a(\varphi_{T/m}(K_i(\theta))) - a(K_{i+1}(\theta+\tfrac\omega{m})) 
	\\
	=&
	\int_{0}^1 {\rm D} a(K_{i+1}(\theta+\tfrac\omega{m})+s E_i(\theta)) E_i(\theta)\ ds
\end{split}
\]
and
\[
\begin{split}
	\Delta^2 a_i(\theta)= &
	a(\varphi_{T/m}(K_i(\theta))) - a(K_{i+1}(\theta+\tfrac\omega{m})) -
	{\rm D}a(K_{i+1}(\theta+\tfrac\omega{m})) E_i(\theta) 
	\\
	=&
	\int_{0}^1 (1-s)\ {\rm D}^2 a(K_{i+1}(\theta+\tfrac\omega{m})+s E_i(\theta)) [E_i(\theta),E_i(\theta)]\ ds.
\end{split}
\]
Hence, using that
\[
{\rm D}\varphi_{T/m}(K_i(\theta)) {\rm D}K_i(\theta) - {\rm D}K_{i+1}(\theta+\tfrac\omega{m}) = {\rm D}E_i(\theta)
\]
and that 
\[
\sum_{i= 0}^{m-1} \langle {\rm D}K_{i+1}( \theta+\tfrac\omega{m})^\top a(K_{i+1}(\theta+\tfrac\omega{m}))\rangle
= \sum_{i= 0}^{m-1} \langle {\rm D}K_{i}(\theta)^\top a(K_{i}(\theta))\rangle
\]
we have:
\[
\begin{split}
 \sum_{i= 0}^{m-1} \langle \eta_i^{31}(\theta) \rangle 
= &
\phantom{+}
\sum_{i= 0}^{m-1}    
\langle {\rm D}K_i(\theta)^\top 
\bigl( {\rm D}\varphi_{T/m}(K_i(\theta))^\top a(\varphi_{T/m}(K_i(\theta))) - a(K_{i}(\theta))\bigr) \rangle\\
& - \sum_{i= 0}^{m-1} 
\langle {\rm D}E_i(\theta)^\top \Delta^1 a_i(\theta) + {\rm D}K_{i+1}(\theta+\tfrac\omega{m})^\top \Delta^2 a_i(\theta)
\rangle 
\\
=&  \phantom{+}
\sum_{i= 0}^{m-1} \langle  {\rm D}\bigl(p_{T/m}(K_i(\theta))\bigr)^\top  \rangle 
\\
& - \sum_{i= 0}^{m-1} 
\langle {\rm D}E_i(\theta)^\top \Delta^1 a_i(\theta) + {\rm D}K_{i+1}(\theta+\tfrac\omega{m})^\top \Delta^2 a_i(\theta)\rangle
\\
=& 
- \sum_{i= 0}^{m-1} 
\langle {\rm D}E_i(\theta)^\top \Delta^1 a_i(\theta) + {\rm D}K_{i+1}(\theta+\tfrac\omega{m})^\top \Delta^2 a_i(\theta)\rangle,
\end{split}
\]
which is quadratically small. We have used the exactness of the Hamiltonian
flow, being $p_t$ is the primitive function of $\varphi_t$.


\bibliographystyle{plain}
\bibliography{references} 

\end{document}